\documentclass[11pt,reqno]{amsart}
\usepackage[T1]{fontenc}
\usepackage{amsmath}
\usepackage{amsthm}
\usepackage{amsfonts}
\usepackage{amssymb}
\usepackage[colorlinks=true]{hyperref}

\newtheorem{theorem}{Theorem}[section]
\newtheorem{lemma}[theorem]{Lemma}
\newtheorem{proposition}[theorem]{Proposition}

\theoremstyle{remark}

\numberwithin{equation}{section}

\catcode`@=11
\def\section{\@startsection{section}{1}%
  \z@{1.5\linespacing\@plus\linespacing}{.5\linespacing}%
  {\normalfont\bfseries\large\centering}}
\catcode`@=12

\newcommand{\be}{\begin{equation}}
\newcommand{\ee}{\end{equation}}
\newcommand{\R}{\mathbb{R}}
\newcommand{\Z}{\mathbb{Z}}
\newcommand{\N}{\mathbb{N}}
\newcommand{\T}{\mathbb{T}}
\newcommand{\Rp}{\mathbb{R}^2_+}
\newcommand{\pa}{\partial}
\newcommand{\eps}{\varepsilon}
\newcommand{\1}{\mathbf{1}}
\newcommand{\cL}{\mathcal{L}}
\newcommand{\cH}{\mathcal{H}}
\newcommand{\cD}{\mathcal{D}}
\def\fref{\eqref}

\DeclareMathOperator{\supp}{supp}
\DeclareMathOperator{\curl}{curl}
\DeclareMathOperator{\dv}{div}
\DeclareMathOperator{\dist}{dist}
\DeclareMathOperator{\sgn}{sgn}

\begin{document}

\title[Vorticity gradient growth in the plane]{Double exponential growth for the Euler equation on the whole plane}

\author[X. Chen]{Xi Chen}
\address[X.Chen]{Department of Mathematics and Computer Science, University of Basel, Spiegelgasse 1, 4051 Basel, Switzerland}
\email{xi01.chen@unibas.ch}

\author[Z. Zhang]{Zhifei Zhang}
\address[Z. Zhang]{School of Mathematical Sciences, Peking University, Beijing 100871, China}
\email{zfzhang@math.pku.edu.cn}

\subjclass[2020]{35Q31, 76B03, 76B47}

\keywords{Euler equation, growth of the vorticity gradient, small scale creation, Lamb dipole, orbital stability, Bahouri--Chemin solution, Serfati solutions}

\begin{abstract}
We consider the incompressible Euler equation on the whole plane $\R^2$, without boundary and without forcing, and we construct two smooth solutions whose vorticity gradient grows fast for all large times. The first solution has finite energy, integrable vorticity and initial data in every Sobolev space. Given any function $\Gamma(t)$ with $\Gamma(t)/t\to 0$, the logarithm of its vorticity gradient is not bounded by any multiple of $\Gamma(t)$, and, at every sufficiently large time $t$, it is larger than any multiple of $t/\log_*t$, where $\log_*$ is the iterated logarithm. The second solution has bounded vorticity, bounded vorticity derivatives on every finite time interval, and a velocity given by the absolutely convergent Biot--Savart integral which grows at most logarithmically in space. Its vorticity gradient grows double exponentially in time for all large times. Its initial vorticity is neither integrable nor square integrable, and its initial velocity and strain are unbounded. Both constructions place infinitely many localized amplification mechanisms far apart from each other in one initial datum. The first one is based on the hyperbolic stagnation point behind a smooth traveling dipole, controlled for long times by an energy--Casimir argument. The second one is based on the hyperbolic point of the Bahouri--Chemin cross on the torus, transported to the plane by spatial dilations and by the uniqueness theory of Serfati solutions.
\end{abstract}

\maketitle
\tableofcontents

\section{Introduction}

\subsection{The two-dimensional Euler equation}

In this paper we consider the incompressible Euler equation on the whole plane, written in vorticity form:
\be
\label{euler}
\left\{\begin{array}{ll}
\pa_t\omega+u\cdot\nabla\omega=0,\\
u=K*\omega,\ \ K(x)=\dfrac{1}{2\pi}\dfrac{x^\perp}{|x|^2},\\
\omega_{|t=0}=\omega_0,
\end{array}\right.
\qquad (t,x)\in\R_+\times\R^2,\ \ \omega(t,x)\in\R ,
\ee
where $x^\perp=(-x_2,x_1)$. The scalar $\omega$ is the vorticity and the vector field $u$ is the velocity. The kernel $K$ is the Biot--Savart kernel, namely $K=\nabla^\perp\left(\frac{1}{2\pi}\log|x|\right)$ with $\nabla^\perp=(-\pa_2,\pa_1)$, so that the velocity is divergence free and its curl is the vorticity:
$$\dv u=0,\ \ \curl u:=\pa_1u_2-\pa_2u_1=\omega.$$
Taking the curl of the velocity formulation
\be
\label{eulervelocity}
\pa_tu+u\cdot\nabla u+\nabla p=0,\ \ \dv u=0,
\ee
where $p$ is the pressure, yields the first equation of \fref{euler}. The vorticity is transported by the flow, and hence all its $L^p$ norms are conserved:
$$\|\omega(t)\|_{L^p}=\|\omega_0\|_{L^p},\ \ 1\le p\le \infty.$$
The global existence of smooth solutions goes back to Wolibner and H\"older \cite{Wolibner1933,Holder1933}; see Kato \cite{Kato1967} for a systematic treatment, Yudovich \cite{Yudovich1963} for the theory of bounded vorticities, and the monographs \cite{MajdaBertozzi2002,MarchioroPulvirenti1994}.

\medskip

The conservation of $\|\omega(t)\|_{L^\infty}$ does not prevent the creation of small scales. Indeed, differentiating the transport equation, the vorticity gradient solves
\be
\label{gradienteq}
\pa_t\nabla\omega+u\cdot\nabla(\nabla\omega)=-(\nabla u)^T\nabla\omega,
\ee
so that $\nabla\omega$ is amplified exactly where the flow stretches fluid elements. Integrating \fref{gradienteq} along the trajectories yields
\be
\label{gradientgronwall}
\|\nabla\omega(t)\|_{L^\infty}\le \|\nabla\omega_0\|_{L^\infty}\exp\left(\int_0^t\|\nabla u(s)\|_{L^\infty}ds\right).
\ee
The velocity gradient is controlled by the vorticity up to a logarithm: for integrable and Lipschitz vorticities,
\be
\label{logestimateintro}
\|\nabla u\|_{L^\infty}\lesssim \|\omega\|_{L^1}+\|\omega\|_{L^\infty}\left(1+\log\left(1+\frac{\|\nabla\omega\|_{L^\infty}}{\|\omega\|_{L^\infty}}\right)\right),
\ee
The singular-integral decomposition in the proof of Lemma \ref{lemmabiotsavart} also gives \fref{logestimateintro}: apply it directly to $\omega$, subtract $\omega(x)$ on $|z|<r$, and choose $r=M/(M+L)$, where $M=\|\omega\|_{L^\infty}>0$ and $L=\|\nabla\omega\|_{L^\infty}$. The region $|z|>1$ contributes at most $C\|\omega\|_{L^1}$ because $|\nabla K(z)|\lesssim |z|^{-2}$; the remaining terms are bounded by $CM(1+\log(1+L/M))$. The case $M=0$ is immediate. Inserting \fref{logestimateintro} into \fref{gradientgronwall} and using the conservation of $\|\omega\|_{L^1\cap L^\infty}$, one obtains the classical double exponential upper bound
\be
\label{doubleexpupper}
\log\left(e+\|\nabla\omega(t)\|_{L^\infty}\right)\le \log\left(e+\|\nabla\omega_0\|_{L^\infty}\right)e^{Ct},
\ee
where $C$ depends only on $\|\omega_0\|_{L^1\cap L^\infty}$. The corresponding argument for a smooth bounded domain is given in \cite[Theorem 2.1]{KiselevSverak2014}; on $\R^2$, the same Gronwall argument applies after using the whole-plane estimate \fref{logestimateintro}. The question we address in this paper is which lower bounds can actually be achieved by one smooth solution in a domain without boundary.

\subsection{Small scale creation}

The question of the sharpness of \fref{doubleexpupper} is classical and is closely related to the inherent instability of two-dimensional inviscid flows, see Yudovich \cite{Yudovich1974,Yudovich2000}, Nadirashvili \cite{Nadirashvili1991}, Koch \cite{Koch2002}, Morgulis, Shnirelman and Yudovich \cite{MSY2008}, and the surveys \cite{Constantin2007,Kiselev2018,DrivasElgindi2023}. Bahouri and Chemin \cite{BahouriChemin1994} exhibited a singular stationary solution on the torus, the {\it cross}, whose velocity is only log-Lipschitz at a hyperbolic stagnation point, and whose flow map loses H\"older regularity in time. This example shows that the logarithm in \fref{logestimateintro} is genuine, and it is the origin of all the double exponential mechanisms.

\medskip

\noindent{\it Domains with boundary}. A smooth solution whose vorticity gradient grows double exponentially for all times was constructed by Kiselev and \v Sver\'ak in the disk \cite{KiselevSverak2014}: the boundary allows one to place a hyperbolic point of the flow on a fluid particle trajectory which stays at the boundary, and to keep the vorticity near this point close to the Bahouri--Chemin cross for all times. This mechanism was extended to symmetric smooth domains by Xu \cite{Xu2016}, and Zlato\v s obtained double exponential growth with the maximal rate on the half-plane \cite{Zlatos2026}. For domains with corners, see \cite{ItohMiuraYoneda2016,KiselevZlatos2015}, and for vortex patches, see \cite{Denisov2015patch,KiselevLi2019}.

\medskip

\noindent{\it The torus}. In the absence of a boundary, Denisov constructed smooth solutions on the torus with superlinear growth of the vorticity gradient \cite{Denisov2009}, and showed that double exponential growth can occur on an arbitrarily long, but finite, time interval \cite{Denisov2015}; in the latter construction the initial datum depends on the length of the time interval. Zlato\v s proved exponential growth of the vorticity gradient for $C^{1,\alpha}$ vorticities on the torus which are smooth away from one point, and exponential growth of the vorticity Hessian for smooth data \cite{Zlatos2015}, and exponential growth of the gradient for $C^1$ data was obtained by Lei and Shi \cite{LeiShi2016}. Hoang and Radosz \cite{HoangRadosz2017} showed that a hyperbolic stagnation point does not produce local double exponential growth under additional assumptions on the flow.

\medskip

\noindent{\it The whole plane}. The first all-time growth result for smooth compactly supported vorticities in $\R^2$ is due to Choi and Jeong \cite{ChoiJeong2022}, who proved linear growth of the vorticity gradient for perturbations of the Lamb dipole, using its orbital stability \cite{AbeChoi2022}. Jeong, Yao and Zhou \cite{JeongYaoZhou2025} recently obtained superlinear growth for smooth compactly supported data. In both works, a stable coherent structure organizes the particle trajectories for all times. Growth mechanisms based on the twisting of particle trajectories near steady states were developed by Drivas, Elgindi and Jeong \cite{DrivasElgindiJeong2024}, and small scale creation for generic data was recently studied by Alazard and Said \cite{AlazardSaid2026}.

\medskip

\noindent{\it Stability of traveling dipoles}. The use of conserved quantities and rearrangements in the stability theory of vortices goes back to Arnold \cite{Arnold1965,Arnold1966}, Holm, Marsden, Ratiu and Weinstein \cite{HMRW1985}, Wan and Pulvirenti \cite{WanPulvirenti1985}, and Burton \cite{Burton1988}; see also the modulation theory of Weinstein \cite{Weinstein1986} for dispersive equations. Variational constructions and the orbital stability of vortex pairs were obtained by Burton, Nussenzveig Lopes and Lopes Filho \cite{BurtonLopes2013} and Burton \cite{Burton2021}, and Abe and Choi proved the orbital stability of the Lamb dipole \cite{AbeChoi2022}. The finite mass assumption and the sign condition were relaxed in \cite{AbeChoiJeong2025,AbeChoiJeongQin2026}, the stability of several Lamb dipoles was studied in \cite{AbeJeongYao2025}, quantitative estimates were obtained in \cite{LiSongZhou2026}, and a variational family of traveling dipoles was constructed in \cite{AbeJeongQinZhan2026}; see also \cite{ChoiLim2022} for radial vorticities. In the present paper, we need a smooth profile with an explicit orbital stability estimate, which we prove directly for the profile we construct, see Section \ref{sectiondipole}.

\medskip

\noindent{\it Remote placement of solutions}. Since the Biot--Savart kernel decays, two solutions placed far apart from each other interact weakly on a fixed time interval. This observation is at the heart of several ill-posedness results, where infinitely many localized solutions are placed in one initial datum, see Bourgain and Li \cite{BourgainLi2015,BourgainLi2015C1}, Misio\l ek and Yoneda \cite{MisiolekYoneda2016}, Elgindi and Jeong \cite{ElgindiJeong2017}, Elgindi and Masmoudi \cite{ElgindiMasmoudi2020}, Jeong \cite{Jeong2021}, and C\'ordoba, Mart\'inez-Zoroa and O\.za\'nski \cite{CordobaMartinezOzanski2024}. In the present paper, the same principle is used in the opposite direction: the solution remains smooth on every finite time interval, and remote placement serves to store infinitely many growth mechanisms, each of which acts on its own time interval.

\medskip

\noindent{\it The question of Kiselev and \v Sver\'ak}. At the end of \cite{KiselevSverak2014}, Kiselev and \v Sver\'ak state that the question of whether double exponential growth of the vorticity gradient can happen in the bulk of the fluid, in the periodic or in the full plane case, remains open, see \cite[Section 5, p. 1218]{KiselevSverak2014}. The aim of this paper is to give two answers to this question in the full plane, in two different function classes.

\subsection{Statement of the results}

We denote by $\log_*$ the iterated logarithm: for $t>1$,
$$\log_*t:=\min\left\{n\in\N\ :\ \log^{\circ n}t\le 1\right\},\ \ \log^{\circ 0}t=t,\ \ \log^{\circ(n+1)}t=\log\left(\log^{\circ n}t\right),$$
where the logarithms are iterated as long as their argument is larger than one. Hence $\log_* t$ is the number of logarithms needed to bring $t$ below one; it tends to infinity more slowly than any iterate $\log^{\circ n}t$. Our first result concerns finite energy solutions.

\begin{theorem}[Finite energy solutions with fast growth]
\label{thmmain}
Let $\Gamma:[1,+\infty)\to [1,+\infty)$ be a function such that $\Gamma(t)/t\to 0$ as $t\to +\infty$. There exists an initial vorticity
$$\omega_0\in L^1(\R^2)\cap\bigcap_{m\ge 0}H^m(\R^2)\ \ \mbox{with}\ \ u_0=K*\omega_0\in\bigcap_{m\ge 0}H^m(\R^2),$$
which is odd in $x_2$ and nonnegative on $\{x_2>0\}$, such that the global smooth solution $\omega$ of \fref{euler} satisfies the following.
\begin{enumerate}
\item {\em Growth beyond the prescribed rate along a sequence}: there exists a sequence $t_n\to +\infty$ such that
\be
\label{growthsequence}
\lim_{n\to +\infty}\frac{\log\left(2+\|\nabla\omega(t_n)\|_{L^\infty}\right)}{\Gamma(t_n)}=+\infty.
\ee
In particular, $\limsup_{t\to +\infty}\|\nabla\omega(t)\|_{L^\infty}e^{-c\Gamma(t)}=+\infty$ for every $c>0$.
\item {\em Growth at every large time}:
\be
\label{growthalltime}
\lim_{t\to +\infty}\frac{\log_*t}{t}\log\left(2+\|\nabla\omega(t)\|_{L^\infty}\right)=+\infty.
\ee
In particular, $\lim_{t\to +\infty}\|\nabla\omega(t)\|_{L^\infty}e^{-ct/\log_* t}=+\infty$ for every $c>0$.
\end{enumerate}
\end{theorem}

The initial vorticity of Theorem \ref{thmmain} has zero mean, since it is odd, and its support is not compact. The initial datum depends on $\Gamma$, but the sequence $(t_n)$ in \fref{growthsequence} does not depend on the constant $c>0$. Since $t/\log_* t$ is larger than $t^\beta$ for every $0<\beta<1$ and $t$ large, \fref{growthalltime} implies in particular that the vorticity gradient grows faster than every stretched exponential $\exp(ct^\beta)$, $0<\beta<1$, $c>0$, at every large time.

\medskip

Our second result is a double exponential growth for all large times, in a class of bounded vorticities which do not decay at infinity. We denote by $C^m_b(\R^2)$ the space of functions whose derivatives up to order $m$ are continuous and bounded on $\R^2$.

\begin{theorem}[Double exponential growth in the plane]
\label{thmdoubleexp}
There exist a solution $(\omega,u,p)\in C^\infty([0,+\infty)\times \R^2)$ of \fref{euler} and \fref{eulervelocity}, and constants $c_0>0$ and $t_0\ge 0$, such that the following holds.
\begin{enumerate}
\item {\em Biot--Savart law}: for every $(t,x)$, the function $y\mapsto K(x-y)\omega(t,y)$ is integrable on $\R^2$ and
$$u(t,x)=\int_{\R^2}K(x-y)\omega(t,y)dy.$$
Moreover, for every $S>0$ there exists $C_S$ such that
\be
\label{loggrowthvelocity}
|u(t,x)|\le C_S\log(2+|x|),\ \ 0\le t\le S,\ \ x\in\R^2.
\ee
\item {\em Regularity on finite time intervals}: $\|\omega(t)\|_{L^\infty}\le 1$ for all $t\ge 0$, and
$$\sup_{0\le t\le S}\|\omega(t)\|_{C^m_b}<+\infty\ \ \mbox{for every}\ \ S>0,\ m\ge 0.$$
\item {\em Double exponential growth}:
\be
\label{doubleexpgrowth}
\|\nabla\omega(t)\|_{L^\infty}\ge \exp\left(c_0e^{\frac{t}{\pi}}\right),\ \ t\ge t_0.
\ee
\item {\em Initial data}: $\omega_0\notin L^1(\R^2)\cup L^2(\R^2)$, $u_0\notin L^2(\R^2)\cup L^\infty(\R^2)$ and $\nabla u_0\notin L^\infty(\R^2)$.
\end{enumerate}
\end{theorem}

\subsection{Comments on the results}

\noindent 1. {\it Time scales and summability}. Both solutions are obtained by placing infinitely many localized solutions, which we call {\it components}, far apart from each other in one initial datum. Each component produces growth of the vorticity gradient on its own time interval. The two theorems correspond to two different ways of making infinitely many components compatible with the regularity of the datum. The Euler equation is invariant under the two scalings
\be
\label{scalings}
\omega(t,x)\mapsto a\,\omega(at,x),\qquad \omega(t,x)\mapsto \omega\left(t,\frac{x}{r}\right),\qquad a,r>0.
\ee
In Theorem \ref{thmmain}, the components must have summable Sobolev norms, and hence small amplitudes $a\ll 1$. By the first scaling in \fref{scalings}, a component of amplitude $a$ evolves on the slow time scale $t\sim a^{-1}$: an exponential growth $e^{\lambda t}$ of the profile becomes $e^{\lambda at}$ in physical time. This is why Theorem \ref{thmmain} produces growth rates $\exp(o(t))$ and not $\exp(ct)$. In Theorem \ref{thmdoubleexp}, the components keep an amplitude of order one, and are spread out by the second scaling in \fref{scalings} with $r\gg 1$, which does not change the time variable. The price to pay is that the $L^1$ and $L^2$ norms of the vorticity, the energy and the velocity of the components become large, which is the origin of the last item of Theorem \ref{thmdoubleexp}.

\medskip

\noindent 2. {\it The function class of Theorem \ref{thmdoubleexp}}. Bounded vorticities which do not decay at infinity have been studied by Serfati \cite{Serfati1995}, Ambrose, Kelliher, Lopes Filho and Nussenzveig Lopes \cite{AKLNL2015}, Taniuchi \cite{Taniuchi2004}, Taniuchi, Tashiro and Yoneda \cite{TTY2010}, Kelliher \cite{Kelliher2015}, Cozzi \cite{Cozzi2015}, Cozzi and Kelliher \cite{CozziKelliher2016,CozziKelliher2019}, and Cobb and Koch \cite{CobbKoch2026}. When the vorticity is not integrable, the velocity is not determined by the vorticity alone, and a normalization at infinity must be specified. In Theorem \ref{thmdoubleexp} the velocity is given at every time by the absolutely convergent Biot--Savart integral, which is the most classical normalization. The velocity grows at most logarithmically at infinity, as in the classes considered in \cite{Cozzi2015,CozziKelliher2019,CobbKoch2026}. We do not discuss the uniqueness of the solution in classes of unbounded velocities.

\medskip

\noindent 3. {\it On the question of Kiselev and \v Sver\'ak}. Theorem \ref{thmdoubleexp} shows that double exponential growth of the vorticity gradient can be sustained for all large times by one smooth solution of the unforced Euler equation without boundary. Since the initial energy, velocity and strain of this solution are infinite, Theorem \ref{thmdoubleexp} does not answer the question of \cite{KiselevSverak2014} for finite energy solutions, or for solutions with bounded velocity, which remains open. Theorem \ref{thmmain} is a step in this direction in the finite energy class, where the double exponential upper bound \fref{doubleexpupper} holds.

\medskip

\noindent 4. {\it On the rate in Theorem \ref{thmmain}}. The lower bound \fref{growthalltime} does not yield a growth $\exp(ct)$ for some $c>0$. On the time interval where the $n$-th component is active, our lower bound for $\log\|\nabla\omega(t)\|_{L^\infty}$ is proportional to $a_nt$, where $a_n\to 0$ is the amplitude of this component, see Comment 1. This is a limitation of the method, not an upper bound for the Euler equation.

\subsection{Strategy of the proof}

We outline in this subsection the main steps of the proofs of Theorems \ref{thmmain} and \ref{thmdoubleexp}.\\

\noindent{\bf Step 1:} {\it A smooth traveling dipole and its orbital stability}. The Lamb dipole $Q_L$ is an explicit vorticity, supported in the unit disk, odd in $x_2$, nonnegative in the upper half-plane, which travels with speed one along the horizontal axis. In the moving frame, the point $(-1,0)$ behind the dipole is a hyperbolic stagnation point of the velocity, with expansion rate $2$ along the symmetry axis. A fluid particle on the axis which remains close to $(-1,0)$ for a long time sees the vertical derivative of the vorticity grow exponentially. This geometric observation was used by Choi and Jeong \cite{ChoiJeong2022}. Two difficulties arise if one wants to use it with constants which do not depend on the length $T$ of the time interval: the Lamb dipole is only Lipschitz, and a perturbation of the dipole must remain under control on $[0,T]$. We therefore construct in Section \ref{sectiondipole} a smooth compactly supported traveling dipole $Q=F(\Psi)$, where $\Psi$ is the stream function in the moving frame, by smoothing the {\it function} $F$ rather than the vorticity, through a fixed point argument based on the nondegeneracy of the linearized operator of the Lamb dipole. The latter is computed explicitly with Bessel functions. We then prove the orbital stability of $Q$ with an energy--Casimir functional, following the ideas of Arnold \cite{Arnold1966} and Burton \cite{Burton1988}: under the constraint of fixed impulse, the vorticity remains $L^2$ close to a translate $Q(\cdot-z(t)e_1)$, with an explicit square root dependence on the size of the initial perturbation, and the speed $z'(t)$ remains close to one.\\

\noindent{\bf Step 2:} {\it Exponential stretching on a time interval of length $T$}. In the rest of this outline, $T\ge 1$ denotes the length of the time interval on which one component produces growth. We perturb $Q$ by $\delta P$, where $P$ is a fixed smooth function whose vertical derivative is equal to one on a neighborhood of $(-1,0)$ on the axis, and $\delta=\delta_0(1+T)^{-6}$ with a small constant $\delta_0>0$. The orbital stability gives an $L^2$ bound of order $\sqrt\delta$ on the perturbation $\eps(t)=\omega(t,\cdot+z(t)e_1)-Q$ of the dipole. Transport then gives an $L^\infty$ bound of order $T\sqrt{\delta|\log\delta|}\lesssim (1+T)^{-1}$ on $[0,T]$. The logarithmic estimate \fref{logestimateintro} for the velocity gradient of $\eps$ involves the logarithm of its Lipschitz norm, which is at most of order $T$, so that the product remains small uniformly in $T$. Hence the hyperbolic structure persists on $[0,T]$, and the vertical derivative of the vorticity on the axis grows like $\delta e^{t}$ for $0\le t\le T$. All the Sobolev norms of the initial data are bounded independently of $T$, see Proposition \ref{propstretching}.\\

\noindent{\bf Step 3:} {\it Remote assembly and choice of the time scales}. We rescale the solutions of Step 2 by \fref{scalings} with amplitudes $a_n\to 0$, sizes $r_n\to 0$, and lengths $T_n\to +\infty$. Two families of parameters are needed: one family produces the sequence of times in \fref{growthsequence}, and the other one produces intervals $[S_{n-1},S_n]$ which cover all large times, for \fref{growthalltime}. The second family requires the choice $T_n=\exp(\sqrt{T_{n-1}})$, which leads to the iterated logarithm. We then place the components at distances $d_n\to +\infty$ along the horizontal axis. The interaction between a component and the previous ones is arbitrarily small on a fixed time interval if $d_n$ is large enough (Lemma \ref{lemmadecoupling}), and the lower bounds for the vorticity gradient are preserved under this assembly (Lemma \ref{lemmaassembly}).\\

\noindent{\bf Step 4:} {\it Double exponential growth on a finite time interval on the torus}. For Theorem \ref{thmdoubleexp}, the local mechanism is the Bahouri--Chemin cross $\Theta(x)=\sgn(\sin x_1)\sgn(\sin x_2)$ on the torus $\T^2=(\R/2\pi\Z)^2$. Near the origin, its velocity is hyperbolic with a rate $\kappa\log(1/x_2)$, where $\kappa=\frac{2}{\pi}$. Hence the logarithm $\ell(t)=\log(1/x_2(t))$ of the inverse height of a fluid particle close to the vertical axis satisfies $\ell'\approx\kappa\ell$. A particle starting at height $e^{-\ell_0}$ therefore satisfies $\ell(t)\approx \ell_0e^{\kappa t}$, and horizontal segments near it are stretched by the factor $\exp(\int_0^t\kappa\ell)\approx\exp(\ell_0(e^{\kappa t}-1))$. We consider the smooth datum $\theta_\sigma+b_0$, where $\theta_\sigma$ is a regularization of $\Theta$ at scale $\sigma$, and $b_0$ is a small bump of height $A$ and radius $h$ placed at height $e^{-\ell_0}$ near the vertical axis. Two errors must be controlled on the time interval $[0,T]$: the effect of the smoothing, which is transported along the coordinate axes and remains far from the bump, and the strain created by the bump itself. The latter is controlled by the interpolation estimate $\|\nabla (K*b)\|_{L^\infty}\lesssim (\|\nabla b\|_{L^\infty}\|b\|_{L^2})^{\frac 12}$, by choosing $A$ as the inverse square root of the maximal stretching factor. An area preservation argument then shows that $\|\nabla\theta(t)\|_{L^\infty}\gtrsim \frac{A}{h}\exp(\ell_0(e^{\kappa t}-1))$, up to lower order factors, at every time $0\le t\le T$, see Proposition \ref{propperiodic}.\\

\noindent{\bf Step 5:} {\it Dilation and small derivatives at early times}. The dilated function $\theta(t,hx)$ is a solution of the Euler equation on a torus of period $2\pi/h$, with the same time variable. Its gradient is larger than $A\exp(\ell_0(e^{\kappa t}-1))$, which is a double exponential at times close to $T$ because $\log A\sim -\frac{\ell_0}{2}e^{\kappa T}$. Its derivatives at times $t\le k$ are bounded by $C_kA+C_{k,\sigma}h$, which is small. These two facts allow us to assign to the $n$-th component the time interval $[\tau_{n-1},\tau_n]$ with $\tau_n=\frac{n+n_0}{2}$, while keeping all the derivatives of the assembled solution bounded on finite time intervals.\\

\noindent{\bf Step 6:} {\it From the torus to the plane, and assembly}. A periodic solution has infinite energy. We cut off its stream function at a large scale $R$, and show that the corresponding solutions in the plane converge to the periodic one on compact sets, with global control of all the derivatives, see Proposition \ref{propcompactification}. The key tools are the Serfati identity and the uniqueness of Serfati solutions proved by Ambrose, Kelliher, Lopes Filho and Nussenzveig Lopes \cite{AKLNL2015}. We finally place the resulting compactly supported components at fast increasing distances, and check that the limit is a smooth solution whose velocity is given by the absolutely convergent Biot--Savart integral, see Proposition \ref{propassemblyinfinite}.\\

This paper is organized as follows. In Section \ref{sectionprelim}, we collect standard estimates on the Biot--Savart law, and we prove that solutions placed far apart from each other decouple. In Section \ref{sectiondipole}, we construct the smooth traveling dipole and prove its orbital stability. Section \ref{sectionfiniteenergy} is devoted to the proof of Theorem \ref{thmmain}. In Section \ref{sectionperiodic}, we construct the periodic solutions with double exponential growth on a finite time interval, and in Section \ref{sectionplane} we prove Theorem \ref{thmdoubleexp}.

\vskip 0.5cm

\noindent{\bf Notations}. We write $x=(x_1,x_2)\in\R^2$, $e_1=(1,0)$, and $\Rp=\{x_2>0\}$ for the upper half-plane. We denote by $B_R(x)$ the open disk of radius $R$ centered at $x$, $B_R=B_R(0)$ and $B=B_1$. We denote by $\1_A$ the indicator function of a set $A$. For a function $f$ on $\R^2$, $\nabla f$ is its gradient and $D^mf$ is the tensor of its partial derivatives of order $m$, measured with the Euclidean norm; for a vector field $v$, $\nabla v=(\pa_jv_i)_{1\le i,j\le 2}$ is measured with the operator norm. We use the $L^p$ norms $\|\cdot\|_{L^p}$ on the domain indicated in the context, and the inhomogeneous Sobolev spaces $H^s$ with norm
$$\|f\|_{H^s}=\|(1+|\xi|^2)^{\frac s2}\widehat f(\xi)\|_{L^2}.$$
For an integer $m\ge 0$ and $0<\gamma<1$, we use the H\"older norms
$$\|f\|_{C^{m,\gamma}}=\sum_{j=0}^m\|D^jf\|_{L^\infty}+\sup_{x\neq y}\frac{|D^mf(x)-D^mf(y)|}{|x-y|^\gamma}.$$
We write $C_0(\R^2)$ for the space of continuous functions vanishing at infinity, and $\log_+s=\max\{0,\log s\}$. For a continuous function $g$ of time, we use the upper Dini derivative
$$\frac{d^+}{dt}g(t)=\limsup_{\tau\downarrow 0}\frac{g(t+\tau)-g(t)}{\tau}.$$
We write $A\lesssim B$ if $A\le CB$ for a universal constant $C>0$, and $A\lesssim_S B$ if the constant $C$ depends on a parameter $S$. We write $A\sim B$ if $A\lesssim B$ and $B\lesssim A$. The letter $C$ denotes a constant which may change from line to line.

\section{Preliminaries}
\label{sectionprelim}

In this section we collect the tools which are used in both parts of the paper: classical estimates on the Biot--Savart law and on transport equations, and the decoupling of solutions which are placed far apart from each other.

\subsection{Biot--Savart law and transport}

Except for the velocity bound \fref{velocitybound}, which involves integrable vorticities, the estimates below are local in space: they only involve disks of radius at most $3$. They therefore apply to bounded vector fields on $\R^2$ which do not decay at infinity, and in particular to periodic vector fields, lifted to $\R^2$, with constants which do not depend on the period.

\begin{lemma}[Estimates on the Biot--Savart law]
\label{lemmabiotsavart}
Let $0<\gamma<1$ and let $m\ge 0$ be an integer.
\begin{enumerate}
\item {\em Velocity}: for $f\in L^1\cap L^\infty(\R^2)$ and $0<\rho\le 1$,
\be
\label{velocitybound}
\|K*f\|_{L^\infty}\lesssim \rho\left(\|f\|_{L^\infty}+\|f\|_{L^1}\right)+\left(\log\frac 2\rho\right)^{\frac 12}\|f\|_{L^2}.
\ee
\end{enumerate}
Let now $v$ be a bounded $C^1$ divergence-free vector field on $\R^2$, with bounded vorticity $f=\curl v$.
\begin{enumerate}
\setcounter{enumi}{1}
\item {\em Log-Lipschitz bound}: for $x,y\in\R^2$,
\be
\label{loglipschitz}
|v(x)-v(y)|\lesssim \left(\|v\|_{L^\infty}+\|f\|_{L^\infty}\right)|x-y|\left(1+\log_+\frac{1}{|x-y|}\right).
\ee
\item {\em Logarithmic estimate}: if $f$ is Lipschitz and not identically zero, then
\be
\label{logestimate}
\|\nabla v\|_{L^\infty}\lesssim \|v\|_{L^\infty}+\|f\|_{L^\infty}\left(1+\log_+\frac{\|\nabla f\|_{L^\infty}}{\|f\|_{L^\infty}}\right).
\ee
\item {\em Schauder estimates}: if $f\in C^{m,\gamma}$, then
\be
\label{schauder}
\|v\|_{C^{m+1,\gamma}}\lesssim_{m,\gamma}\|v\|_{L^\infty}+\|f\|_{C^{m,\gamma}},\ \ \ \|\nabla v\|_{C^{m,\gamma}}\lesssim_{m,\gamma}\|\nabla v\|_{L^\infty}+\|f\|_{C^{m,\gamma}}.
\ee
\item {\em Local version}: for every $x_0\in\R^2$, the estimates \fref{logestimate} and \fref{schauder} remain true if their left-hand sides are computed on $B_1(x_0)$, and if $\|v\|_{L^\infty}$ and the norms of $f$ are computed on $B_3(x_0)$.
\end{enumerate}
\end{lemma}

The two estimates in \fref{schauder} are used in different situations: the first one when the velocity is bounded uniformly with respect to the parameters of the problem, and the second one when only the velocity gradient is bounded uniformly. For instance, the velocity of a periodic solution dilated by a factor $h^{-1}$ is multiplied by $h^{-1}$, while its gradient is unchanged.

\begin{proof}
Proof of \fref{velocitybound}. Since $|K(z)|=(2\pi|z|)^{-1}$, we split the integral defining $K*f(x)$ into the regions $|z|<\rho$, $\rho\le |z|\le \rho^{-1}$ and $|z|>\rho^{-1}$, where $z=x-y$. The first region contributes at most $\rho\|f\|_{L^\infty}$, the last one at most $\rho\|f\|_{L^1}/(2\pi)$, and by the Cauchy--Schwarz inequality the middle one contributes at most
$$\frac{1}{2\pi}\left(\int_{\rho\le|z|\le \rho^{-1}}\frac{dz}{|z|^2}\right)^{\frac12}\|f\|_{L^2}=\frac{1}{2\pi}\left(4\pi\log\frac1\rho\right)^{\frac 12}\|f\|_{L^2}.$$

\noindent Proof of \fref{loglipschitz}, \fref{logestimate} and \fref{schauder}. We fix $x_0\in\R^2$ and a cutoff function $\chi_{x_0}=\chi(\cdot-x_0)$, where $\chi\in C^\infty_c(B_3)$ is equal to one on $B_2$. We split
$$v=w+(v-w),\ \ w:=K*(\chi_{x_0}f).$$
The log-Lipschitz bound \fref{loglipschitz} follows directly from \cite[Lemma B.3, pp.~4558--4559]{AKLNL2015}, applied to $v$. For the logarithmic estimate, write $g=\chi_{x_0}f$ and set $M=\|f\|_{L^\infty(B_3(x_0))}$ and $L=\|\nabla f\|_{L^\infty(B_3(x_0))}$. Then $\|g\|_{L^\infty}\lesssim M$, $\|\nabla g\|_{L^\infty}\lesssim M+L$, and $\|w\|_{L^\infty}\lesssim M$, since $g$ is supported in $B_3(x_0)$. For $x\in B_2(x_0)$, the distributional formula for $\nabla w$ consists of a fixed matrix times $g(x)$ and the principal value convolution of $g$ with $\nabla K$. The latter kernel has size $O(|z|^{-2})$ and zero average on circles. Subtracting $g(x)$ on $|z|<r\le 1$ gives a contribution at most $Cr\|\nabla g\|_{L^\infty}$; the part $r\le|z|\le 1$ is at most $CM\log(1/r)$; and the remaining part is at most $CM$, since $g(x-z)=0$ for $|z|\ge 5$. Consequently, for $M>0$, the choice $r=M/(M+L)$ yields
\begin{align*}
\|\nabla w\|_{L^\infty(B_2(x_0))}
&\lesssim M\left(1+\log\left(1+\frac LM\right)\right)\\
&\lesssim M\left(1+\log_+\frac LM\right).
\end{align*}
If $M=0$, then $w=0$, and the corresponding vorticity term is interpreted as zero.

The H\"older estimate for the second derivatives of the Newtonian potential \cite[Lemma 4.4]{GilbargTrudinger2001}, applied to $g$ and its derivatives, gives
$$\|\nabla w\|_{C^{m,\gamma}(B_2(x_0))}\lesssim_{m,\gamma}\|f\|_{C^{m,\gamma}(B_3(x_0))}.$$
On the other hand, $h=v-w$ is divergence free and curl free on $B_2(x_0)$, so its two components are harmonic there, and
$$\|h\|_{L^\infty(B_2(x_0))}\lesssim\|v\|_{L^\infty(B_3(x_0))}+M.$$
The interior estimates for harmonic functions \cite[Theorem 2.10]{GilbargTrudinger2001} bound all the derivatives of $h$ on $B_1(x_0)$ by this quantity. Together with the preceding bounds for $w$, this proves the local logarithmic estimate and the first local Schauder estimate. For the second Schauder estimate, we apply the same harmonic estimates to $\nabla h$, using
\begin{align*}
\|\nabla h\|_{L^\infty(B_2(x_0))}
&\le\|\nabla v\|_{L^\infty(B_2(x_0))}+\|\nabla w\|_{L^\infty(B_2(x_0))}\\
&\lesssim\|\nabla v\|_{L^\infty(B_2(x_0))}+\|f\|_{C^{0,\gamma}(B_3(x_0))}.
\end{align*}
Thus this estimate uses the velocity gradient, rather than the velocity, on its right-hand side. All these local bounds have constants independent of $x_0$. Taking the supremum over $x_0$ gives the global estimates: H\"older quotients at distances less than $1$ are controlled in a unit disk, while those at distances at least $1$ are bounded by twice the corresponding supremum norm. This also proves the local version stated in the lemma.
\end{proof}

We shall also use the following classical estimates for transport equations, see \cite[Proposition 3.10, pp.~129--130, Section 3.2, p.~132, and (3.20), p.~135]{BahouriCheminDanchin2011}.

\begin{lemma}[Transport estimates]
\label{lemmatransport}
Let $v$ be a $C^1$ divergence-free vector field on $[0,S]\times\R^2$ with $v$ and $\nabla v$ bounded, and let $g$ be a $C^1$ solution of
$$\pa_tg+v\cdot\nabla g=F,\ \ 0\le t\le S.$$
Let $X(t,\alpha)$ be the flow of $v$: $\pa_tX(t,\alpha)=v(t,X(t,\alpha))$, $X(0,\alpha)=\alpha$. Then for $0\le t\le S$:
\begin{enumerate}
\item {\em Lagrangian formula}: $g(t,X(t,\alpha))=g(0,\alpha)+\int_0^tF(s,X(s,\alpha))ds$, and $X(t,\cdot)$ is a measure preserving diffeomorphism of $\R^2$.
\item {\em Lebesgue norms}: for $p\in\{1,+\infty\}$,
$$\|g(t)\|_{L^p}\le\|g(0)\|_{L^p}+\int_0^t\|F(s)\|_{L^p}ds.$$
\item {\em H\"older norms}: for $0<\gamma<1$,
\be
\label{transportholder}
\|g(t)\|_{C^{0,\gamma}}\le \left(\|g(0)\|_{C^{0,\gamma}}+\int_0^t\|F(s)\|_{C^{0,\gamma}}ds\right)\exp\left(\gamma\int_0^t\|\nabla v(s)\|_{L^\infty}ds\right).
\ee
\end{enumerate}
\end{lemma}

The same estimates hold, with the same proof, for periodic vector fields, the norms being computed on a period cell.

\begin{proof}
The first item follows by differentiating $g(t,X(t,\alpha))$, and the Jacobian determinant of $X(t,\cdot)$ equals one because $\dv v=0$. The second item follows from the Lagrangian formula and this measure preservation. For the third item, let $\Phi(s;t,x)$ denote the position at time $s$ of the trajectory which is at $x$ at time $t$. Gronwall's lemma gives
$$\|D_x\Phi(s;t,\cdot)\|_{L^\infty}\le\exp\left(\int_s^t\|\nabla v(\tau)\|_{L^\infty}d\tau\right),\qquad 0\le s\le t.$$
Writing the Lagrangian formula at the terminal point gives
$$g(t,x)=g(0,\Phi(0;t,x))+\int_0^tF(s,\Phi(s;t,x))ds.$$
If $[\cdot]_{\gamma}$ denotes the H\"older seminorm, the composition estimate $[q\circ\Phi]_{\gamma}\le[q]_{\gamma}\|D\Phi\|_{L^\infty}^{\gamma}$ therefore yields
\begin{align*}
[g(t)]_{\gamma}
&\le[g(0)]_{\gamma}\exp\left(\gamma\int_0^t\|\nabla v(\tau)\|_{L^\infty}d\tau\right)\\
&\quad+\int_0^t[F(s)]_{\gamma}\exp\left(\gamma\int_s^t\|\nabla v(\tau)\|_{L^\infty}d\tau\right)ds.
\end{align*}
Adding the $L^\infty$ estimate and enlarging the time interval in each exponential to $[0,t]$ proves \fref{transportholder}.
\end{proof}

\subsection{Solutions placed far apart}

The Euler equation is nonlocal, but the Biot--Savart kernel decays. Hence two smooth finite energy solutions placed at a large distance $d$ from each other interact weakly on a fixed time interval. The following lemma quantifies this fact in Sobolev norms.

\begin{lemma}[Decoupling of distant solutions]
\label{lemmadecoupling}
Let $u_A$ and $u_B$ be global smooth solutions of \fref{eulervelocity} on $\R^2$ such that $u_A,u_B\in C([0,S];H^m(\R^2))$ for all $S>0$ and all integers $m\ge 0$. For $d>0$, let $u_d$ be the global smooth solution of \fref{eulervelocity} with initial velocity
$$u_d(0,x)=u_A(0,x)+u_B(0,x-de_1).$$
Then for every $S>0$ and every integer $s\ge 3$,
\be
\label{decoupling}
\lim_{d\to+\infty}\ \sup_{0\le t\le S}\|u_d(t)-u_A(t)-u_B(t,\cdot-de_1)\|_{H^s}=0.
\ee
\end{lemma}

\begin{proof}
We denote by $\mathbb P$ the Leray projection onto divergence-free vector fields, so that a smooth solution of \fref{eulervelocity} in $H^s$ satisfies $\pa_tu+\mathbb P(u\cdot\nabla u)=0$. The global existence of $u_d$ in $C([0,S];H^m)$ for all $S>0$ and all integers $m\ge 0$ follows from \cite[Theorems 7.1 and 7.25]{BahouriCheminDanchin2011}: for each $m\ge 3$, take $p=r=a=2$ and use $B^m_{2,2}=H^m$; the hypothesis $\nabla u_d(0)\in L^2$ holds, and the lower Sobolev orders follow by embedding. Uniqueness identifies the solutions obtained at different Sobolev orders.\\

\noindent{\bf Step 1:} The error of the sum. The approximate solution $v_d(t,x)=u_A(t,x)+u_B(t,x-de_1)$ satisfies
\begin{align*}
&\pa_tv_d+\mathbb P(v_d\cdot\nabla v_d)=R_d,\\
&R_d:=\mathbb P\left[(u_A\cdot\nabla)u_B(\cdot-de_1)+(u_B(\cdot-de_1)\cdot\nabla)u_A\right].
\end{align*}
We claim that $\sup_{0\le t\le S}\|R_d(t)\|_{H^s}\to 0$ as $d\to+\infty$. Indeed, the sets $\{u_A(t):0\le t\le S\}$ and $\{u_B(t):0\le t\le S\}$ are compact in $H^{s+1}$. Let $\chi_\rho=\chi(\cdot/\rho)$, where $\chi\in C^\infty_c(B_2)$ is equal to one on $B_1$. Since $\chi_\rho g\to g$ in $H^{s+1}$ for every $g\in H^{s+1}$ as $\rho\to +\infty$, and the multiplication by $\chi_\rho$ is bounded on $H^{s+1}$ uniformly in $\rho\ge 1$, this convergence is uniform on compact sets. Given $\eta>0$, we fix $\rho$ such that $\|(1-\chi_\rho)u_A(t)\|_{H^{s+1}}+\|(1-\chi_\rho)u_B(t)\|_{H^{s+1}}\le\eta$ for $0\le t\le S$. For $d>4\rho$, the functions $\chi_\rho u_A$ and $(\chi_\rho u_B)(\cdot-de_1)$ have disjoint supports, so that their products vanish. Using $\|f\cdot\nabla g\|_{H^s}\lesssim_s \|f\|_{H^{s+1}}\|g\|_{H^{s+1}}$ for $s\ge 2$, and the boundedness of $\mathbb P$ on $H^s$, we infer $\|R_d(t)\|_{H^s}\lesssim \eta$ for $0\le t\le S$, and the claim follows.\\

\noindent{\bf Step 2:} Energy estimate. The difference $q_d=u_d-v_d$ vanishes at $t=0$ and solves
$$\pa_tq_d+\mathbb P\left[v_d\cdot\nabla q_d+q_d\cdot\nabla v_d+q_d\cdot\nabla q_d\right]=-R_d.$$
Let $J=(1-\Delta)^{1/2}$. The commutator estimate \cite[(1.1)--(1.2), p.~891, and Remark 1.1(a), p.~892]{KatoPonce1988}, with $p=2$, gives
\begin{align*}
\|[J^s,w\cdot\nabla]q\|_{L^2}
\lesssim_s{}&\|\nabla w\|_{L^\infty}\|q\|_{H^s}
+\|w\|_{H^s}\|\nabla q\|_{L^\infty}.
\end{align*}
The Leray projection is self-adjoint and commutes with $J^s$, and $\mathbb Pq_d=q_d$. For each of the divergence-free fields $w=v_d,q_d$, the highest-order transport term cancels:
$$((w\cdot\nabla)J^sq_d,J^sq_d)_{L^2}=0.$$
Combining these facts with $H^s(\R^2)\subset W^{1,\infty}(\R^2)$ and $\|q_d\cdot\nabla v_d\|_{H^s}\lesssim_s\|q_d\|_{H^s}\|v_d\|_{H^{s+1}}$ yields, for $s\ge 3$:
$$\frac{d}{dt}\|q_d\|_{H^s}\lesssim_s \left(\|v_d\|_{H^{s+1}}+\|q_d\|_{H^s}\right)\|q_d\|_{H^s}+\|R_d\|_{H^s}.$$
The norm $\|v_d(t)\|_{H^{s+1}}$ is bounded on $[0,S]$ independently of $d$. As long as $\|q_d\|_{H^s}\le 1$, Gronwall's lemma gives
$$\sup_{0\le t\le S}\|q_d(t)\|_{H^s}\le C_S\sup_{0\le t\le S}\|R_d(t)\|_{H^s},$$
with $C_S$ independent of $d$. For $d$ large, the right-hand side is smaller than $\frac12$, and a continuity argument shows that this bound holds on the whole interval $[0,S]$. This proves \fref{decoupling}.
\end{proof}

We now place infinitely many solutions far apart from each other and show that lower bounds on the vorticity gradient of each of them, on given compact sets of space and time, survive the assembly. The point is that the distances are chosen one after the other: when the $n$-th solution is placed, the first $n-1$ ones are already fixed, and only finitely many conditions must be satisfied.

\begin{lemma}[Assembly of distant solutions]
\label{lemmaassembly}
Let $(V_n)_{n\ge 1}$ be global smooth solutions of \fref{eulervelocity} with $V_n\in C([0,S];H^m(\R^2))$ for all $S>0$ and $m\ge 0$, and let $\omega_n=\curl V_n$. Let $0<S_1\le S_2\le\dots$, let $s_n\ge 4$ be integers, let $\eps_n\in(0,1]$, and let $\Sigma_n\subset[0,S_n]\times\R^2$ be compact sets. For positive numbers $d_1,d_2,\dots$, let $u^{(N)}$ be the global smooth solution of \fref{eulervelocity} with initial velocity
\be
\label{defprefix}
u^{(N)}(0,x)=\sum_{n=1}^NV_n(0,x-d_ne_1),
\ee
let $\omega^{(N)}=\curl u^{(N)}$, $u^{(0)}=0$, and let
\be
\label{definteraction}
w_n:=u^{(n)}-u^{(n-1)}-V_n(\cdot,\cdot-d_ne_1)
\ee
be the interaction error created by the $n$-th solution. There exist thresholds $\overline d_n$, depending only on $(V_j,S_j,s_j,\eps_j,\Sigma_j)_{1\le j\le n}$ and on $d_1,\dots,d_{n-1}$, such that if $d_n\ge \overline d_n$ for all $n\ge 1$, then the following holds for all $n\ge 1$.
\begin{enumerate}
\item {\em Small interaction}:
\be
\label{smallinteraction}
\sup_{0\le t\le S_n}\|w_n(t)\|_{H^{s_n}}\le \eps_n.
\ee
\item {\em Preservation of the vorticity gradient}: for all $N\ge n$ and $(t,x)\in\Sigma_n$,
\be
\label{preservation}
\left|\nabla\omega^{(N)}(t,x+d_ne_1)-\nabla\omega_n(t,x)\right|\le C_0\sum_{j=n}^N\eps_j,
\ee
where $C_0$ is a universal constant.
\end{enumerate}
\end{lemma}

\begin{proof}
Let $C_{\rm S}$ be the constant of the Sobolev embedding $\|\nabla\curl f\|_{L^\infty}\le C_{\rm S}\|f\|_{H^4}$ for vector fields $f\in H^4(\R^2)$, and let $C_0=2+C_{\rm S}$. Assume that $d_1,\dots,d_{n-1}$ are fixed. We choose $\overline d_n$ such that, for every $d_n\ge\overline d_n$, the three following conditions hold:
\begin{align*}
&{\rm (a)}\ \ \sup_{0\le t\le S_n}\|w_n(t)\|_{H^{s_n}}\le\eps_n,\\
&{\rm (b)}\ \ |\nabla\omega^{(n-1)}(t,x+d_ne_1)|\le\eps_n\ \ \mbox{for}\ \ (t,x)\in\Sigma_n,\\
&{\rm (c)}\ \ |\nabla\omega_n(t,x+(d_j-d_n)e_1)|\le\eps_n\ \ \mbox{for}\ \ (t,x)\in\Sigma_j,\ \ 1\le j\le n-1.
\end{align*}
Condition (a) holds for $d_n$ large by Lemma \ref{lemmadecoupling} applied to $u_A=u^{(n-1)}$, $u_B=V_n$, $S=S_n$ and $s=s_n$. For (b) and (c), we observe that a solution $u\in C([0,S];H^4)$ satisfies $\nabla\curl u\in C([0,S];H^2)\subset C([0,S];C_0(\R^2))$; the set $\{\nabla\curl u(t):0\le t\le S\}$ is then compact in $C_0(\R^2)$, and hence $\sup_{0\le t\le S}|\nabla\curl u(t,x)|\to 0$ as $|x|\to+\infty$. Since $\Sigma_j$ is compact, its projection on $\R^2$ is bounded, and (b), (c) hold for $d_n$ large. The definition of $\overline d_n$ is complete, and (a) is \fref{smallinteraction}.

We now prove \fref{preservation}. By \fref{definteraction}, for $N\ge n$,
$$\omega^{(N)}=\omega^{(n-1)}+\omega_n(\cdot,\cdot-d_ne_1)+\sum_{j=n+1}^N\omega_j(\cdot,\cdot-d_je_1)+\sum_{j=n}^N\curl w_j.$$
Let $(t,x)\in\Sigma_n$, so that $t\le S_n\le S_j$ for $j\ge n$. We evaluate the gradient of this identity at $(t,x+d_ne_1)$. The first term is bounded by $\eps_n$ by (b); the terms of the first sum are bounded by $\eps_j$ by (c) applied at the step $j$; and by \fref{smallinteraction} and the choice of $C_{\rm S}$, $|\nabla\curl w_j(t,\cdot)|\le C_{\rm S}\|w_j(t)\|_{H^4}\le C_{\rm S}\eps_j$. This proves \fref{preservation}.
\end{proof}

\section{A smooth traveling dipole and its orbital stability}
\label{sectiondipole}

This section is devoted to the construction of a smooth compactly supported traveling wave $Q$ of \fref{euler}, close to the Lamb dipole, and to the proof of its orbital stability with explicit constants. These are the two properties of the reference profile which make the growth mechanism of Section \ref{sectionfiniteenergy} uniform with respect to the length of the time interval.

\subsection{The Lamb dipole and its linearized operator}

\noindent{\it Odd vorticities}. All the vorticities of this section and of Section \ref{sectionfiniteenergy} are odd in $x_2$. We identify such a function with its restriction to the upper half-plane $\Rp$. We denote by $(f,g)=\int_{\Rp}fg\,dx$ the scalar product of $L^2(\Rp)$, and by $G=(-\Delta_{\Rp})^{-1}$ the inverse of the Dirichlet Laplacian of $\Rp$:
$$Gf(x)=\int_{\Rp}G(x,y)f(y)dy,\ \ G(x,y)=\frac{1}{2\pi}\log\frac{|x-\overline y|}{|x-y|},\ \ \overline y=(y_1,-y_2).$$
The kernel $G(x,y)$ is positive on $\Rp\times\Rp$. If $f$ is extended oddly to $\R^2$, its stream function $\frac{1}{2\pi}\log|\cdot|*f$ is equal to $-Gf$ on $\Rp$, and hence its velocity is
\be
\label{velocityodd}
K*f=-\nabla^\perp Gf\ \ \mbox{on}\ \ \Rp.
\ee
We use the impulse and the energy
\be
\label{defimpulseenergy}
I(f)=\int_{\Rp}x_2f\,dx,\qquad E(f)=\frac12(Gf,f)=\frac 12\int_{\Rp}|K*f|^2=\frac14\int_{\R^2}|K*f|^2.
\ee
For a smooth solution of \fref{euler} whose vorticity is odd in $x_2$ and compactly supported, $I(\omega(t))$ and $E(\omega(t))$ are conserved. Indeed, on every finite time interval the vorticity supports lie in a common compact set, by transport and the boundedness of $K*\omega$. Oddness gives $\int_{\R^2}\omega=0$ and $u_2=0$ on $\{x_2=0\}$. Hence, setting $\phi=G\omega$ on $\Rp$, we have $-\Delta\phi=\omega$, $\phi=0$ on the axis, $u=-\nabla^\perp\phi$, and $\phi=O(|x|^{-1})$, $\nabla\phi=O(|x|^{-2})$ at infinity. Integration by parts therefore gives $E(\omega)=\frac12\int_{\Rp}|\nabla\phi|^2=\frac14\int_{\R^2}|u|^2$. For the impulse, the transport equation and the oddness of the Biot--Savart kernel give
\begin{align*}
2\frac{d}{dt}I(\omega)
&=\frac{d}{dt}\int_{\R^2}x_2\omega(x)dx=\int_{\R^2}u_2\omega\\
&=\iint_{\R^2\times\R^2}K_2(x-y)\omega(x)\omega(y)dxdy=0,
\end{align*}
since interchanging $x$ and $y$ changes the sign of the last integral. This double integral is absolutely convergent because $|K(z)|\lesssim |z|^{-1}$ and $\omega$ is bounded with compact support. The symmetry of $G$ and the transport equation similarly yield
\begin{align*}
\frac{d}{dt}E(\omega)
&=\int_{\Rp}\phi\,\pa_t\omega=-\int_{\Rp}\phi\,u\cdot\nabla\omega\\
&=\int_{\Rp}\omega\,u\cdot\nabla\phi=0.
\end{align*}
Here the boundary terms vanish because $u_2=0$ on the axis and $\omega$ is compactly supported; the last equality follows from $u=-\nabla^\perp\phi$. Smoothness and the common compact support justify differentiation of the Green-kernel quadratic form.

\medskip

\noindent{\it The Lamb dipole}. We denote by $J_\nu$ the Bessel function of the first kind of order $\nu$, and by $j_{\nu,\ell}$ its $\ell$-th positive zero. We let
$$k:=j_{1,1}\approx 3.8317,$$
and we use polar coordinates $x=(r\cos\vartheta,r\sin\vartheta)$. The Lamb dipole \cite{Lamb1906,Chaplygin1903,MeleshkoVanHeijst1994} is the vorticity
\be
\label{deflamb}
Q_L(x)=\left\{\begin{array}{ll}-\dfrac{2k}{J_0(k)}J_1(kr)\sin\vartheta, & r<1,\\[2mm] 0, & r\ge 1.\end{array}\right.
\ee
Since $J_0(k)<0$ and $J_1>0$ on $(0,k)$, $Q_L$ is positive on $B\cap\Rp$. It is odd in $x_2$, even in $x_1$, Lipschitz and supported in the closed unit disk. Its stream function in the frame moving with speed one,
$$\Psi_L:=GQ_L-x_2,$$
is explicitly given on $\Rp$ by
\be
\label{lambstream}
\Psi_L(x)=\left\{\begin{array}{ll}-\dfrac{2}{kJ_0(k)}J_1(kr)\sin\vartheta, & r<1,\\[2mm] \left(\dfrac1{r^2}-1\right)x_2, & r\ge 1.\end{array}\right.
\ee
Indeed, the right-hand side of \fref{lambstream} plus $x_2$ vanishes on $\{x_2=0\}$ and at infinity, it is $C^1$ across the circle $\{r=1\}$ because $J_1(k)=0$ and $kJ_1'(k)=kJ_0(k)$, and, by the Bessel equation $\Delta(J_1(kr)\sin\vartheta)=-k^2J_1(kr)\sin\vartheta$, minus its Laplacian is equal to $Q_L$. Comparing \fref{deflamb} and \fref{lambstream}, $\Psi_L>0$ on $B\cap\Rp$, $\Psi_L<0$ on $\Rp\setminus\overline{B}$, and
\be
\label{lambrelation}
Q_L=k^2(\Psi_L)_+\ \ \mbox{on}\ \ \Rp.
\ee
By \fref{velocityodd} and $\nabla^\perp x_2=-e_1$, the velocity $U_L=K*Q_L$ satisfies $U_L-e_1=-\nabla^\perp\Psi_L$, so that \fref{lambrelation} implies $(U_L-e_1)\cdot\nabla Q_L=0$: the Lamb dipole travels with speed one in the direction $e_1$. On the symmetry axis, \fref{lambstream} gives
\be
\label{lambaxis}
U_{L,1}(x_1,0)=\left\{\begin{array}{ll}x_1^{-2}, & |x_1|\ge 1,\\[1mm] 1-\dfrac{2J_1(k|x_1|)}{kJ_0(k)|x_1|}>1, & 0<|x_1|<1.\end{array}\right.
\ee
Hence, in the moving frame, the horizontal velocity $U_{L,1}-1$ on the axis vanishes at the rear point $(-1,0)$, is negative on its left and positive on its right, and $\pa_1U_{L,1}(-1,0)=2$. The point $(-1,0)$ is a hyperbolic stagnation point of the relative flow, and fluid particles on the axis are pushed away from it at the exponential rate $2$. This is the stretching mechanism used in Section \ref{sectionfiniteenergy}.

\medskip

\noindent{\it The linearized operator}. The stability of the dipole is governed by the quadratic form $\frac{1}{2k^2}\|f\|_{L^2}^2-E(f)$, which is the second variation of the energy--Casimir functional introduced in Subsection \ref{subsectionorbital}. Its restriction to functions supported in the disk $B_R$ is represented by the bounded self-adjoint operator on $L^2(\Rp)$
\be
\label{deflinearized}
\cL_R:=\frac{1}{2k^2}{\rm Id}-\frac 12\1_{B_R}G\1_{B_R},\ \ R\ge 1,\qquad \cL:=\cL_1.
\ee
By \fref{lambrelation}, $\1_BGQ_L=\1_B(\Psi_L+x_2)=\1_B(k^{-2}Q_L+x_2)$, and hence
\be
\label{lambnegative}
\cL Q_L=-\frac12x_2\1_B,\qquad (\cL Q_L,Q_L)=-\frac12I(Q_L)<0.
\ee
Thus $\cL$ has a negative direction, which is related to the impulse. The following lemma shows that it is the only one, and that the constraint of fixed impulse removes it.

\begin{lemma}[Energy of odd vorticities]
\label{lemmaenergy}
For $f\in L^2(\Rp)$ with $x_2f\in L^1(\Rp)$,
\be
\label{energyinterpolation}
E(f)\lesssim \|x_2f\|_{L^1(\Rp)}\|f\|_{L^2(\Rp)}.
\ee
\end{lemma}

\begin{proof}
It suffices to consider $f\in C^\infty_c(\Rp)$ and to conclude by density. The Fourier transform of the odd extension of $f$ is
$$\widehat f(\xi)=\frac{1}{2\pi}\int_{\R^2}e^{-ix\cdot\xi}f(x)dx=-\frac{i}{\pi}\int_{\Rp}e^{-ix_1\xi_1}\sin(x_2\xi_2)f(x)dx,$$
and $|\sin(x_2\xi_2)|\le x_2|\xi|$ yields $|\widehat f(\xi)|\le \pi^{-1}|\xi|\|x_2f\|_{L^1(\Rp)}$. By \fref{defimpulseenergy} and Plancherel, for every $\rho>0$,
\begin{align*}
E(f)&=\frac14\int_{\R^2}\frac{|\widehat f(\xi)|^2}{|\xi|^2}d\xi\le \frac{1}{4\pi^2}\int_{|\xi|\le\rho}\|x_2f\|_{L^1(\Rp)}^2d\xi+\frac{1}{4\rho^2}\int_{\R^2}|\widehat f(\xi)|^2d\xi\\
&\lesssim \rho^2\|x_2f\|^2_{L^1(\Rp)}+\rho^{-2}\|f\|^2_{L^2(\Rp)}.
\end{align*}
The choice $\rho^2=\|f\|_{L^2(\Rp)}/\|x_2f\|_{L^1(\Rp)}$ yields \fref{energyinterpolation}.
\end{proof}

\begin{lemma}[Spectral structure of the Lamb dipole]
\label{lemmaspectral}
The operators $\cL_R$ defined in \fref{deflinearized} have the following properties.
\begin{enumerate}
\item {\em Spectrum}: the operator $\cL$ has exactly one negative eigenvalue, which is simple; its kernel is spanned by $\pa_1Q_L$; and there exists $\mu_0>0$ such that the rest of its spectrum is contained in the half-line $[\mu_0,+\infty)$. In particular, every solution $f\in L^2(\Rp)$ of the equation $f=k^2\1_BGf$ is a multiple of $\pa_1Q_L$.
\item {\em Coercivity under constraints}: let $\varphi$ be a fixed smooth function with compact support in $B\cap\Rp$, odd in $x_1$, such that $(\pa_1Q_L,\varphi)\neq 0$. There exist $R_0\in(1,2]$ and $\mu_1>0$ such that for all $R\in[1,R_0]$ and all $f\in L^2(\Rp)$,
\be
\label{coercivityconstraint}
(f,x_2\1_{B_R})=(f,\varphi)=0\ \ \Longrightarrow\ \ (\cL_Rf,f)\ge\mu_1\|f\|^2_{L^2(\Rp)}.
\ee
\end{enumerate}
\end{lemma}

The first constraint in \fref{coercivityconstraint} is a localized version of the conservation of the impulse, and the second one fixes the translation parameter of the dipole. Such a function $\varphi$ exists because $\pa_1Q_L$ is odd in $x_1$ and does not vanish identically on $B\cap\Rp$.

\begin{proof}
\noindent{\bf Step 1:} Eigenfunctions. The operator $\mathcal T=\1_BG\1_B$ is self-adjoint, nonnegative, and compact on $L^2(\Rp)$: since $|x-\overline y|\le 2$ for $x,y\in B$, its kernel is bounded by $\frac{1}{2\pi}\log\frac{2}{|x-y|}$, which is square integrable on $B\times B$. Let $Tg=\mu g$ with $\mu>0$ and $g\neq 0$. Then $g$ is supported in $B$, and $\phi=Gg$ satisfies
\begin{align*}
&-\Delta\phi=\mu^{-1}\phi\ \ \mbox{in}\ \ B\cap\Rp,\qquad \Delta\phi=0\ \ \mbox{in}\ \ \Rp\setminus \overline B,\\
&\phi=0\ \ \mbox{on}\ \ \{x_2=0\},\qquad \phi(x)\to 0\ \ \mbox{as}\ \ |x|\to+\infty,
\end{align*}
and $\phi$ is $C^1$ across the circle $\{r=1\}$. We expand $\phi=\sum_{m\ge 1}\phi_m(r)\sin(m\vartheta)$ for $\vartheta\in(0,\pi)$. Regularity at the origin and decay at infinity impose $\phi_m(r)=\alpha_mJ_m(\mu^{-\frac12}r)$ for $r<1$ and $\phi_m(r)=\beta_mr^{-m}$ for $r>1$. The continuity of $\phi_m$ and $\phi_m'$ at $r=1$, together with the identity $sJ_m'(s)+mJ_m(s)=sJ_{m-1}(s)$, gives $\alpha_mJ_{m-1}(\mu^{-\frac12})=0$. Hence the positive eigenvalues of $\mathcal T$ are the numbers $j_{m-1,\ell}^{-2}$, $m,\ell\ge 1$, with eigenfunctions $\1_BJ_m(j_{m-1,\ell}r)\sin(m\vartheta)$.

The eigenvalues of $\cL$ on these eigenfunctions are $\frac12(k^{-2}-j_{m-1,\ell}^{-2})$, and $\cL=\frac{1}{2k^2}{\rm Id}$ on the kernel of $\mathcal T$. We now use two classical properties of the zeros of Bessel functions \cite[(10.21.2) and \S 10.21(iv)]{DLMF}: the interlacing $j_{0,1}<j_{1,1}<j_{0,2}$, and the fact that $j_{\nu,1}$ is increasing in $\nu$. They show that $j_{m-1,\ell}<k=j_{1,1}$ only for $(m,\ell)=(1,1)$, and that $j_{m-1,\ell}=k$ only for $(m,\ell)=(2,1)$. Since the eigenvalues of $\mathcal T$ accumulate only at zero, $\cL$ has exactly one negative eigenvalue $-\mu_-$, simple, with eigenfunction $e_-=\1_BJ_1(j_{0,1}r)\sin\vartheta$, its kernel is spanned by $\1_BJ_2(kr)\sin(2\vartheta)$, and its remaining spectrum is contained in $[\mu_0,+\infty)$ for some $\mu_0>0$. Finally, the identity $\pa_1(J_1(kr)\sin\vartheta)=-\frac k2J_2(kr)\sin(2\vartheta)$ shows that $\pa_1Q_L=\frac{k^2}{J_0(k)}\1_BJ_2(kr)\sin(2\vartheta)$ spans the kernel. If $f=k^2\1_BGf$, then $f$ is supported in $B$ and $\cL f=0$, which proves the last claim of (1).\\

\noindent{\bf Step 2:} Coercivity under the constraints $x_2\1_B$ and $\pa_1Q_L$. Let $f\in L^2(\Rp)$ with $(f,\pa_1Q_L)=0$ and $(f,x_2\1_B)=0$. By \fref{lambnegative}, the second condition reads $(\cL f,Q_L)=0$. Since $(Q_L,\pa_1Q_L)=\frac12\int_{\Rp}\pa_1(Q_L^2)=0$, both $f$ and $Q_L$ are orthogonal to the kernel, and we decompose
$$f=a e_-+f_+,\ \ Q_L=b e_-+Q_+,$$
where $f_+,Q_+$ belong to the spectral subspace of $[\mu_0,+\infty)$, and where $e_-$ is now normalized in $L^2$. The constraint becomes $\mu_-ab=(\cL f_+,Q_+)$, and the Cauchy--Schwarz inequality for the nonnegative form $(\cL\cdot,\cdot)$ on this subspace yields
$$\mu_-a^2\le \vartheta_0(\cL f_+,f_+),\ \ \vartheta_0:=\frac{(\cL Q_+,Q_+)}{\mu_-b^2}.$$
By \fref{lambnegative}, $-\mu_-b^2+(\cL Q_+,Q_+)=(\cL Q_L,Q_L)<0$, so that $b\neq 0$ and $\vartheta_0<1$. Therefore
$$(\cL f,f)=-\mu_-a^2+(\cL f_+,f_+)\ge (1-\vartheta_0)(\cL f_+,f_+)\ge (1-\vartheta_0)\mu_0\|f_+\|^2_{L^2},$$
and, since $a^2\le\mu_-^{-1}\vartheta_0\|\cL\|\|f_+\|^2_{L^2}$, we obtain $(\cL f,f)\ge \mu_2\|f\|_{L^2}^2$ for some $\mu_2>0$.\\

\noindent{\bf Step 3:} Replacing $\pa_1Q_L$ by $\varphi$. Let $f\in L^2(\Rp)$ with $(f,x_2\1_B)=(f,\varphi)=0$. We decompose $f=a\pa_1Q_L+g$ with $(g,\pa_1Q_L)=0$. Since $x_2\1_B$ is even in $x_1$ and $\pa_1Q_L$ is odd in $x_1$, $(g,x_2\1_B)=(f,x_2\1_B)=0$. The condition $(f,\varphi)=0$ gives $|a||(\pa_1Q_L,\varphi)|=|(g,\varphi)|\le\|\varphi\|_{L^2}\|g\|_{L^2}$, so that $\|f\|_{L^2}\lesssim \|g\|_{L^2}$. Since $\cL\pa_1Q_L=0$, Step 2 yields $(\cL f,f)=(\cL g,g)\ge \mu_2\|g\|^2_{L^2}\gtrsim\|f\|^2_{L^2}$. This is \fref{coercivityconstraint} for $R=1$.\\

\noindent{\bf Step 4:} Disks of radius close to one. By dominated convergence, the Hilbert--Schmidt norm of the operator $\1_{B_R}G\1_{B_R}-\1_BG\1_B$ tends to zero as $R\downarrow 1$, so that
$$\lim_{R\downarrow 1}\left(\|\cL_R-\cL\|_{L^2\to L^2}+\|x_2\1_{B_R}-x_2\1_B\|_{L^2}\right)=0.$$ Let $f$ satisfy the constraints of \fref{coercivityconstraint} for some $R\ge 1$, and let
$$g=f-\frac{(f,x_2\1_B)}{\|x_2\1_B\|^2_{L^2}}x_2\1_B.$$
Then $(g,x_2\1_B)=0$, $(g,\varphi)=(f,\varphi)=0$ since $(x_2\1_B,\varphi)=0$ by parity, and $\|f-g\|_{L^2}\lesssim |(f,x_2\1_B-x_2\1_{B_R})|\le o_{R\to 1}(1)\|f\|_{L^2}$. Step 3 applied to $g$ and the convergence of $\cL_R$ to $\cL$ yield \fref{coercivityconstraint} for $R$ close to one.
\end{proof}

\subsection{A smooth traveling dipole}

The Lamb dipole is only Lipschitz: its vorticity has a jump of its gradient across the unit circle, which comes from the corner of the function $s\mapsto k^2s_+$ in \fref{lambrelation}. Mollifying $Q_L$ would destroy the traveling wave structure. Instead, we keep the structure $Q=F(\Psi)$, which automatically produces a traveling wave since $U-e_1=-\nabla^\perp\Psi$, and we replace $k^2s_+$ by a smooth function $F_\sigma$ which is flat at $s=0$. The equation $Q=F_\sigma(GQ-x_2)$ is a semilinear elliptic problem, which we solve by a fixed point argument near $Q_L$. Its linearization at $Q_L$ is ${\rm Id}-k^2\1_BG$, whose kernel is spanned by $\pa_1Q_L$ by Lemma \ref{lemmaspectral}: this degeneracy comes from the invariance under horizontal translations, and it is removed by imposing that the solution is even in $x_1$.

\begin{proposition}[Smooth traveling dipoles]
\label{propdipole}
Let $R_0$ and $\varphi$ be as in Lemma \ref{lemmaspectral}. There exist $\sigma_0>0$ and $\nu>0$ such that for every $\sigma\in(0,\sigma_0]$ there exist $Q_\sigma\in C^\infty_c(\R^2)$ and $F_\sigma\in C^\infty(\R)$ with the following properties, where $U_\sigma=K*Q_\sigma$ and $\Psi_\sigma=GQ_\sigma-x_2$.
\begin{enumerate}
\item {\em Symmetry and support}: $Q_\sigma$ is odd in $x_2$, even in $x_1$, nonnegative on $\Rp$, and $\supp Q_\sigma\subset B_{R_0}$.
\item {\em Traveling wave}: $Q_\sigma=F_\sigma(\Psi_\sigma)$ on $\Rp$, where
\be
\label{propertiesF}
F_\sigma=0\ \ \mbox{on}\ \ (-\infty,0],\ \ 0<F_\sigma'\le k^2\ \ \mbox{on}\ \ (0,+\infty),\ \ F_\sigma^{(j)}(0)=0\ \ \mbox{for all}\ \ j\ge 0.
\ee
In particular $(U_\sigma-e_1)\cdot\nabla Q_\sigma=0$, and $Q_\sigma(x-te_1)$ is a smooth solution of \fref{euler}.
\item {\em Exterior bound}: $-\Psi_\sigma(x)\ge \nu x_2$ for all $x\in\Rp\setminus B_{R_0}$.
\item {\em Convergence to the Lamb dipole}: as $\sigma\to 0$,
\be
\label{convergencelamb}
\|Q_\sigma-Q_L\|_{L^\infty}+\|U_\sigma-U_L\|_{L^\infty}+\|\nabla U_\sigma-\nabla U_L\|_{L^\infty}\to 0,
\ee
and $(\pa_1Q_\sigma,\varphi)\to(\pa_1Q_L,\varphi)$.
\end{enumerate}
\end{proposition}

\begin{proof}
\noindent{\bf Step 1:} The function $F_\sigma$. We fix a nondecreasing function $\zeta\in C^\infty(\R)$ which vanishes on $(-\infty,0]$, is positive on $(0,+\infty)$, is equal to one on $[1,+\infty)$, and whose derivatives of all orders vanish at $0$. We let
$$F_\sigma(s)=k^2\int_0^{s}\zeta\left(\frac\tau\sigma\right)d\tau\ \ \mbox{for}\ \ s>0,\ \ F_\sigma(s)=0\ \ \mbox{for}\ \ s\le 0.$$
Then $F_\sigma$ satisfies \fref{propertiesF}, $F_\sigma'=k^2$ on $[\sigma,+\infty)$, and, since $0\le 1-\zeta(\tau/\sigma)\le \1_{\tau<\sigma}$,
\be
\label{Fclose}
|F_\sigma(s)-k^2s_+|\le k^2\sigma,\ \ s\in\R.
\ee

\noindent{\bf Step 2:} The linearized operator. We fix $p>2$, let $B_2^+=B_2\cap\Rp$, and let $L^p_{\rm ev}$ be the space of functions $g\in L^p(B_2^+)$ which are even in $x_1$, extended by zero outside $B_2^+$. For $g\in L^p_{\rm ev}$,
\be
\label{greenbounds}
\|Gg\|_{L^\infty(\Rp)}+\|\nabla Gg\|_{L^\infty(\Rp)}+\left\|\frac{Gg}{x_2}\right\|_{L^\infty(\Rp)}\le C_G\|g\|_{L^p(B_2^+)}.
\ee
Indeed, for $y\in B_2^+$, $|x-\overline y|\le |x-y|+4$ yields $0\le G(x,y)\lesssim 1+|\log|x-y||$, and $|\nabla_xG(x,y)|\lesssim |x-y|^{-1}$; both kernels belong to $L^{p'}(B_2^+)$ uniformly in $x$, since $p'<2$. The last term in \fref{greenbounds} is bounded by the second one because $Gg=0$ on $\{x_2=0\}$. By \fref{greenbounds} and the Arzel\`a--Ascoli theorem, the operator $g\mapsto \1_BGg$ is compact on $L^p_{\rm ev}$, and it preserves the parity in $x_1$. If $g\in L^p_{\rm ev}$ satisfies $g=k^2\1_BGg$, then $g\in L^2(\Rp)$, and Lemma \ref{lemmaspectral} shows that $g$ is a multiple of $\pa_1Q_L$, which is odd in $x_1$; hence $g=0$. By the Fredholm alternative, the operator $\mathcal A={\rm Id}-k^2\1_BG$ is invertible on $L^p_{\rm ev}$ with a bounded inverse.\\

\noindent{\bf Step 3:} Fixed point. We look for $Q_\sigma=Q_L+g$ on $B_2^+$, with $g\in L^p_{\rm ev}$. By \fref{lambrelation} and $GQ_L=\Psi_L+x_2$, the equation $Q_\sigma=F_\sigma(GQ_\sigma-x_2)$ on $B_2^+$ is equivalent to $g=F_\sigma(\Psi_L+Gg)-k^2(\Psi_L)_+$, that is,
\be
\label{fixedpoint}
\mathcal Ag=\mathcal N_\sigma(g),\qquad \mathcal N_\sigma(g):=F_\sigma(\Psi_L+Gg)-k^2(\Psi_L)_+-k^2\1_BGg\ \ \mbox{on}\ \ B_2^+.
\ee
Let $g_1,g_2\in L^p_{\rm ev}$ with $\|g_i\|_{L^p(B_2^+)}\le\rho$, so that $|Gg_i|\le C_G\rho$ by \fref{greenbounds}. At a point of $B_2^+$ where $\Psi_L>C_G\rho+\sigma$, we are in $B$ and both arguments $\Psi_L+Gg_i$ are larger than $\sigma$, so that $F_\sigma(\Psi_L+Gg_1)-F_\sigma(\Psi_L+Gg_2)=k^2G(g_1-g_2)$ and $\mathcal N_\sigma(g_1)=\mathcal N_\sigma(g_2)$. At a point where $\Psi_L<-C_G\rho-\sigma$, we are outside $B$, both arguments are negative, and again $\mathcal N_\sigma(g_1)=\mathcal N_\sigma(g_2)$. On the remaining set $\mathcal S_{\rho,\sigma}=\{x\in B_2^+:|\Psi_L(x)|\le C_G\rho+\sigma\}$, $|\mathcal N_\sigma(g_1)-\mathcal N_\sigma(g_2)|\le 2k^2|G(g_1-g_2)|$. Therefore
$$\|\mathcal N_\sigma(g_1)-\mathcal N_\sigma(g_2)\|_{L^p(B_2^+)}\le 2k^2C_G|\mathcal S_{\rho,\sigma}|^{\frac1p}\|g_1-g_2\|_{L^p(B_2^+)}.$$
The function $\Psi_L$ is continuous on $\overline{B_2^+}$ and vanishes only on the half circle $\{r=1\}$ and on the axis, which have zero measure, so that $|\mathcal S_{\rho,\sigma}|\to 0$ as $\rho+\sigma\to 0$. We fix $\rho_0>0$ and $\sigma_0>0$ such that $2k^2C_G\|\mathcal A^{-1}\||\mathcal S_{\rho_0,\sigma_0}|^{\frac1p}\le\frac12$. Moreover, $\|\mathcal N_\sigma(0)\|_{L^\infty}\le k^2\sigma$ by \fref{Fclose}, so that $\|\mathcal A^{-1}\mathcal N_\sigma(0)\|_{L^p(B_2^+)}\le C_A\sigma$. Decreasing $\sigma_0$ so that $2C_A\sigma_0\le\rho_0$, the map $g\mapsto\mathcal A^{-1}\mathcal N_\sigma(g)$ is a contraction of the closed ball of $L^p_{\rm ev}$ of radius $2C_A\sigma$, for every $\sigma\in(0,\sigma_0]$. Its fixed point $g_\sigma$ satisfies
\be
\label{boundgsigma}
\|g_\sigma\|_{L^p(B_2^+)}\le 2C_A\sigma.
\ee

\noindent{\bf Step 4:} Equation on the half-plane and support. We define $Q_\sigma=Q_L+g_\sigma$ on $B_2^+$ and $Q_\sigma=0$ on $\Rp\setminus B_2^+$, and we extend $Q_\sigma$ oddly to $\R^2$. Then $Q_\sigma=F_\sigma(\Psi_\sigma)$ on $B_2^+$. Let $R_1=\frac{1+R_0}{2}$. For $x\in\Rp\setminus B_{R_1}$, \fref{lambstream}, \fref{greenbounds} and \fref{boundgsigma} yield
$$\frac{\Psi_\sigma(x)}{x_2}=\frac{\Psi_L(x)}{x_2}+\frac{Gg_\sigma(x)}{x_2}\le \frac{1}{R_1^2}-1+2C_GC_A\sigma\le -\nu,\ \ \nu:=\frac12\left(1-\frac{1}{R_1^2}\right),$$
after decreasing $\sigma_0$. Hence $F_\sigma(\Psi_\sigma)=0$ on $\Rp\setminus B_{R_1}$. In particular $Q_\sigma=0$ on $B_2^+\setminus B_{R_1}$, the identity $Q_\sigma=F_\sigma(\Psi_\sigma)$ holds on the whole half-plane, $\supp Q_\sigma\subset\overline{B_{R_1}}\subset B_{R_0}$, and the exterior bound (3) holds.\\

\noindent{\bf Step 5:} Regularity and traveling wave. By \fref{greenbounds} and the explicit formula for $\Psi_L$, the function $\Psi_\sigma=\Psi_L+Gg_\sigma$ is Lipschitz on $\overline{\Rp}$. Since $\|F_\sigma'\|_{L^\infty}\le k^2$, the identity $Q_\sigma=F_\sigma(\Psi_\sigma)$ shows that $Q_\sigma$ is Lipschitz as well. The interior regularity of the Newtonian potential and the local boundary Schauder estimate \cite[Lemma 4.2 and Theorem 4.11]{GilbargTrudinger2001}, applied to $-\Delta\Psi_\sigma=Q_\sigma$ with $\Psi_\sigma=0$ on $\{x_2=0\}$, give $\Psi_\sigma\in C^{2,\gamma}_{\rm loc}(\overline{\Rp})$ for every $0<\gamma<1$. Composition with $F_\sigma$ and repeated application of local boundary regularity \cite[Theorem 6.19 and the subsequent local boundary remark, pp.~111--112]{GilbargTrudinger2001} then give $Q_\sigma,\Psi_\sigma\in C^\infty_{\rm loc}(\overline{\Rp})$. Since $\Psi_\sigma=0$ on $\{x_2=0\}$ and all the derivatives of $F_\sigma$ vanish at $0$, the chain rule shows that all the derivatives of $Q_\sigma$ vanish on $\{x_2=0\}$, so that its odd extension is smooth on $\R^2$. It is even in $x_1$ because $g_\sigma\in L^p_{\rm ev}$, and nonnegative on $\Rp$ because $F_\sigma\ge 0$. By \fref{velocityodd}, $U_\sigma-e_1=-\nabla^\perp\Psi_\sigma$ on $\Rp$, and hence
$$(U_\sigma-e_1)\cdot\nabla Q_\sigma=-F'_\sigma(\Psi_\sigma)\nabla^\perp\Psi_\sigma\cdot\nabla\Psi_\sigma=0\ \ \mbox{on}\ \ \Rp,$$
and on $\R^2$ by symmetry. Since the velocity of $Q_\sigma(\cdot-te_1)$ is $U_\sigma(\cdot-te_1)$, the function $Q_\sigma(x-te_1)$ solves \fref{euler}.\\

\noindent{\bf Step 6:} Convergence. On $B_2^+$, $g_\sigma=F_\sigma(\Psi_L+Gg_\sigma)-k^2(\Psi_L+Gg_\sigma)_++k^2\left[(\Psi_L+Gg_\sigma)_+-(\Psi_L)_+\right]$, so that \fref{Fclose}, \fref{greenbounds} and \fref{boundgsigma} yield $\|Q_\sigma-Q_L\|_{L^\infty}=\|g_\sigma\|_{L^\infty}\lesssim\sigma$. By \fref{velocityodd} and \fref{greenbounds}, $\|U_\sigma-U_L\|_{L^\infty}=\|\nabla Gg_\sigma\|_{L^\infty}\lesssim\sigma$. Moreover $|\nabla Q_\sigma|=F_\sigma'(\Psi_\sigma)|\nabla\Psi_\sigma|\le k^2(1+\|\nabla GQ_\sigma\|_{L^\infty})\lesssim 1$, so that $\|\nabla g_\sigma\|_{L^\infty}\lesssim 1$. We apply the logarithmic estimate \fref{logestimate} to the field $U_\sigma-U_L$, whose vorticity is the odd extension of $g_\sigma$. Since $s\mapsto s(1+\log_+(L/s))$ is nondecreasing on $(0,+\infty)$ for every $L>0$, we obtain
$$\|\nabla U_\sigma-\nabla U_L\|_{L^\infty}\lesssim \sigma+\sigma\left(1+\log_+\frac{C}{\sigma}\right)\to 0.$$
Finally, $(\pa_1Q_\sigma,\varphi)=-(Q_\sigma,\pa_1\varphi)\to-(Q_L,\pa_1\varphi)=(\pa_1Q_L,\varphi)$.
\end{proof}

\subsection{Orbital stability}
\label{subsectionorbital}

From now on, we fix $\sigma\in(0,\sigma_0]$ such that
$$|(\pa_1Q_\sigma,\varphi)|\ge\frac12|(\pa_1Q_L,\varphi)|;$$
it will be further decreased once in Section \ref{sectionfiniteenergy}, see \fref{choicesigma}. We drop the index $\sigma$ and write $Q$, $F$, $U$, $\Psi$. All the constants below depend on this fixed profile.

\medskip

\noindent{\it The energy--Casimir functional}. The orbital stability of $Q$ relies on a conserved functional for which $Q$ is a critical point with a coercive second variation, as in the energy--Casimir method of Arnold \cite{Arnold1965,Arnold1966}, see also \cite{HMRW1985,WanPulvirenti1985}. Since $F$ is a smooth increasing bijection of $[0,+\infty)$, we may define the convex Casimir function
\be
\label{defcasimir}
\Phi(s)=\int_0^sF^{-1}(\tau)d\tau,\ \ s\ge 0,
\ee
and, for bounded compactly supported functions $f\ge 0$ on $\Rp$, the energy--Casimir functional
\be
\label{defenergycasimir}
\cH(f)=\int_{\Rp}\Phi(f)dx-E(f)+I(f).
\ee
The function $\Phi$ is $C^1$, $\Phi(0)=\Phi'(0)=0$, and $\Phi(s)\le F^{-1}(M)s$ for $0\le s\le M$. Since $(F^{-1})'=1/F'(F^{-1})\ge k^{-2}$ on $(0,+\infty)$ by \fref{propertiesF}, the function $s\mapsto\Phi(s)-\frac{s^2}{2k^2}$ is convex, which means
\be
\label{casimirconvex}
\Phi(s_1)-\Phi(s_2)-\Phi'(s_2)(s_1-s_2)\ge\frac{(s_1-s_2)^2}{2k^2},\ \ s_1,s_2\ge 0.
\ee
For a smooth solution of \fref{euler} which is odd in $x_2$, nonnegative on $\Rp$ and compactly supported, $\int_{\Rp}\Phi(\omega(t))dx$ is conserved because $\omega(t)$ is the composition of $\omega_0$ with a measure preserving map, and hence $\cH(\omega(t))=\cH(\omega_0)$. Moreover $\cH$ is invariant under horizontal translations. The profile $Q$ is a critical point of $\cH$ in the following one-sided sense: on $\{Q>0\}$, the first variation $\Phi'(Q)-GQ+x_2=F^{-1}(Q)-\Psi$ vanishes, and on $\{Q=0\}$ it is equal to $-\Psi\ge 0$, which penalizes the vorticity which leaves the support of $Q$. The exterior bound (3) of Proposition \ref{propdipole} makes this penalization quantitative.

\begin{proposition}[Orbital stability of the smooth dipole]
\label{proporbital}
For $\omega_0\in C^\infty_c(\R^2)$, odd in $x_2$, nonnegative on $\Rp$, with $I(\omega_0)=I(Q)$, let
$$\cD(\omega_0):=\|\omega_0-Q\|_{L^2(\Rp)}^2+|\cH(\omega_0)-\cH(Q)|.$$
There exist $\cD_0>0$ and $C>0$, depending only on $Q$ and $\varphi$, such that the following holds. If $\cD(\omega_0)\le\cD_0$, then the solution $\omega$ of \fref{euler} admits a translation parameter $z\in C^1([0,+\infty))$ such that the perturbation
$$\eps(t,x):=\omega(t,x+z(t)e_1)-Q(x)$$
satisfies, for all $t\ge 0$:
\begin{enumerate}
\item {\em Orthogonality}: $(\eps(t),\varphi)=0$.
\item {\em Smallness}:
\be
\label{orbitalbound}
\|\eps(t)\|_{L^2(\Rp)}+\|x_2\eps(t)\|_{L^1(\Rp)}+\|K*\eps(t)\|_{L^2(\R^2)}\le C\cD(\omega_0)^{\frac12}.
\ee
\item {\em Translation}: $|z'(t)-1|\le C\cD(\omega_0)^{\frac12}$ and $|z(0)|\le C\|\omega_0-Q\|_{L^2(\Rp)}$.
\item {\em Equation}: $\eps$ is a smooth solution of
\be
\label{eqeps}
\pa_t\eps+(U+K*\eps-z'e_1)\cdot\nabla\eps=-\left(K*\eps+(1-z')e_1\right)\cdot\nabla Q.
\ee
\end{enumerate}
\end{proposition}

The constants in Proposition \ref{proporbital} do not depend on the support or on the higher norms of $\omega_0$: this is the reason why the profile $Q$ is fixed once and for all, and it is the key to the uniformity in Section \ref{sectionfiniteenergy}.

\begin{proof}
\noindent{\bf Step 1:} Expansion of the functional. Let $f$ be bounded, compactly supported, odd in $x_2$, with $Q+f\ge 0$ on $\Rp$. Since $E(Q+f)-E(Q)=(GQ,f)+E(f)$, $I(f)=(x_2,f)$ and $GQ-x_2=\Psi$,
$$\cH(Q+f)-\cH(Q)=\int_{\Rp}\left[\Phi(Q+f)-\Phi(Q)-\Psi f\right]dx-E(f).$$
On $\{Q>0\}$, $\Psi=F^{-1}(Q)=\Phi'(Q)$, and the integrand is bounded from below by $\frac{f^2}{2k^2}$ by \fref{casimirconvex}. On $\{Q=0\}$, we have $f\ge 0$, $\Psi\le 0$ since $F(\Psi)=0$, and the integrand $\Phi(f)-\Psi f$ is bounded from below by $\frac{f^2}{2k^2}+(-\Psi)f$. Therefore
\be
\label{expansionH}
\cH(Q+f)-\cH(Q)\ge\frac{1}{2k^2}\|f\|^2_{L^2}-E(f)+\mathcal J(f),
\ee
where
$$\mathcal J(f):=\int_{\Rp\cap\{Q=0\}}(-\Psi)f\,dx\ge\nu\int_{\Rp\setminus B_{R_0}}x_2f\,dx,$$
and the last inequality follows from $-\Psi f\ge 0$ on $\{Q=0\}$ and the exterior bound (3) of Proposition \ref{propdipole}.\\

\noindent{\bf Step 2:} Coercivity. Let $f$ be as in Step 1, with moreover $I(f)=0$, $(f,\varphi)=0$ and $\rho:=\|f\|_{L^2}\le 1$. We split $f=f_i+f_e$ with $f_i=\1_{B_{R_0}}f$ and $f_e=\1_{\Rp\setminus B_{R_0}}f$. Since $Q=0$ outside $B_{R_0}$, $f_e\ge 0$, and we let $N:=I(f_e)=\|x_2f_e\|_{L^1}\ge 0$. The constraint $I(f)=0$ reads $(f_i,x_2\1_{B_{R_0}})=-N$, so that $N\lesssim\rho$, and $(f_i,\varphi)=0$ since $\varphi$ is supported in $B$. The function $g=f_i+N\|x_2\1_{B_{R_0}}\|^{-2}_{L^2}x_2\1_{B_{R_0}}$ satisfies both constraints of \fref{coercivityconstraint} with $R=R_0$, because $(x_2\1_{B_{R_0}},\varphi)=0$ by parity. Since $\|f_i-g\|_{L^2}\lesssim N$, Lemma \ref{lemmaspectral} yields $(\cL_{R_0}f_i,f_i)\ge\frac{\mu_1}2\|f_i\|^2_{L^2}-CN^2$. We now expand
$$\frac{\|f\|^2_{L^2}}{2k^2}-E(f)=(\cL_{R_0}f_i,f_i)+\frac{\|f_e\|^2_{L^2}}{2k^2}-E(f_e)-(Gf_i,f_e),$$
since $f_i$ is supported in $B_{R_0}$.
By Lemma \ref{lemmaenergy}, $E(f_e)\lesssim N\rho$. Since $G$ is positive, $|(Gf_i,f_e)|\le 2E(f_i)^{\frac12}E(f_e)^{\frac12}$, and $E(f_i)\lesssim\|f_i\|^2_{L^2}$ because $f_i$ is supported in the fixed disk $B_{R_0}$. Hence $|(Gf_i,f_e)|\lesssim\rho(N\rho)^{\frac12}\le\eta\rho^2+C_\eta N\rho$ for every $\eta>0$. Choosing $\eta$ small and using $N\lesssim\rho$, we obtain $\frac{1}{2k^2}\|f\|^2_{L^2}-E(f)\ge\mu_3\rho^2-C_4N\rho$ for some $\mu_3>0$. Since $\mathcal J(f)\ge\nu N$ by \fref{expansionH}, we have $C_4N\rho\le\frac12\mathcal J(f)$ if $\rho\le\rho_1:=\min\{1,\frac{\nu}{2C_4}\}$, and \fref{expansionH} yields
\be
\label{coercivityH}
\cH(Q+f)-\cH(Q)\ge\mu_3\|f\|^2_{L^2}+\frac\nu2\int_{\Rp\setminus B_{R_0}}x_2f\,dx\ \ \mbox{if}\ \ \|f\|_{L^2}\le\rho_1.
\ee

\noindent{\bf Step 3:} Modulation. The map $(z,f)\mapsto(f(\cdot+ze_1),\varphi)-(Q,\varphi)=(f,\varphi(\cdot-ze_1))-(Q,\varphi)$ is $C^1$ on $\R\times L^2(\Rp)$, and its derivative with respect to $z$ at a point where $f(\cdot+ze_1)=Q$ is $-(Q,\pa_1\varphi)=(\pa_1Q,\varphi)\neq 0$. By the implicit function theorem, there exist $r_1>0$ and $C>0$ such that, if $\|f(\cdot+z_0e_1)-Q\|_{L^2}\le r_1$ for some $z_0\in\R$, there exists a unique $z$ with $|z-z_0|\le Cr_1$ and $(f(\cdot+ze_1)-Q,\varphi)=0$; moreover $|z-z_0|\le C\|f(\cdot+z_0e_1)-Q\|_{L^2}$, and $z$ is a $C^1$ function of $f$.\\

\noindent{\bf Step 4:} Proof of (1), (2) and (4). Let $\cD=\cD(\omega_0)\le\cD_0$. Step 3 with $z_0=0$ defines $z(0)$ with $|z(0)|\le C\|\omega_0-Q\|_{L^2}$, and then
$$\|\eps(0)\|_{L^2}\le\|\omega_0-Q\|_{L^2}+|z(0)|\|\pa_1Q\|_{L^2}\lesssim\cD^{\frac12}.$$
Since $t\mapsto\omega(t)\in L^2(\Rp)$ is $C^1$, Step 3 defines $z\in C^1$ on a maximal interval $[0,t^*)$ on which $\|\eps(t)\|_{L^2}<\rho_2:=\min\{\rho_1,r_1/2\}$. On this interval, $Q+\eps(t)\ge 0$, $I(\eps(t))=I(\omega(t))-I(Q)=0$ by the conservation and the translation invariance of the impulse, and $(\eps(t),\varphi)=0$. The conservation and the translation invariance of $\cH$ give $\cH(Q+\eps(t))-\cH(Q)=\cH(\omega_0)-\cH(Q)\le\cD$, and \fref{coercivityH} yields
\be
\label{boundorbitalproof}
\mu_3\|\eps(t)\|^2_{L^2}+\frac\nu2\int_{\Rp\setminus B_{R_0}}x_2\eps(t)\,dx\le\cD.
\ee
Choosing $\cD_0$ such that $(\cD_0/\mu_3)^{\frac12}\le\rho_2/2$ and $\|\eps(0)\|_{L^2}<\rho_2$, a continuity argument based on Step 3 and on the bound on $z'$ proved in Step 5 shows that $t^*=+\infty$. Next, $\|x_2\eps\|_{L^1}\le\int_{B_{R_0}}x_2|\eps|+\int_{\Rp\setminus B_{R_0}}x_2\eps\lesssim\|\eps\|_{L^2}+\cD\lesssim\cD^{\frac12}$, and by \fref{defimpulseenergy} and Lemma \ref{lemmaenergy}, $\|K*\eps\|^2_{L^2(\R^2)}=4E(\eps)\lesssim\|x_2\eps\|_{L^1}\|\eps\|_{L^2}\lesssim\cD$. This proves (1) and (2). The function $\widetilde\omega(t,x)=\omega(t,x+z(t)e_1)=Q+\eps$ solves
$$\pa_t\widetilde\omega+\left(K*\widetilde\omega-z'e_1\right)\cdot\nabla\widetilde\omega=0,$$
and $K*\widetilde\omega=U+K*\eps$. Subtracting $(U-e_1)\cdot\nabla Q=0$ yields \fref{eqeps}.\\

\noindent{\bf Step 5:} Speed of translation. We differentiate $(\eps,\varphi)=0$ in time. Since $\varphi$ is compactly supported in $\Rp$ and $K*\widetilde\omega-z'e_1$ is divergence free, an integration by parts gives $0=(\pa_t\widetilde\omega,\varphi)=(\widetilde\omega,(U+K*\eps-z'e_1)\cdot\nabla\varphi)$. Using $(Q,(U-e_1)\cdot\nabla\varphi)=-((U-e_1)\cdot\nabla Q,\varphi)=0$, we obtain
$$(z'-1)(\widetilde\omega,\pa_1\varphi)=(\eps,(U-e_1)\cdot\nabla\varphi)+(\widetilde\omega\,K*\eps,\nabla\varphi).$$
The factor $(\widetilde\omega,\pa_1\varphi)=-(\pa_1Q,\varphi)+(\eps,\pa_1\varphi)$ is larger than $\frac12|(\pa_1Q,\varphi)|$ in absolute value for $\cD_0$ small. The right-hand side is bounded by $C\|\eps\|_{L^2}+C\|K*\eps\|_{L^2}(1+\|\eps\|_{L^2})\lesssim\cD^{\frac12}$ by \fref{orbitalbound}. This proves (3).
\end{proof}

\section{Finite energy solutions: proof of Theorem \ref{thmmain}}
\label{sectionfiniteenergy}

In this section, we first turn the orbital stability of $Q$ into an exponential growth of the vorticity gradient on a time interval $[0,T]$, with initial data whose norms are bounded independently of $T$. We then rescale these solutions, place them far apart from each other, and choose the parameters so as to prove Theorem \ref{thmmain}.

\subsection{Exponential stretching behind the dipole}

\noindent{\it The geometry of the rear stagnation point}. By \fref{lambaxis}, the function $x_1\mapsto U_{L,1}(x_1,0)-1$ vanishes at $x_1=-1$ with derivative $2$, and $\nabla U_L$ is continuous since $Q_L$ is Lipschitz. We fix $r_*\in(0,\frac14]$ such that $\pa_1U_{L,1}(x_1,0)\ge\frac32$ for $x_1\in[-1-r_*,-1+r_*]$, and we decrease $\sigma$ in Subsection \ref{subsectionorbital} so that, by \fref{convergencelamb},
\be
\label{choicesigma}
\|U-U_L\|_{L^\infty}+\|\nabla U-\nabla U_L\|_{L^\infty}\le\frac{r_*}{4}.
\ee
Integrating $\pa_1U_{L,1}\ge\frac32$ from $-1$, we obtain on the segment $J_*=[-1-r_*,-1+r_*]$ of the axis:
\be
\label{reargeometry}
\pa_1U_1(x_1,0)\ge\frac54\ \ \mbox{for}\ \ x_1\in J_*,\qquad \pm\left(U_1(-1\pm r_*,0)-1\right)\ge\frac54r_*.
\ee
In the frame moving with the dipole, the horizontal velocity on the axis therefore points outwards at both ends of $J_*$, and it expands $J_*$ at a rate at least $\frac54$.

\medskip

\noindent{\it The initial data}. We fix a function $\chi\in C^\infty_c(B_3)$, even in $x_2$, with $0\le\chi\le 1$ and $\chi=1$ on $[-2,0]\times[-1,1]$, and we let $P(x)=x_2\chi(x)$. Then $P$ is odd in $x_2$, nonnegative on $\Rp$, and $\pa_2P(x_1,0)=1$ for $x_1\in[-2,0]$. For $T\ge 1$ and a constant $\delta_0\in(0,1]$ to be fixed, we let
\be
\label{defdataT}
\delta_T=\delta_0(1+T)^{-6},\qquad \omega_{0,T}=m_T(Q+\delta_TP),\qquad m_T=\frac{I(Q)}{I(Q)+\delta_TI(P)}.
\ee
Note that $m_T\in[\frac12,1]$ for $\delta_0$ small.
The data $\omega_{0,T}$ are smooth, odd in $x_2$, nonnegative on $\Rp$, supported in $B_3$, and they have the same impulse as $Q$. The size $\delta_T$ of the perturbation is chosen so small that the perturbation of the dipole remains under control on $[0,T]$, see Step 2 below, but only polynomially small in $T$, so that the exponential growth $e^{t}$ wins.

\begin{proposition}[Exponential stretching on a time interval of length $T$]
\label{propstretching}
There exist $\delta_0\in(0,1]$ and $c_*>0$, and for every integer $m\ge 0$ a constant $C_m$, such that for every $T\ge 1$ the solution $\omega_T$ of \fref{euler} with initial datum $\omega_{0,T}$ satisfies the following.
\begin{enumerate}
\item {\em Bounds on the data uniform in $T$}: $\|\omega_{0,T}\|_{L^1}+\|\omega_{0,T}\|_{H^m}+\|K*\omega_{0,T}\|_{H^m}\le C_m$.
\item {\em Growth along a curve}: there exists a $C^1$ curve $q_T:[0,T]\to\R\times\{0\}$ such that
\be
\label{growthcurve}
\pa_2\omega_T(s,q_T(s))\ge c_*(1+T)^{-6}e^{s},\ \ 0\le s\le T.
\ee
\end{enumerate}
\end{proposition}

\begin{proof}
\noindent{\bf Step 1:} Orbital stability. Since $1-m_T\lesssim\delta_T$, we have $\|\omega_{0,T}-Q\|_{L^1\cap L^\infty}\lesssim\delta_T$. The function $\Phi$ is Lipschitz on the bounded range of $\omega_{0,T}$ and $Q$, and $G$ is bounded on $L^2(B_3\cap\Rp)$, so that $|\cH(\omega_{0,T})-\cH(Q)|\lesssim\delta_T$ by \fref{defenergycasimir}. Hence $\cD(\omega_{0,T})\le C\delta_T\le C\delta_0\le\cD_0$ for $\delta_0$ small, and Proposition \ref{proporbital} applies: the solution $\omega_T$ admits a translation parameter $z_T$, which we denote by $z$, and a perturbation $\eps(t,x)=\omega_T(t,x+z(t)e_1)-Q(x)$ such that, for all $t\ge 0$,
\be
\label{boundsT}
\|\eps(t)\|_{L^2(\Rp)}+\|K*\eps(t)\|_{L^2(\R^2)}+|z'(t)-1|\lesssim\delta_T^{\frac12},\qquad |z(0)|\lesssim\delta_T.
\ee

\noindent{\bf Step 2:} $L^\infty$ bounds. Since the vorticity is transported,
$$\|\eps(t)\|_{L^1}+\|\eps(t)\|_{L^\infty}\le\|\omega_{0,T}\|_{L^1\cap L^\infty}+\|Q\|_{L^1\cap L^\infty}\lesssim 1.$$ The velocity bound \fref{velocitybound} with $\rho=\delta_T^{1/2}$ and \fref{boundsT} yield
\be
\label{boundvelocityeps}
\|K*\eps(t)\|_{L^\infty}\lesssim\delta_T^{\frac12}\left(1+|\log\delta_T|\right)^{\frac12}.
\ee
The perturbation $\eps$ is transported by a divergence-free field with the source term $-(K*\eps+(1-z')e_1)\cdot\nabla Q$ in \fref{eqeps}. Moreover $\|\eps(0)\|_{L^\infty}\le\|\omega_{0,T}-Q\|_{L^\infty}+|z(0)|\|\nabla\omega_{0,T}\|_{L^\infty}\lesssim\delta_T$. Lemma \ref{lemmatransport}, \fref{boundsT} and \fref{boundvelocityeps} yield, for $0\le t\le T$,
$$\|\eps(t)\|_{L^\infty}\lesssim\delta_T+T\delta_T^{\frac12}(1+|\log\delta_T|)^{\frac12}.$$
We now insert $\delta_T=\delta_0(1+T)^{-6}$ and use $1+|\log\delta_T|\le(1+|\log \delta_0|)(1+6\log(1+T))$ and $\sup_{T\ge 1}(1+T)^{-1}(1+6\log(1+T))^{\frac12}<+\infty$. We obtain
\be
\label{boundepsinfty}
\|\eps(t)\|_{L^\infty}\le\frac{d(\delta_0)}{1+T},\ \ 0\le t\le T,\qquad d(\delta_0):=C \delta_0^{\frac12}\left(1+|\log \delta_0|\right)^{\frac12},
\ee
where $C$ is independent of $T$ and $\delta_0$, and $d(\delta_0)\to 0$ as $\delta_0\to 0$. In words: the $L^2$ smallness given by the orbital stability is converted into an $L^\infty$ smallness which beats the length $T$ of the time interval.\\

\noindent{\bf Step 3:} Bootstrap on the velocity gradient of the perturbation. We claim that, for $\delta_0$ small enough,
\be
\label{bootstrapgradient}
\|\nabla(K*\eps)(t)\|_{L^\infty}\le\frac18,\ \ 0\le t\le T.
\ee
Assume that $\|\nabla(K*\eps)(s)\|_{L^\infty}\le\frac14$ for $s\in[0,t]$, with $t\le T$. Then the velocity $u_T(s,x+z(s)e_1)=U(x)+K*\eps(s,x)$ satisfies $\|\nabla u_T(s)\|_{L^\infty}\le M_0:=\|\nabla U\|_{L^\infty}+\frac14$, a constant independent of $T$ and $\delta_0$. By \fref{gradientgronwall}, $\|\nabla\omega_T(s)\|_{L^\infty}\lesssim e^{M_0s}$, and hence $\|\nabla\eps(s)\|_{L^\infty}\lesssim e^{M_0T}$. We apply the logarithmic estimate \fref{logestimate} to $v=K*\eps(s)$. Since $s\mapsto s(1+\log_+(L/s))$ is nondecreasing, \fref{boundepsinfty} yields
\begin{align*}
\|\nabla(K*\eps)(s)\|_{L^\infty}&\lesssim\|K*\eps(s)\|_{L^\infty}+\frac{d(\delta_0)}{1+T}\left(1+\log_+\frac{Ce^{M_0T}(1+T)}{d(\delta_0)}\right)\\
&\lesssim\delta_T^{\frac12}(1+|\log\delta_T|)^{\frac12}+d(\delta_0)\left(1+|\log d(\delta_0)|\right),
\end{align*}
because $\frac{1+M_0T+\log(1+T)}{1+T}\lesssim 1$. The logarithm of the Lipschitz norm of $\eps$ grows linearly in $T$, but it is multiplied by the factor $(1+T)^{-1}$ of \fref{boundepsinfty}. The right-hand side tends to zero as $\delta_0\to 0$, uniformly in $T\ge 1$, so that it is smaller than $\frac18$ for $\delta_0$ small. The same computation at $s=0$ shows that the bootstrap assumption holds at $t=0$, and a continuity argument proves \fref{bootstrapgradient}. Decreasing $\delta_0$ further, \fref{boundsT} and \fref{boundvelocityeps} also give
\be
\label{smallvelocityaxis}
\|K*\eps(t)\|_{L^\infty}+|z'(t)-1|\le\frac{r_*}{8},\ \ t\ge 0.
\ee

\noindent{\bf Step 4:} Trajectories near the rear stagnation point. Since $\omega_T$ is odd in $x_2$, it vanishes on the axis, and the vertical component of $u_T$ vanishes there. We denote by
$$v(t,x_1)=U_1(x_1,0)+(K*\eps)_1(t,x_1,0)-z'(t)$$
the horizontal velocity on the axis in the frame moving with the dipole. By \fref{reargeometry}, \fref{bootstrapgradient} and \fref{smallvelocityaxis}, for $0\le t\le T$,
\be
\label{axisvelocity}
\pa_1v(t,x_1)\ge 1\ \ \mbox{for}\ \ x_1\in J_*,\qquad v(t,-1-r_*)<0<v(t,-1+r_*).
\ee
Let $s\in[0,T]$, and let $y_s$ be the solution of $y_s'(t)=v(t,y_s(t))$ on $[0,s]$ with $y_s(s)=-1$. By \fref{axisvelocity}, backwards in time the field $-v$ points into $J_*$ at both endpoints, so that $y_s(t)\in J_*$ for $0\le t\le s$. The point $X_s(t)=(y_s(t)+z(t),0)$ is a trajectory of $u_T$, since $\frac{d}{dt}(y_s+z)=U_1+(K*\eps)_1=u_{T,1}(t,X_s(t))$. Moreover, $|X_s(0)-(y_s(0),0)|=|z(0)|\lesssim\delta_T\le r_*$ for $\delta_0$ small, so that $X_s(0)\in[-2,0]\times\{0\}$, where $\pa_2P=1$.\\

\noindent{\bf Step 5:} Growth of the vertical derivative on the axis. Let $g(t,x_1)=\pa_2\omega_T(t,x_1,0)$. We differentiate the transport equation \fref{euler} with respect to $x_2$ and evaluate it on the axis, where $\omega_T=\pa_1\omega_T=0$, $u_{T,2}=0$ and $\pa_2u_{T,2}=-\pa_1u_{T,1}$ by incompressibility. We obtain
$$\pa_tg+u_{T,1}\pa_1g=(\pa_1u_{T,1})g\ \ \mbox{on}\ \ \{x_2=0\}.$$
Along the trajectory $X_s$, $\pa_1u_{T,1}(t,X_s(t))=\pa_1v(t,y_s(t))\ge 1$ by \fref{axisvelocity}. Since $Q\ge 0$ on $\Rp$ and $Q=0$ on the axis, $\pa_2Q\ge 0$ on the axis, and hence $g(0,X_s(0))\ge m_T\delta_T\pa_2P(X_s(0))\ge\frac{\delta_T}{2}$. Therefore
$$g(s,X_s(s))=g(0,X_s(0))\exp\left(\int_0^s\pa_1u_{T,1}(t,X_s(t))dt\right)\ge\frac{\delta_T}2e^{s}.$$
This proves \fref{growthcurve} with $q_T(s)=X_s(s)=(z(s)-1,0)$, which is $C^1$, and $c_*=\frac{\delta_0}{2}$.\\

\noindent{\bf Step 6:} Bounds on the data. The bound on $\|\omega_{0,T}\|_{L^1}+\|\omega_{0,T}\|_{H^m}$ follows from \fref{defdataT}. By Lemma \ref{lemmaenergy},
$$\|K*\omega_{0,T}\|^2_{L^2(\R^2)}=4E(\omega_{0,T})\lesssim\|x_2\omega_{0,T}\|_{L^1}\|\omega_{0,T}\|_{L^2}\lesssim 1,$$
and the identity $\|\nabla(K*f)\|_{L^2}=\|f\|_{L^2}$ gives the bounds on the higher derivatives of the velocity.
\end{proof}

\subsection{Rescaled components}

We now use the two scalings \fref{scalings}. For $T\ge 1$ and parameters $a,r\in(0,1]$, the function
\be
\label{defcomponent}
\omega_{T,a,r}(t,x)=a\,\omega_T\left(at,\frac xr\right)
\ee
is a global smooth solution of \fref{euler}, with velocity $V_{T,a,r}(t,x)=ar\,u_T(at,\frac xr)$. The amplitude $a$ slows down the time, and the size $r$ makes the support small. By Proposition \ref{propstretching}, using $\|f(\cdot/r)\|_{L^p}=r^{2/p}\|f\|_{L^p}$:
\be
\label{componentbounds}
\left\{\begin{array}{l}
\supp\omega_{T,a,r}(0)\subset B_{3r},\ \ \|\omega_{T,a,r}(0)\|_{L^1}+\|V_{T,a,r}(0)\|_{L^2}\lesssim ar^2,\\[1mm]
\|D^m\omega_{T,a,r}(0)\|_{L^2}\lesssim_m ar^{1-m},\ \ m\ge 0,
\end{array}\right.
\ee
and, along the curve $t\mapsto rq_T(at)$,
\be
\label{componentgrowth}
|\nabla\omega_{T,a,r}(t,rq_T(at))|\ge c_*\frac ar(1+T)^{-6}e^{at},\ \ 0\le t\le\frac Ta.
\ee
Moreover $V_{T,a,r}\in C([0,S];H^m)$ for all $S>0$ and $m\ge 0$, by \cite[Theorems 7.1 and 7.25]{BahouriCheminDanchin2011}, since its initial velocity is divergence free and belongs to every Sobolev space.

\subsection{Proof of Theorem \ref{thmmain}}

We fix the function $\Gamma$. The proof consists in choosing two families of components: the first family produces the sequence of times in \fref{growthsequence}, and the second one produces time intervals which cover all large times, for \fref{growthalltime}.\\

\noindent{\bf Step 1:} The first family. For $j\ge 1$, we choose the amplitude $e^{-j}$, the size $1$, and a length $\tau_j\ge j^2$ such that
\be
\label{choicetauj}
\log\left(c_*e^{-j}(1+\tau_j)^{-6}\right)\ge-\frac{\tau_j}{4},\qquad\Gamma\left(e^j\tau_j\right)\le\frac{\tau_j}{j^2},\qquad e^j\tau_j>e^{j-1}\tau_{j-1}.
\ee
The last inequality is imposed only for $j\ge 2$. This is possible: for fixed $j$, the first left-hand side is $O_j(1)-6\log(1+\tau)$, and $\Gamma(e^j\tau)/\tau=e^j\Gamma(e^j\tau)/(e^j\tau)\to 0$ as $\tau\to+\infty$ by assumption on $\Gamma$. We let $t_j=e^j\tau_j$, which increases to $+\infty$. By \fref{componentgrowth} at time $t=t_j$, the component $\omega_{\tau_j,e^{-j},1}$ satisfies, at the point $x^{(j)}=q_{\tau_j}(\tau_j)$,
\be
\label{growthfirstfamily}
|\nabla\omega_{\tau_j,e^{-j},1}(t_j,x^{(j)})|\ge c_*e^{-j}(1+\tau_j)^{-6}e^{\tau_j}\ge e^{\frac34\tau_j}.
\ee

\noindent{\bf Step 2:} The second family. We fix an integer $j_0\ge 4$ and, for $j\ge 0$, we let
\be
\label{choicearj}
m_j=j+j_0,\qquad a_j=\frac{e^{\sqrt{\log m_j}}}{m_j},\qquad r_j=e^{-\sqrt{\log m_j}}.
\ee
The amplitudes $a_j\le 1$ decrease to zero, slower than $m_j^{-1}$, and $a_j/a_{j-1}\to 1$. We fix $E_0\ge 10$, and we define the towers $E_{j+1}=e^{E_j}$, and the lengths
\be
\label{choiceTj}
T_0=9E_0^2,\ \ T_j=\exp\left(T_{j-1}^{\frac12}\right),\ \ S_j=\frac{T_j}{a_j},\ \ J_j=[S_{j-1},S_j],
\ee
where the recursion for $T_j$ is used for $j\ge 1$, $S_j$ is defined for $j\ge 0$, and $J_j$ for $j\ge 1$.
The times $S_j$ increase to $+\infty$, so that the intervals $J_j$ cover $[S_0,+\infty)$. By \fref{componentgrowth}, the component $\omega_{T_j,a_j,r_j}$ satisfies
\be
\label{growthsecondfamily}
|\nabla\omega_{T_j,a_j,r_j}(t,r_jq_{T_j}(a_jt))|\ge \frac{c_*a_j}{r_j(1+T_j)^{6}}e^{a_jt},\ \ t\in J_j.
\ee
The choice \fref{choicearj} is dictated by two competing requirements. By \fref{componentbounds},
\begin{align*}
&\|D^m\omega_{T_j,a_j,r_j}(0)\|_{L^2}\lesssim_m \frac{e^{m\sqrt{\log m_j}}}{m_j},\\
&\|\omega_{T_j,a_j,r_j}(0)\|_{L^1}+\|V_{T_j,a_j,r_j}(0)\|_{L^2}\lesssim \frac{e^{-\sqrt{\log m_j}}}{m_j}.
\end{align*}
Since $e^{m\sqrt{\log m_j}}=m_j^{o(1)}$, and since the change of variables $x=e^{v^2}$ shows that $\int^{+\infty}\frac{dx}{xe^{\sqrt{\log x}}}=\int^{+\infty}2ve^{-v}dv<+\infty$, we obtain
\be
\label{summability}
\sum_{j\ge 1}\|D^m\omega_{T_j,a_j,r_j}(0)\|^2_{L^2}+\sum_{j\ge 1}\left(\|\omega_{T_j,a_j,r_j}(0)\|_{L^1}+\|V_{T_j,a_j,r_j}(0)\|_{L^2}\right)<+\infty
\ee
for all $m\ge 0$.
On the other hand, the amplitude $a_j$, which is the rate of growth in physical time, decays more slowly than $m_j^{-1}$. The Sobolev norms are only square summable, which is enough because the supports will be disjoint.\\

\noindent{\bf Step 3:} Assembly. We order the components as follows: for $j\ge 1$, the component of index $n=2j-1$ is $\omega_{\tau_j,e^{-j},1}$, and the component of index $n=2j$ is $\omega_{T_j,a_j,r_j}$. We denote by $\omega_n$ and $V_n$ the vorticity and the velocity of the $n$-th component. We apply Lemma \ref{lemmaassembly} with $s_n=4$, $\eps_n=2^{-n}C_0^{-1}$, the nondecreasing horizons $S_{2j-1}^\sharp=S_{2j}^\sharp=\max\{t_1,\dots,t_j,S_1,\dots,S_j\}$, and the compact sets
$$\Sigma_{2j-1}=\{(t_j,x^{(j)})\},\qquad \Sigma_{2j}=\left\{(t,r_jq_{T_j}(a_jt)):t\in J_j\right\}.$$
We choose the distances $d_n\ge\overline d_n$ so that, in addition, $d_n\ge d_{n-1}+6$, which makes the disks $B_3(d_ne_1)$ disjoint. By \fref{preservation} with $C_0\sum_{j\ge 1}\eps_j\le 1$, for all $N\ge n$ and $(t,x)\in\Sigma_n$,
\be
\label{prefixbound}
|\nabla\omega^{(N)}(t,x+d_ne_1)|\ge|\nabla\omega_n(t,x)|-1.
\ee

\noindent{\bf Step 4:} The initial datum and the solution. We let
$$\omega_0=\sum_{n\ge 1}\omega_n(0,\cdot-d_ne_1),\qquad u_0=\sum_{n\ge 1}V_n(0,\cdot-d_ne_1).$$
For the first family, all the norms in \fref{componentbounds} are bounded by $C_me^{-j}$; together with \fref{summability}, the $L^1$ norms of the vorticities and the $L^2$ norms of the velocities are summable, and the $H^m$ norms of the vorticities are square summable. Since the supports of the vorticities are disjoint, the squared $H^m$ norm of a finite sum of them is the sum of their squared $H^m$ norms, for integer $m$. Hence the series defining $\omega_0$ converges in $L^1$ and in every $H^m$. For the velocities, the identity $\|\nabla v\|_{L^2}=\|\curl v\|_{L^2}$ for divergence-free fields gives $\|v\|_{H^m}\lesssim_m\|v\|_{L^2}+\|\curl v\|_{H^{m-1}}$, so that the series defining $u_0$ converges in every $H^m$. The partial sums satisfy $u_0^{(N)}=K*\omega_0^{(N)}$, and \fref{velocitybound} with $\rho=1$ shows that $K*\omega^{(N)}_0\to K*\omega_0$ uniformly, so that $u_0=K*\omega_0$. The datum $\omega_0$ is odd in $x_2$ and nonnegative on $\Rp$, and its support is not compact.

Let $u$ be the global smooth solution of \fref{eulervelocity} with initial velocity $u_0$, and $\omega=\curl u$. Since $u^{(N)}(0)\to u_0$ in $H^4$, the $H^4$ energy estimate of Step 2 of the proof of Lemma \ref{lemmadecoupling}, with $v_d$ replaced by $u$ and $R_d=0$, shows that $u^{(N)}\to u$ in $C([0,S];H^4)$ for every $S>0$. Since $H^4\subset C^2$, we may pass to the limit $N\to+\infty$ in \fref{prefixbound}, and we obtain
\be
\label{limitbound}
\|\nabla\omega(t)\|_{L^\infty}\ge|\nabla\omega_n(t,x)|-1\ \ \mbox{for all}\ \ n\ge 1\ \ \mbox{and}\ \ (t,x)\in\Sigma_n.
\ee
The solution $\omega(t)$ remains odd in $x_2$, nonnegative on $\Rp$, in $L^1$, and $u(t)=K*\omega(t)$ for all $t$.\\

\noindent{\bf Step 5:} Proof of \fref{growthsequence}. By \fref{limitbound} with $n=2j-1$, \fref{growthfirstfamily} and \fref{choicetauj},
$$\frac{\log\left(2+\|\nabla\omega(t_j)\|_{L^\infty}\right)}{\Gamma(t_j)}\ge\frac{\frac34\tau_j}{\tau_j/j^2}=\frac34j^2\to+\infty.$$
Since $\Gamma\ge 1$, this implies $\|\nabla\omega(t_j)\|_{L^\infty}e^{-c\Gamma(t_j)}\to+\infty$ for every $c>0$.\\

\noindent{\bf Step 6:} Proof of \fref{growthalltime}. Let $t\in J_j$. By \fref{limitbound} with $n=2j$ and \fref{growthsecondfamily},
\begin{align*}
\log\left(2+\|\nabla\omega(t)\|_{L^\infty}\right)&\ge a_jt+\log c_*+\log\frac{a_j}{r_j}-6\log(1+T_j)\\
&\ge a_jt-6T_{j-1}^{\frac12}-\log m_j-C,
\end{align*}
where we used $\log\frac{a_j}{r_j}=2\sqrt{\log m_j}-\log m_j$ and $\log(1+T_j)\le T_{j-1}^{\frac12}+1$. Since $t\ge S_{j-1}$, $a_jt\ge\frac{a_j}{a_{j-1}}T_{j-1}\ge\frac12T_{j-1}$ for $j$ large, and $T_{j-1}$ grows much faster than $T_{j-1}^{1/2}$ and $\log m_j$. Hence, for $j$ large and $t\in J_j$,
\be
\label{lowerJj}
\log\left(2+\|\nabla\omega(t)\|_{L^\infty}\right)\ge\frac{a_jt}{2}.
\ee
It remains to compare $j$ with $\log_*t$. We claim that $T_j\ge 9E_j^2$ for $j\ge 0$. This holds for $j=0$ by definition, and if it holds for $j$, then
$$T_{j+1}=e^{T_j^{1/2}}\ge e^{3E_j}\ge 9e^{2E_j}=9E^2_{j+1}.$$
Since $a_{j-1}\le 1$, we get $t\ge S_{j-1}\ge T_{j-1}\ge E_{j-1}$ for $t\in J_j$. The function $\log_*$ is nondecreasing, and $\log_*E_i=i+\log_*E_0$ since $\log E_{i}=E_{i-1}$ and $E_0>1$. Therefore $\log_*t\ge j-1$, and \fref{lowerJj} yields, for $t\in J_j$ with $j$ large,
$$\frac{\log_*t}{t}\log\left(2+\|\nabla\omega(t)\|_{L^\infty}\right)\ge\frac{(j-1)a_j}{2}=\frac{(j-1)e^{\sqrt{\log(j+j_0)}}}{2(j+j_0)}\to+\infty.$$
This proves \fref{growthalltime}, and the proof of Theorem \ref{thmmain} is complete. \hfill$\Box$

\section{Periodic solutions with double exponential growth on a finite time interval}
\label{sectionperiodic}

This section is devoted to the construction of smooth periodic solutions whose vorticity gradient grows double exponentially at every time of a prescribed interval $[0,T]$, and whose dilations have small derivatives at early times. The mechanism is the hyperbolic point of the Bahouri--Chemin cross \cite{BahouriChemin1994}, which is also at the origin of the constructions of Denisov \cite{Denisov2015} on the torus and of Kiselev and \v Sver\'ak \cite{KiselevSverak2014} in the disk.

\subsection{The Bahouri--Chemin cross}

\noindent{\it Biot--Savart law on tori}. For $L>0$, let $\T^2_L=(\R/L\Z)^2$, and $\T^2=\T^2_{2\pi}$. We identify functions on $\T^2_L$ with $L$-periodic functions on $\R^2$. For a function $f$ with zero mean on $\T^2_L$, $K_L*f$ denotes the unique divergence-free vector field with zero mean on $\T^2_L$ whose curl is $f$, where $K_L=\nabla^\perp G_L$ and $G_L$ is the Green function of the Laplacian on $\T^2_L$ with zero mean. We write $K_{\T}=K_{2\pi}$.

\begin{lemma}[Biot--Savart law on tori]
\label{lemmatorus}
The periodic Biot--Savart kernels have the following properties.
\begin{enumerate}
\item {\em Kernel}: for all integers $j\ge 0$ and all $L\ge 2\pi$,
$$|D^jK_L(z)|\lesssim_j|z|^{-1-j},\ \ z\in\left[-\frac L2,\frac L2\right]^2\setminus\{0\},$$
and $K_L-K$ is smooth on $B_{L/4}$.
\item {\em Velocity}: for $f$ with zero mean on $\T^2_L$, $L\ge 2\pi$,
\be
\label{velocitytorus}
\|K_L*f\|_{L^\infty}\lesssim\|f\|_{L^\infty}+\|f\|_{L^1(\T^2_L)},\qquad \|K_{\T}*f\|_{L^2(\T^2)}\le\|f\|_{L^2(\T^2)}.
\ee
\item {\em Interpolation on $\T^2$}: for $f\in C^1(\T^2)$ with zero mean,
\be
\label{interpolationtorus}
\|\nabla(K_{\T}*f)\|_{L^\infty}\lesssim\left(\|\nabla f\|_{L^\infty}\|f\|_{L^2(\T^2)}\right)^{\frac12}+\|f\|_{L^2(\T^2)}.
\ee
\end{enumerate}
\end{lemma}

The interpolation estimate \fref{interpolationtorus} is the key to control the strain created by a small and steep bump: its vorticity gradient may be huge, but if its $L^2$ norm is small enough, its velocity gradient remains small.

\begin{proof}
On $\T^2$, the Green function is $G_{2\pi}(z)=\frac{1}{2\pi}\log|z|+H(z)$ with $H$ smooth on $[-\pi,\pi]^2$, and hence $|D^jK_{2\pi}(z)|\lesssim_j|z|^{-1-j}$ on $[-\pi,\pi]^2\setminus\{0\}$. The scaling $K_L(z)=\frac{2\pi}{L}K_{2\pi}(\frac{2\pi z}{L})$ yields (1). The first bound in \fref{velocitytorus} follows from (1) by splitting the integral at $|z|=1$, and the second one from the Fourier series, since the multiplier of $K_{\T}*$ is bounded by $|m|^{-1}\le 1$ for $m\in\Z^2\setminus\{0\}$. For \fref{interpolationtorus}, we write $\nabla(K_{\T}*f)(x)$ as the principal value integral of $\nabla K_{\T}(z)f(x-z)$ plus a multiple of $f(x)$. For $0<r\le 1$, the part $|z|<r$ is equal to $\int_{|z|<r}\nabla K(z)(f(x-z)-f(x))dz+O(r\|f\|_{L^\infty})$, since $\nabla K$ has zero mean on circles, and it is bounded by $Cr\|\nabla f\|_{L^\infty}$ plus the local term; the part $|z|\ge r$ is bounded by $\|\nabla K_{\T}\|_{L^2(\{|z|\ge r\})}\|f\|_{L^2}\lesssim r^{-1}\|f\|_{L^2}$; and the local term satisfies $|f(x)|\le|\frac{1}{|B_r|}\int_{B_r(x)}f|+r\|\nabla f\|_{L^\infty}\lesssim r^{-1}\|f\|_{L^2}+r\|\nabla f\|_{L^\infty}$. The choice $r=\min\{1,(\|f\|_{L^2}/\|\nabla f\|_{L^\infty})^{\frac12}\}$ yields \fref{interpolationtorus}.
\end{proof}

\noindent{\it The cross and its regularization}. The Bahouri--Chemin cross is the stationary weak solution of \fref{euler} on $\T^2$ given by
\be
\label{defcross}
\Theta(x)=\sgn(\sin x_1)\sgn(\sin x_2),\qquad U_\Theta=K_{\T}*\Theta.
\ee
It is odd in $x_1$ and in $x_2$, equal to $1$ on the square $(0,\pi)^2$, and its velocity is tangent to the coordinate axes. We let
$$\kappa=\frac2\pi.$$
To regularize $\Theta$, we fix an even function $\varrho\in C^\infty_c(-1,1)$ with $\varrho\ge 0$ and $\int\varrho=1$, we let $s_\sigma=\varrho_\sigma*\sgn(\sin\cdot)$ with $\varrho_\sigma=\sigma^{-1}\varrho(\cdot/\sigma)$, and
\be
\label{defthetasigma}
\theta_\sigma(x)=s_\sigma(x_1)s_\sigma(x_2),\ \ 0<\sigma\le\frac12.
\ee
Then $\theta_\sigma$ is smooth, odd in $x_1$ and in $x_2$, with values in $[0,1]$ on $[0,\pi]^2$; it is equal to $\Theta$ outside the strips $\{\dist(x_1,\pi\Z)<\sigma\}\cup\{\dist(x_2,\pi\Z)<\sigma\}$, and
\be
\label{boundsthetasigma}
\|\theta_\sigma-\Theta\|_{L^1(\T^2)}\lesssim\sigma,\qquad\|D^j\theta_\sigma\|_{L^\infty}\lesssim_j\sigma^{-j}.
\ee

\begin{lemma}[Velocity of the cross]
\label{lemmacross}
There exist $r_0\in(0,\frac14]$ and $C_1\ge 1$ such that, on the region $\mathcal O=\{x:0<x_1<x_2^2,\ 0<x_2<r_0\}$, the following holds.
\begin{enumerate}
\item {\em Hyperbolic velocity}:
\be
\label{crossvalues}
\left|U_{\Theta,1}(x)-\kappa x_1\log\frac1{x_2}\right|\le C_1x_1,\qquad\left|U_{\Theta,2}(x)+\kappa x_2\log\frac1{x_2}\right|\le C_1x_2.
\ee
\item {\em Velocity gradient}:
\be
\label{crossgradient}
\left|\pa_1U_{\Theta,1}(x)-\kappa\log\frac1{x_2}\right|\le C_1,\qquad|\pa_2U_{\Theta,1}(x)|\le C_1x_2,\qquad|\pa_1U_{\Theta,2}(x)|\le C_1.
\ee
\item {\em One-sided smoothness}: for every $y_0\in(0,r_0)$ and $j\ge 1$, $|D^jU_\Theta|\le C_{j,y_0}$ on $\{0<x_1<r_0,\ y_0<x_2<r_0\}$.
\item {\em Strain of the regularized cross}: for $0<\sigma\le r_0$,
\be
\label{strainorigin}
\left|\pa_1(K_{\T}*\theta_\sigma)_1(0)-\kappa\log\frac1\sigma\right|\le C_1.
\ee
\end{enumerate}
\end{lemma}

Let us comment on these formulas. By \fref{crossvalues}, near the origin the flow of the cross is hyperbolic: it pushes the fluid away from the vertical axis and towards the horizontal axis, with the rate $\kappa\log\frac{1}{x_2}$, which is unbounded. A fluid particle at height $x_2(t)$ satisfies $\frac{d}{dt}\log\frac1{x_2}\approx\kappa\log\frac1{x_2}$, so that $\log\frac{1}{x_2(t)}$ grows exponentially, and the stretching rate grows exponentially as well. This is the origin of the double exponential growth. The restriction $x_1<x_2^2$ in $\mathcal O$ ensures that the off-diagonal entry $\pa_2U_{\Theta,1}$ of the velocity gradient is small, so that the flow is close to a pure stretching along the horizontal direction. The estimate \fref{strainorigin} shows that the strain of the regularized cross at the origin is of size $\log\frac1\sigma$, which is the source of the unbounded initial strain in Theorem \ref{thmdoubleexp}.

\begin{proof}
\noindent{\bf Step 1:} The velocity gradient near the origin. Let $\psi=G_{2\pi}*\Theta$ be the stream function, so that $U_\Theta=\nabla^\perp\psi$ and $\pa_1U_{\Theta,1}=-\pa_2G_{2\pi}*\pa_1\Theta$ in the sense of distributions. Near the origin, $\pa_1\Theta=2\sgn(\sin x_2)\delta_{\{x_1=0\}}$ plus singular measures on the lines $\{x_1=\pm\pi\}$, whose contributions are smooth near the origin. Using $G_{2\pi}=\frac{1}{2\pi}\log|z|+H$ and $\int_{-a}^a\frac{x_2-s}{x_1^2+(x_2-s)^2}\sgn(s)ds=\log(x_1^2+x_2^2)-\frac12\log\left((x_1^2+(x_2-a)^2)(x_1^2+(x_2+a)^2)\right)$ for $a=\frac12$, we obtain
$$\pa_1U_{\Theta,1}(x)=-\frac1\pi\log(x_1^2+x_2^2)+H_1(x),\qquad\pa_2U_{\Theta,2}=-\pa_1U_{\Theta,1},$$
where $H_1$ is smooth on $B_{4r_0}$ for $r_0$ small enough.\\

\noindent{\bf Step 2:} Proof of (1) and (2). Since $\Theta$ is odd in $x_1$ and in $x_2$, $U_{\Theta,1}(0,x_2)=U_{\Theta,2}(x_1,0)=0$, and Step 1 yields after integration
\be
\label{crossintegrated}
\left\{\begin{array}{l}
U_{\Theta,1}(x)=-\dfrac{1}{\pi}\displaystyle\int_0^{x_1}\log(s^2+x_2^2)ds+x_1R_1(x),\\[3mm]
U_{\Theta,2}(x)=\dfrac1\pi\displaystyle\int_0^{x_2}\log(x_1^2+s^2)ds+x_2R_2(x),
\end{array}\right.
\ee
with $R_1,R_2$ smooth on $B_{4r_0}$. On $\mathcal O$, $0\le s\le x_1<x_2^2$ gives $\log(s^2+x_2^2)=2\log x_2+O(x_2^2)$, which yields the first bound in \fref{crossvalues} and, with Step 1, the first bound in \fref{crossgradient}. Next, $\int_0^{x_2}\log(x_1^2+s^2)ds=2x_2\log x_2-2x_2+\int_0^{x_2}\log(1+\frac{x_1^2}{s^2})ds$ and $\int_0^{+\infty}\log(1+\frac{x_1^2}{s^2})ds=\pi x_1\le\pi x_2$, which yields the second bound in \fref{crossvalues}. Differentiating \fref{crossintegrated},
$$\pa_2U_{\Theta,1}=-\frac2\pi\arctan\frac{x_1}{x_2}+x_1\pa_2R_1,\qquad\pa_1U_{\Theta,2}=\frac2\pi\arctan\frac{x_2}{x_1}+x_2\pa_1R_2,$$
and $\arctan\frac{x_1}{x_2}\le\frac{x_1}{x_2}\le x_2$ on $\mathcal O$, which proves the last two bounds in \fref{crossgradient}.\\

\noindent{\bf Step 3:} Proof of (3). On the square $(0,\pi)^2$, $\Theta=1$, so that $\phi=\psi-\frac{x_1^2}{2}$ is harmonic there, and $\phi=0$ on $\{x_1=0\}$ since $\psi$ is odd in $x_1$. By the Schwarz reflection principle, the odd extension of $\phi$ in $x_1$ is harmonic on $(-\pi,\pi)\times(0,\pi)$ and bounded. The interior estimates for harmonic functions \cite[Theorem 2.10]{GilbargTrudinger2001} bound its derivatives on $\{|x_1|<r_0,\ y_0<x_2<r_0\}$ by constants depending on $y_0$, and $U_\Theta=\nabla^\perp\phi+(0,x_1)$ on $(0,\pi)^2$.\\

\noindent{\bf Step 4:} Proof of (4). Since mixed second derivatives of $\log|z|$ do not produce a Dirac mass, Lemma \ref{lemmatorus} (1) and $\pa_1\pa_2\log|z|=-\frac{2z_1z_2}{|z|^4}$ yield
$$\pa_1(K_{\T}*\theta_\sigma)_1(0)=-\pa_1\pa_2(G_{2\pi}*\theta_\sigma)(0)={\rm p.v.}\int_{|y|<a}\frac{y_1y_2}{\pi|y|^4}\theta_\sigma(y)dy+O(1).$$
On $\{|y|<2\sigma\}$, $|\theta_\sigma(y)|\lesssim\frac{|y_1||y_2|}{\sigma^2}$ because $s_\sigma$ is odd with $|s'_\sigma|\lesssim\sigma^{-1}$, so that this region contributes $O(1)$. On $\{2\sigma\le|y|<a\}$, the difference $\theta_\sigma-\Theta$ is supported in the strips $\{|y_1|<\sigma\}\cup\{|y_2|<\sigma\}$, where $\frac{|y_1y_2|}{|y|^4}\le\frac{\sigma}{|y|^3}$, and the area of the strips at distance $\sim\rho$ from the origin is $O(\sigma\rho)$, so that the difference contributes $O(1)$. Finally, since $y_1y_2\Theta(y)=|y_1y_2|$,
\begin{align*}
\int_{2\sigma\le|y|<a}\frac{|y_1y_2|}{\pi|y|^4}dy&=\frac{1}{\pi}\int_{2\sigma}^a\frac{d\rho}{\rho}\int_0^{2\pi}|\cos\phi\sin\phi|d\phi\\
&=\frac{2}{\pi}\log\frac{a}{2\sigma}=\kappa\log\frac1\sigma+O(1).
\end{align*}
\end{proof}

\subsection{Growth on a finite time interval}

We now perturb the regularized cross $\theta_\sigma$ by a small bump placed close to the vertical axis, at height $e^{-\ell_0}$. The bump is the passive part of the construction: it is transported and stretched by the flow of the cross, and its vorticity gradient grows double exponentially. We first introduce the parameters.

\medskip

\noindent{\it Constants}. Let $C_1$ and $r_0$ be given by Lemma \ref{lemmacross}. We fix
\be
\label{defC2}
C_2=\max\left\{1+\frac{C_1+2}{\kappa},\ 3C_1+6\right\},\qquad K_0=12.
\ee
The large constant $\ell_0\ge 7C_2$ is the logarithm of the inverse of the initial height of the bump introduced below, and $a_*\in(0,1]$ is a small constant measuring its amplitude; both are fixed in the proof of Proposition \ref{propperiodic} below, depending only on $C_1$ and $r_0$. We introduce the two stretching factors
\be
\label{deflambda}
\left\{\begin{array}{l}
\Lambda_+(t)=\exp\left((\ell_0+C_2)(e^{\kappa t}-1)+C_2t\right),\\[1mm]
\Lambda_-(t)=\frac12\exp\left((\ell_0-C_2)(e^{\kappa t}-1)-C_2t\right).
\end{array}\right.
\ee
The factor $\Lambda_+(t)$ is an upper bound for the deformation of the bump by the flow on $[0,t]$, and $\Lambda_-(t)$ is a lower bound for its horizontal stretching. The ratio of their logarithms is close to one, since $\ell_0\gg C_2$.

\medskip

\noindent{\it The bump}. For $T\ge 1$, we let
\be
\label{defAT}
A_T=a_*\Lambda_+(T)^{-\frac12},\qquad\eps_T=\exp\left(-K_0\ell_0e^{2\kappa T}\right),\qquad x_T=(2\eps_T,e^{-\ell_0}).
\ee
The height $A_T$ of the bump is chosen so that the strain created by the bump remains small on $[0,T]$, see \fref{strainbump} below. The distance $2\eps_T$ to the vertical axis is chosen so small that the bump remains in the region $\mathcal O$ of Lemma \ref{lemmacross} until time $T$. We fix a radial function $\eta_0\in C^\infty_c(B_1)$ with $0\le\eta_0\le 1$ and $\eta_0=1$ on $B_{3/5}$. For $0<h\le\frac{\eps_T}{2}$, we let $b_{T,h}$ be the function on $\T^2$ which is odd in $x_1$ and in $x_2$, and which is equal to
$$-A_T\,\eta_0\left(\frac{x-x_T}{h}\right)\ \ \mbox{on}\ \ (0,\pi)^2.$$
We consider the smooth solution $\theta$ of \fref{euler} on $\T^2$ with initial datum $\theta_\sigma+b_{T,h}$.

\begin{proposition}[Periodic solutions with double exponential growth on $\lbrack 0,T\rbrack$]
\label{propperiodic}
For suitable constants $\ell_0\ge 7C_2$ and $a_*\in(0,1]$, the following holds. For every $T\ge 1$ and every integer $k\in[1,T]$, there exists $\sigma_{T,k}\in(0,r_0]$ and, for every $\sigma\in(0,\sigma_{T,k}]$, there exists $h_{T,k,\sigma}\in(0,\frac{\eps_T}2]$ such that, for all $\sigma\in(0,\sigma_{T,k}]$ and $h\in(0,h_{T,k,\sigma}]$, the solution $\theta$ satisfies the following.
\begin{enumerate}
\item {\em Symmetry}: $\theta(t)$ is odd in $x_1$ and in $x_2$, and $|\theta(t,x)|\le 1$.
\item {\em Growth at every time}:
\be
\label{growthperiodic}
\|\nabla\theta(t)\|_{L^\infty}\ge\frac{A_T}{2h}\Lambda_-(t),\ \ 0\le t\le T.
\ee
\item {\em Early times}: there exist $a_k>0$ and $\overline C_k$ depending only on $k$, and $\overline C_{k,\sigma}$ depending only on $k$ and $\sigma$, such that, if $A_T\le a_k$,
\be
\label{earlytimes}
\sup_{0\le t\le k}\sum_{j=1}^kh^j\|D^j\theta(t)\|_{L^\infty}\le \overline C_kA_T+\overline C_{k,\sigma}h.
\ee
\item {\em Initial data}: the initial velocity $V_0=K_{\T}*(\theta_\sigma+b_{T,h})$ and vorticity satisfy
\be
\label{initialperiodic}
\pa_1V_{0,1}(0)\ge\kappa\log\frac1\sigma-C_2,\qquad \|V_0\|_{L^2(\T^2)}\ge c_V,
\ee
and $\|\theta_\sigma+b_{T,h}\|^p_{L^p(\T^2)}\ge 2\pi^2$ for $p=1,2$,
where $c_V=\frac12\|U_\Theta\|_{L^2(\T^2)}>0$.
\end{enumerate}
\end{proposition}

Note that $h^j\|D^j\theta(t)\|_{L^\infty}=\|D^j[\theta(t,h\,\cdot)]\|_{L^\infty}$: the left-hand side of \fref{earlytimes} measures the derivatives of the dilated solution $\theta(t,hx)$, which will be the actual component in the plane.

\begin{proof}
\noindent{\bf Step 1:} Decomposition of the solution. Since $\theta$ is odd in $x_1$ and in $x_2$, its velocity $u=K_{\T}*\theta$ is tangent to the lines $\{x_1\in\pi\Z\}$ and $\{x_2\in\pi\Z\}$, so that the flow of $u$ preserves each open square of the grid $\pi\Z^2$. The function $\Theta$, which is constant on these squares, is therefore transported by $u$. We decompose
\be
\label{decompositiontheta}
\theta=\Theta+p+b,\qquad\left\{\begin{array}{ll}\pa_tp+u\cdot\nabla p=0, & p(0)=\theta_\sigma-\Theta,\\ \pa_tb+u\cdot\nabla b=0, & b(0)=b_{T,h}.\end{array}\right.
\ee
The function $p$ is the error due to the regularization, and $b$ is the transported bump. Since $\theta_\sigma=\Theta=1$ on the disk $B_h(x_T)$ as soon as $\sigma<\eps_T$, the initial datum takes its values in $[0,1]$ on $(0,\pi)^2$, and so does $\theta(t)$ by transport. This proves (1). By transport, $\|p(t)\|_{L^1(\T^2)}=\|\theta_\sigma-\Theta\|_{L^1}\lesssim\sigma$ and $\|b(t)\|_{L^2(\T^2)}=\|b_{T,h}\|_{L^2}\le 2\sqrt\pi A_Th$. Moreover $\|u(t)\|_{L^\infty}\lesssim 1$ by \fref{velocitytorus}.\\

\noindent{\bf Step 2:} The regularization error. We claim that, for every integer $J\ge 1$, if $\sigma$ is small enough depending on $T$ and $J$, then on the region $\mathcal A_T=\{x:\frac{\eps_T}{2}<x_1<x_2^2,\ 0<x_2<3e^{-\ell_0}\}$,
\be
\label{boundp}
\sum_{j=0}^J|D^j(K_{\T}*p)(t,x)|\le\eps_T^3,\ \ x\in\mathcal A_T,\ \ 0\le t\le T.
\ee
Indeed, by the log-Lipschitz bound \fref{loglipschitz}, $|u(t,x)-u(t,y)|\le C_5|x-y|(1+\log_+\frac{1}{|x-y|})$ with $C_5$ independent of all the parameters. Let $X(t)$ be a trajectory starting in the strip $\{\dist(x_1,\pi\Z)<\sigma\}$, and let $\delta(t)\le 1$ be its distance to the closest line $\{x_1=m\pi\}$. The horizontal velocity vanishes on this line, so that $|\delta'|\le C_5\delta(1+\log\frac1\delta)$, and hence $1+\log\frac{1}{\delta(t)}\ge(1+\log\frac1\sigma)e^{-C_5t}$. The same holds for the horizontal strips. Therefore $p(t)$ is supported in the strips of width $e(\sigma/e)^{e^{-C_5T}}$ around the grid lines for $t\le T$, which is smaller than $\frac{\eps_T}{4}$ for $\sigma$ small. Every point of $\mathcal A_T$ has distance at least $\frac{\eps_T}{2}$ to the lines $\{x_1=0\}$ and $\{x_2=0\}$, and hence at least $\frac{\eps_T}{4}$ to the support of $p(t)$. Lemma \ref{lemmatorus} (1) yields $|D^j(K_{\T}*p)(t,x)|\lesssim_j\eps_T^{-1-j}\|p(t)\|_{L^1}\lesssim_j\eps_T^{-1-j}\sigma$ on $\mathcal A_T$, and \fref{boundp} follows for $\sigma$ small. We let $\sigma_{T,k}$ be such that \fref{boundp} holds with $J=k+1$ for $\sigma\le\sigma_{T,k}$, and such that $\sigma_{T,k}<\eps_T$.\\

\noindent{\bf Step 3:} Bootstrap. We let $\mathcal B_T=\{x:\eps_T<x_1<\frac{x_2^2}{8},\ 0<x_2<2e^{-\ell_0}\}\subset\mathcal A_T$, and we let $h\le\frac{e^{-\ell_0}}{2}$. The initial bump on $(0,\pi)^2$ is supported in the closed disk $\overline{B_h(x_T)}\subset\mathcal B_T$, and $\supp b(t)\cap(0,\pi)^2\subset X(t,\overline{B_h(x_T)})$, where $X$ is the flow of $u$. We claim that for $t\in[0,T]$:
\be
\label{bootstrapbump}
X(t,\overline{B_h(x_T)})\subset\mathcal B_T,\qquad\|\nabla(K_{\T}*b)(t)\|_{L^\infty}\le\frac12.
\ee
We assume that $X(s,\overline{B_h(x_T)})\subset\overline{\mathcal B_T}$ and $\|\nabla(K_{\T}*b)(s)\|_{L^\infty}\le 1$ for $s\in[0,t]$, $t\le T$, and we improve these bounds. In this regime, we have $\log\frac{1}{x_2}\ge \ell_0-1$ on $\mathcal B_T$, and the following estimates hold on $\overline{\mathcal B_T}$. Since $(K_{\T}*b)_i$ is odd in $x_i$ and even in the other variable for $i=1,2$, $|(K_{\T}*b)_i(x)|\le x_i$ for $i=1,2$; since $x_1,x_2>\frac{\eps_T}2$ on $\mathcal A_T$, \fref{boundp} yields $|(K_{\T}*p)_i(x)|\le 2\eps_T^2x_i\le x_i$. Together with Lemma \ref{lemmacross} and $u=U_\Theta+K_{\T}*p+K_{\T}*b$, we obtain
\be
\label{velocityboot}
\left|u_1-\kappa x_1\log\frac1{x_2}\right|\le(C_1+2)x_1,\qquad\left|u_2+\kappa x_2\log\frac1{x_2}\right|\le(C_1+2)x_2,
\ee
and, writing $\nabla u=(f_{ij})$ with $f_{ij}=\pa_ju_i$,
\be
\label{jacobianboot}
\left|f_{11}-\kappa\log\frac1{x_2}\right|\le C_1+2,\qquad f_{22}=-f_{11},\qquad|f_{12}|\le 3,\qquad|f_{21}|\le C_1+2,
\ee
where we used $C_1x_2\le 2C_1e^{-\ell_0}\le 1$ for $\ell_0$ large.

Let $X(t)=X(t,\alpha)$ be a trajectory with $\alpha\in\overline{B_h(x_T)}$, and let $\ell(t)=\log\frac{1}{X_2(t)}$. By \fref{velocityboot}, $|\ell'-\kappa\ell|\le C_1+2$, and $|\ell(0)-\ell_0|\le 1$ since $h\le\frac{e^{-\ell_0}}{2}$. Integrating these differential inequalities and using \fref{defC2}, we obtain
\be
\label{boundell}
(\ell_0-C_2)e^{\kappa t}\le\ell(t)\le(\ell_0+C_2)e^{\kappa t}.
\ee
We now improve the location of the support. By \fref{velocityboot}, $u_2<0$ on $\overline{\mathcal B_T}$ since $\kappa(\ell_0-1)>C_1+2$, so that $X_2(t)<X_2(0)<2e^{-\ell_0}$. At a point with $x_1=\eps_T$, $u_1\ge x_1(\kappa(\ell_0-1)-C_1-2)>0$, so that $X_1(t)>\eps_T$ since $X_1(0)>\eps_T$. Next, $X_1'\le X_1(\kappa\ell+C_1+2)\le X_1(\kappa(\ell_0+C_2)e^{\kappa t}+C_2)$, so that $X_1(t)\le 3\eps_T\Lambda_+(t)$. Since $X_2(t)^2\ge\exp(-2(\ell_0+C_2)e^{\kappa t})$ by \fref{boundell}, the bound $X_1(t)<\frac{X_2(t)^2}{8}$ holds if
$$-K_0\ell_0e^{2\kappa T}+3(\ell_0+C_2)e^{\kappa T}+C_2T+\log 24<0.$$
Using $C_2\le\frac{\ell_0}{7}$, $C_2T\le C_2\kappa^{-1}e^{2\kappa T}\le \ell_0e^{2\kappa T}$, $3(\ell_0+C_2)e^{\kappa T}\le 6\ell_0e^{2\kappa T}$ and $\log 24\le 4\ell_0e^{2\kappa T}$, the left-hand side is bounded by $(-K_0+6+1+4)\ell_0e^{2\kappa T}<0$ since $K_0=12$. Therefore $X(t,\overline{B_h(x_T)})\subset\mathcal B_T$.

We now improve the strain of the bump. By \fref{jacobianboot} and \fref{boundell}, the velocity gradient along $X(t)$ satisfies
\be
\label{normjacobian}
\|\nabla u(t,X(t))\|\le\kappa(\ell_0+C_2)e^{\kappa t}+C_2,
\ee
where we used that the matrix $\nabla u$ is the sum of ${\rm diag}(f_{11},-f_{11})$, whose norm is $|f_{11}|$, and of an off-diagonal matrix of norm at most $C_1+3$. Since
\begin{align*}
&\nabla b(t,X(t,\alpha))=(\nabla_\alpha X(t,\alpha))^{-T}\nabla b(0,\alpha),\\
&\pa_t(\nabla_\alpha X)^{-1}=-(\nabla_\alpha X)^{-1}\nabla u(t,X),
\end{align*}
the bound \fref{normjacobian} yields
\be
\label{gradientbump}
\|\nabla b(t)\|_{L^\infty}\le\|\nabla b_{T,h}\|_{L^\infty}\Lambda_+(t)\le\|\nabla\eta_0\|_{L^\infty}\frac{A_T}{h}\Lambda_+(t).
\ee
The interpolation estimate \fref{interpolationtorus}, $\|b(t)\|_{L^2}\le 2\sqrt\pi A_Th$ and the definition \fref{defAT} of $A_T$ yield
\be
\label{strainbump}
\|\nabla(K_{\T}*b)(t)\|_{L^\infty}\lesssim A_T\Lambda_+(t)^{\frac12}+A_Th\lesssim a_*\le\frac12
\ee
for $a_*$ small enough. Hence both bounds in \fref{bootstrapbump} are strictly improved, and a continuity argument proves \fref{bootstrapbump} on $[0,T]$. In particular, \fref{velocityboot}, \fref{jacobianboot}, \fref{boundell} and \fref{normjacobian} hold for all trajectories starting in $\overline{B_h(x_T)}$ and all $t\in[0,T]$.\\

\noindent{\bf Step 4:} Horizontal stretching. Let $\alpha\in\overline{B_h(x_T)}$, and let $\xi(t)=\pa_{\alpha_1}X(t,\alpha)$, which solves $\xi'=\nabla u(t,X(t,\alpha))\xi$, $\xi(0)=e_1$. We claim that
\be
\label{stretchingxi}
\xi_1(t)\ge\Lambda_-(t),\ \ 0\le t\le T.
\ee
Let $\xi_1(t)=e^{\int_0^tf_{11}}\zeta(t)$. Solving the equation $\xi_2'=f_{21}\xi_1+f_{22}\xi_2$ with $\xi_2(0)=0$, and inserting the result in $\xi_1'=f_{11}\xi_1+f_{12}\xi_2$, we obtain the exact identity
$$\zeta(t)=1+\int_0^tf_{21}(s)\zeta(s)\left(\int_s^tf_{12}(\tau)\exp\left(-\int_s^\tau(f_{11}-f_{22})\right)d\tau\right)ds.$$
By \fref{jacobianboot} and \fref{boundell}, $f_{11}-f_{22}=2f_{11}\ge 2\kappa(\ell_0-C_2)e^{\kappa s}-2C_2\ge\kappa \ell_0e^{\kappa s}$ on $[s,\tau]$ for $\ell_0$ large. Hence the integral operator acting on $\zeta$ has norm on $L^\infty(0,T)$ at most
$$\int_0^{+\infty}(C_1+2)\frac{3}{\kappa \ell_0e^{\kappa s}}ds=\frac{3(C_1+2)}{\kappa^2\ell_0}\le\frac14$$
for $\ell_0$ large. Thus $\|\zeta\|_{L^\infty}\le\frac43$ and $\zeta\ge 1-\frac14\cdot\frac43=\frac23$. Since $f_{11}\ge\kappa(\ell_0-C_2)e^{\kappa t}-C_2$, \fref{stretchingxi} follows. We now fix $\ell_0\ge 7C_2$ such that all the previous requirements hold, and then $a_*$ such that \fref{strainbump} holds.\\

\noindent{\bf Step 5:} Growth at every time. Let $0\le t\le T$. On the initial segment $\Gamma_0=\{x_T+(s,0):|s|\le\frac h2\}$, the initial vorticity is equal to $1-A_T$, and on the circle $\pa B_h(x_T)$ it is equal to $1$. By \fref{stretchingxi}, the map $s\mapsto X_1(t,x_T+(s,0))$ is increasing, so that the curve $X(t,\Gamma_0)$ is the graph of a function over an interval of length at least $h\Lambda_-(t)$. Let $d>0$ be the distance between $X(t,\Gamma_0)$ and $X(t,\pa B_h(x_T))$. For every point of the graph, the vertical segment of length $2d$ centered at this point does not meet $X(t,\pa B_h(x_T))$, and hence lies in the domain $X(t,B_h(x_T))$, whose area is $\pi h^2$. These segments are disjoint, so that $2d\,h\Lambda_-(t)\le\pi h^2$. Let $y\in X(t,\Gamma_0)$ and $y'\in X(t,\pa B_h(x_T))$ with $|y-y'|=d$. By transport, $\theta(t,y')-\theta(t,y)=A_T$, and the mean value theorem yields
$$\|\nabla\theta(t)\|_{L^\infty}\ge\frac{A_T}{d}\ge\frac{2}{\pi}\frac{A_T}{h}\Lambda_-(t),$$
which implies \fref{growthperiodic}. This lower bound holds at every time $t\in[0,T]$, for the same solution.\\

\noindent{\bf Step 6:} Early times: the bump. Let $k\in[1,T]$ be an integer. We now assume that $\sigma\le\sigma_{T,k}$, where $\sigma_{T,k}$ is chosen in Step 2 with $J=k+1$, and we study the dilated functions $f^h(x)=f(hx)$, which are periodic of period $\frac{2\pi}{h}$. The dilated bump $b^h(t,x)=b(t,hx)$ is transported by the dilated velocity $u^h(t,x)=h^{-1}u(t,hx)$, which does not change the time variable. We split $u^h=u_e^h+u_b^h$, where $u_e=U_\Theta+K_{\T}*p$ is the velocity of the regularized cross and $u_b=K_{\T}*b$ is the velocity of the bump, so that $u_b^h=K_{2\pi/h}*b^h$ and $D^ju^h_e(t,x)=h^{j-1}(D^ju_e)(t,hx)$.

Let $t\le k$. By \fref{boundell}, the support of $b(t)$ in $(0,\pi)^2$ lies in $\mathcal B_T\cap\{x_2\ge 2y_k\}$, where $2y_k=\exp(-(\ell_0+C_2)e^{\kappa k})$ depends only on $k$. On this set, Lemma \ref{lemmacross} (3) and \fref{boundp} yield $|D^ju_e|\le C_{j,y_k}+1$ for $1\le j\le k+1$, and hence, by symmetry and since $h\le 1$,
\be
\label{externalfield}
\sup\left\{|D^ju_e^h(t,x)|:x\in\supp b^h(t)\right\}\le C_k',\ \ 1\le j\le k+1,\ \ 0\le t\le k,
\ee
where $C'_k$ depends only on $k$. On the other hand, by \fref{velocitytorus} applied to $D^jb^h$, which has zero mean,
$$\|D^ju_b^h\|_{L^\infty}\lesssim\|D^jb^h\|_{L^\infty}+\|D^jb^h\|_{L^1(\T^2_{2\pi/h})}.$$
Let $Y(t)=\sum_{j=1}^{k+1}\left(\|D^jb^h(t)\|_{L^\infty}+\|D^jb^h(t)\|_{L^1(\T^2_{2\pi/h})}\right)$. Differentiating $j$ times the transport equation of $b^h$, we obtain $(\pa_t+u^h\cdot\nabla)D^jb^h=F_j$, where $F_j$ is a linear combination of contractions of $D^lu^h$ with $D^{j+1-l}b^h$, $1\le l\le j$. The function $F_j$ vanishes outside $\supp b^h(t)$, so that only the values of $D^lu^h$ on this support enter. By Lemma \ref{lemmatransport}, $\frac{d^+}{dt}\|D^jb^h\|_{L^\infty}\le\|F_j\|_{L^\infty}$ and $\frac{d^+}{dt}\|D^jb^h\|_{L^1}\le\|F_j\|_{L^1}$, and \fref{externalfield} yields
$$\frac{d^+}{dt}Y\le C''_k(1+Y)Y,\qquad Y(0)\le C_kA_T,$$
where the bound on $Y(0)$ follows from the explicit form of $b^h(0)$, which consists of four bumps of radius one and height $A_T$. Hence, if $A_T\le a_k$ with $a_k$ small depending only on $k$,
\be
\label{boundbumpdilated}
\sup_{0\le t\le k}\sum_{j=1}^{k+1}\left(\|D^jb^h(t)\|_{L^\infty}+\|D^ju^h_b(t)\|_{L^\infty}\right)\le C_kA_T\le 1.
\ee
In words: after dilation, the bump has size one and amplitude $A_T$; on the fixed time interval $[0,k]$, it is deformed by a flow whose derivatives are bounded on its support, and its derivatives remain of order $A_T$. Note that the large strain of the regularized cross near the axes does not enter, since $F_j$ only involves the velocity on the support of the bump.\\

\noindent{\bf Step 7:} Early times: the regularized cross. The function $\theta_c=\Theta+p=\theta-b$ is smooth, transported by $u$, with $\theta_c(0)=\theta_\sigma$ and $|\theta_c|\le 2$. Its velocity is $u_c=K_{\T}*\theta_c=u-u_b$. We first bound the strain $M(t)=\|\nabla u(t)\|_{L^\infty}$ on $[0,k]$ on the fixed torus. By \fref{bootstrapbump}, $\|\nabla u_b\|_{L^\infty}\le\frac12$. By the velocity bound \fref{velocitytorus}, $\|u_c\|_{L^\infty}\lesssim\|\theta_c\|_{L^\infty}\lesssim 1$, and by the logarithmic estimate \fref{logestimate} and the monotonicity of $s\mapsto s(1+\log_+(L/s))$, $\|\nabla u_c\|_{L^\infty}\lesssim 1+\log(1+\|\nabla\theta_c\|_{L^\infty})$. By \fref{gradientgronwall} and \fref{boundsthetasigma}, $\|\nabla\theta_c(t)\|_{L^\infty}\lesssim\sigma^{-1}\exp(\int_0^tM)$. Hence $M(t)\le C(1+\log\frac1\sigma)+C\int_0^tM$, and Gronwall's lemma yields
\be
\label{strainbody}
\sup_{0\le t\le k}M(t)\le M_{k,\sigma}:=C\left(1+\log\frac1\sigma\right)e^{Ck}.
\ee
This bound depends on $\sigma$, but not on $h$, and it is invariant under the dilation: $\|\nabla u^h(t)\|_{L^\infty}=M(t)$.

We now estimate the dilated function $\theta_c^h(t,x)=\theta_c(t,hx)$, which is transported by $u^h=u_c^h+u_b^h$, with $u_c^h=K_{2\pi/h}*\theta_c^h$. Initially, $\|D^{j}\theta_c^h(0)\|_{C^{0,\gamma}}\lesssim_{j,\sigma}h$ for $j\ge 1$ and $h\le 1$ by \fref{boundsthetasigma}, where we fix $\gamma=\frac12$. Let $Z_m(t)=\|D^{m+1}\theta_c^h(t)\|_{C^{0,\gamma}}$ for $0\le m\le k-1$. The second Schauder estimate \fref{schauder}, which applies on the torus of period $\frac{2\pi}h\ge 4$ with uniform constants, gives
$$\|\nabla u_c^h\|_{C^{m,\gamma}}\lesssim_m\|\nabla u_c^h\|_{L^\infty}+\|\theta_c^h\|_{C^{m,\gamma}}\lesssim_m M_{k,\sigma}+1+\sum_{i<m}Z_i+Z_0\1_{m=0},$$
where we used $|\theta_c^h|\le 2$ and, for $m=0$, $\|\theta_c^h\|_{C^{0,\gamma}}\lesssim\|\theta_c^h\|_{L^\infty}+\|\nabla\theta_c^h\|_{L^\infty}$,
and \fref{boundbumpdilated} gives $\|\nabla u_b^h\|_{C^{m,\gamma}}\le C_k$ for $m\le k-1$. The function $D^{m+1}\theta_c^h$ solves a transport equation by $u^h$ whose source term is a linear combination of contractions of $D^{l+1}u^h$ with $D^{m+1-l}\theta_c^h$, $0\le l\le m$; each of these terms contains a derivative of $\theta_c^h$ of positive order. By the H\"older transport estimate \fref{transportholder}, whose exponential factor is bounded by $e^{\gamma kM_{k,\sigma}}$ by \fref{strainbody}, we obtain for $t\le k$
\begin{align*}
Z_m(t)\le C_{k,\sigma}\Big(Z_m(0)&+\int_0^t(1+Z_0)Z_m\,ds\\
&+\1_{m\ge 1}\int_0^t\Big(1+\sum_{i<m}Z_i\Big)\sum_{i<m}Z_i\,ds\Big),
\end{align*}
where the quadratic term $Z_0Z_m$ comes from the term $\nabla u^h_c\cdot D^{m+1}\theta^h_c$ for $m=0$. For $m=0$, this inequality and $Z_0(0)\lesssim_\sigma h$ yield $Z_0\le C_{k,\sigma}h$ on $[0,k]$ for $h$ small enough, and by induction on $m$, the inequality is linear in $Z_m$ with small lower order terms, so that $Z_m\le C_{k,\sigma}h$ on $[0,k]$ for $m\le k-1$. Together with \fref{boundbumpdilated} and $h^j\|D^j\theta(t)\|_{L^\infty}\le\|D^j\theta^h_c(t)\|_{L^\infty}+\|D^jb^h(t)\|_{L^\infty}$, this proves \fref{earlytimes}.\\

\noindent{\bf Step 8:} Initial data. By Lemma \ref{lemmacross} (4) and \fref{bootstrapbump} at $t=0$, $\pa_1V_{0,1}(0)\ge\kappa\log\frac1\sigma-C_1-\frac12$. Next, by \fref{velocitytorus}, $\|V_0\|_{L^2}\ge\|U_\Theta\|_{L^2}-\|\theta_\sigma-\Theta\|_{L^2}-\|b_{T,h}\|_{L^2}\ge 2c_V-C\sigma^{\frac12}-CA_Th$, where $c_V>0$ since $\curl U_\Theta=\Theta\neq 0$. Similarly, $\|\theta_\sigma+b_{T,h}\|_{L^1}\ge 4\pi^2-C\sigma-CA_Th^2$ and $\|\theta_\sigma+b_{T,h}\|_{L^2}\ge 2\pi-C\sigma^{\frac12}-CA_Th$. Decreasing $\sigma_{T,k}$ and $h_{T,k,\sigma}$ yields \fref{initialperiodic}.
\end{proof}

\subsection{A sequence of periodic components}

We now choose the parameters of Proposition \ref{propperiodic} along a sequence, so that the dilated solutions are active on consecutive time intervals of length $\frac12$, and have summable derivatives at early times.

\begin{proposition}[The sequence of periodic components]
\label{propsequence}
There exist an integer $n_0\ge 1$, a constant $c_1>0$, integers $k_n\to+\infty$, and sequences $\sigma_n\to 0$ and $h_n\in(0,2^{-n}]$ such that the following holds. Let $\tau_n=\frac{n+n_0}{2}$, let $\theta_n$ be the solution of Proposition \ref{propperiodic} with $T=\tau_n$, $k=k_n$, $\sigma=\sigma_n$ and $h=h_n$, and let
$$\omega_n^{\rm per}(t,x)=\theta_n(t,h_nx),\qquad V_n^{\rm per}(t,x)=h_n^{-1}(K_{\T}*\theta_n)(t,h_nx),$$
which is a smooth solution of \fref{euler} on the torus of period $\frac{2\pi}{h_n}$. Then for all $n\ge 1$:
\begin{enumerate}
\item {\em Symmetry}: $\omega_n^{\rm per}$ is odd in $x_1$ and in $x_2$, and $|\omega_n^{\rm per}|\le 1$.
\item {\em Growth on the $n$-th time interval}:
\be
\label{growthwindow}
\|\nabla\omega_n^{\rm per}(t)\|_{L^\infty}\ge 8\exp\left(c_1e^{\kappa t}\right),\ \ t\in[\tau_{n-1},\tau_n].
\ee
\item {\em Small derivatives at early times}: $k_n\le\tau_n$ and
\be
\label{earlysequence}
\sup_{0\le t\le k_n}\sum_{j=1}^{k_n}\|D^j\omega_n^{\rm per}(t)\|_{L^\infty}\le 2^{-n-4}.
\ee
\item {\em Initial data}: on the cell $\mathcal C_n=[-\frac{\pi}{h_n},\frac{\pi}{h_n}]^2$,
\be
\label{initialsequence}
\pa_1V^{\rm per}_{n,1}(0,0)\ge\kappa\log\frac{1}{\sigma_n}-C_2,\ \ \ \int_{\mathcal C_n}|V_n^{\rm per}(0)|^2\ge\frac{c_V^2}{h_n^4},
\ee
and $\int_{\mathcal C_n}|\omega_n^{\rm per}(0)|^p\ge 2\pi^2h_n^{-2}$ for $p=1,2$.
\end{enumerate}
\end{proposition}

\begin{proof}
\noindent{\bf Step 1:} Growth on consecutive intervals. The dilation multiplies the gradient by $h_n$, so that \fref{growthperiodic} yields $\|\nabla\omega_n^{\rm per}(t)\|_{L^\infty}\ge\frac{A_{\tau_n}}{2}\Lambda_-(t)$ for $0\le t\le\tau_n$. By \fref{deflambda} and \fref{defAT}, for $t\in[\tau_n-\frac12,\tau_n]$, using $e^{\kappa\tau_n}\le e^{\frac\kappa2}e^{\kappa t}$,
$$\log\left(\frac{A_{\tau_n}}{2}\Lambda_-(t)\right)\ge\left[(\ell_0-C_2)-\frac{e^{\kappa/2}}{2}(\ell_0+C_2)\right]e^{\kappa t}-C_2(t+\tau_n)-C.$$
Since $e^{\kappa/2}=e^{1/\pi}\le 1.375$ and $\ell_0\ge 7C_2$, the bracket is at least $0.3125\ell_0-1.6875C_2\ge\frac{C_2}{2}$. Hence, with $c_1=\frac{C_2}{4}$ and $n_0$ large, \fref{growthwindow} holds for all $n\ge 1$. This is where the choice $\ell_0\ge 7C_2$ is used: the cost $A_{\tau_n}$ of the small amplitude is a double exponential at the end time $\tau_n$, which is compensated by the stretching $\Lambda_-(t)$ at times $t\ge\tau_n-\frac12$.\\

\noindent{\bf Step 2:} Diagonal choice of the early times. Since $A_{\tau_n}\le a_*\exp(-\frac{\ell_0}{2}(e^{\kappa\tau_n}-1))$ decays double exponentially in $n$, for every integer $k\ge 1$ there exists $N_k$ such that $\tau_n\ge k$, $A_{\tau_n}\le a_k$ and $\overline C_kA_{\tau_n}\le 2^{-n-5}$ for $n\ge N_k$, where $a_k$ and $\overline C_k$ are given by Proposition \ref{propperiodic}. Increasing $n_0$, we may assume $N_1=1$, and we may choose $N_k$ strictly increasing. We let $k_n=\max\{k\ge 1:N_k\le n\}$, so that $k_n\to+\infty$ and $k_n\le\tau_n$. We then choose $\sigma_n=\min\{\sigma_{\tau_n,k_n},2^{-n}\}$, and $h_n\le\min\{h_{\tau_n,k_n,\sigma_n},2^{-n}\}$ such that $\overline C_{k_n,\sigma_n}h_n\le 2^{-n-5}$. Then \fref{earlytimes} yields \fref{earlysequence}.\\

\noindent{\bf Step 3:} Initial data. Since $V^{\rm per}_n(0,x)=h_n^{-1}V_0(h_nx)$, where $V_0$ is the initial velocity of $\theta_n$, the strain at the origin is unchanged by the dilation, and $\int_{\mathcal C_n}|V_n^{\rm per}(0)|^2=h_n^{-4}\|V_0\|^2_{L^2(\T^2)}$, $\int_{\mathcal C_n}|\omega_n^{\rm per}(0)|^p=h_n^{-2}\|\theta_n(0)\|^p_{L^p(\T^2)}$. Hence \fref{initialperiodic} yields \fref{initialsequence}.
\end{proof}

\section{Double exponential growth in $\R^2$: proof of Theorem \ref{thmdoubleexp}}
\label{sectionplane}

In this section we transport the periodic solutions of Proposition \ref{propsequence} to the plane, and we assemble them. Since the periodic solutions have infinite energy, we cut off their stream function at a large scale, and we need to show that the resulting solutions in the plane remain close to the periodic ones on compact sets, while their derivatives remain controlled everywhere. The natural framework for this limit is the theory of Serfati solutions.

\subsection{Serfati solutions and compactification}

Following Serfati \cite{Serfati1995} and Ambrose, Kelliher, Lopes Filho and Nussenzveig Lopes \cite{AKLNL2015}, we fix a radial function $a\in C^\infty_c(\R^2)$ equal to one near the origin, and we consider the Serfati identity
\be
\label{serfatiidentity}
u_j(t)-u_j(0)=(aK_j)*\left(\omega(t)-\omega(0)\right)-\int_0^t\nabla\nabla^\perp\left[(1-a)K_j\right]*\cdot(u\otimes u)(s)ds,
\ee
for $j=1,2$,
where $K_j$ is the $j$-th component of $K$, and $(\nabla\nabla^\perp g)*\cdot M=\sum_{l,m}(\pa_l\pa^\perp_mg)*M_{lm}$ with $\pa^\perp=(-\pa_2,\pa_1)$. The identity \fref{serfatiidentity} is obtained by splitting the Biot--Savart kernel into a local part $aK$ and a nonlocal part $(1-a)K$, and by writing the time derivative of the nonlocal part with the velocity equation; the nonlocal kernel $\nabla\nabla^\perp[(1-a)K_j]$ is $O(|z|^{-3})$ at infinity, and hence integrable. It makes sense for bounded velocities and vorticities, and it fixes the behavior of the velocity at infinity: for instance, the solution $u(t,x)=(t,0)$ of \fref{eulervelocity} with $p=-x_1$ does not satisfy it, see \cite[Section 7]{AKLNL2015}. We recall the following results, which are proved in \cite{AKLNL2015} both for the whole plane and for exterior domains; we only use the case of the whole plane.

\begin{lemma}[Serfati solutions \cite{AKLNL2015}]
\label{lemmaserfati}
Let $S>0$. A Serfati solution on $[0,S]$ is a vector field $u\in L^\infty((0,S)\times\R^2)\cap C([0,S]\times\R^2)$, divergence free, whose vorticity $\omega=\curl u$ is bounded, which solves the vorticity equation $\pa_t\omega+u\cdot\nabla\omega=0$ in the sense of distributions, and which satisfies \fref{serfatiidentity} for every radial $a\in C^\infty_c(\R^2)$ equal to one near the origin.
\begin{enumerate}
\item {\em Uniqueness} \cite[Theorem 2.5]{AKLNL2015}: two Serfati solutions on $[0,S]$ with the same initial velocity coincide.
\item {\em Smooth solutions} \cite[Proposition 4.1]{AKLNL2015}: a smooth solution of \fref{eulervelocity} on $[0,S]$ with compactly supported initial vorticity, whose velocity is given by the Biot--Savart law, satisfies \fref{serfatiidentity}.
\item {\em Uniform bound} \cite[Proposition 3.1 and Section 4.2]{AKLNL2015}: with $a_\lambda=a(\cdot/\lambda)$, $\|a_\lambda K\|_{L^1}\lesssim\lambda$ and $\|\nabla\nabla^\perp[(1-a_\lambda)K]\|_{L^1}\lesssim\lambda^{-1}$. As a consequence, a smooth solution $u$ satisfying \fref{serfatiidentity} for all $a_\lambda$, $\lambda>0$, satisfies for $0\le t\le S$
\be
\label{serfatibound}
\|u(t)\|^2_{L^\infty}\lesssim\|u(0)\|^2_{L^\infty}+\|\omega(0)\|_{L^\infty}\int_0^t\|u(s)\|_{L^\infty}^2ds.
\ee
\end{enumerate}
\end{lemma}

The bound \fref{serfatibound} follows as in \cite[Section 4.2, Step 2]{AKLNL2015}: since $\|\omega(t)\|_{L^\infty}=\|\omega(0)\|_{L^\infty}$, the identity \fref{serfatiidentity} with the cutoff $a_\lambda$ gives $\|u(t)\|_{L^\infty}\le\|u(0)\|_{L^\infty}+C\lambda\|\omega(0)\|_{L^\infty}+C\lambda^{-1}\int_0^t\|u\|^2_{L^\infty}$, and we choose $\lambda=(\|\omega(0)\|^{-1}_{L^\infty}\int_0^t\|u\|_{L^\infty}^2)^{\frac12}$. Its main feature is that it involves only $L^\infty$ norms, and no integrability of the vorticity. We now state the compactification result.

\begin{proposition}[From the torus to the plane]
\label{propcompactification}
Let $L\ge 4$, and let $\omega^{\rm per}$ be a smooth solution of \fref{euler} on $\T^2_L$, with velocity $V^{\rm per}$, lifted to $\R^2$, such that $\omega^{\rm per}$ and $V^{\rm per}$ have zero mean. Let $\psi^{\rm per}$ be the periodic stream function, so that $V^{\rm per}=\nabla^\perp\psi^{\rm per}$. Let $\chi\in C^\infty_c(B_2)$ be a radial function with values in $[0,1]$, equal to one on $B_1$. For $R\ge 1$, let $(\omega_R,V_R)$ be the global smooth solution of \fref{euler} on $\R^2$ with initial velocity
$$V_R(0)=\nabla^\perp\psi_R,\qquad\psi_R=\chi\left(\frac{\cdot}{R}\right)\psi^{\rm per}(0)\in C^\infty_c(\R^2).$$
Then, for every $S>0$, the following holds.
\begin{enumerate}
\item {\em Local convergence}: for every compact set $\mathcal K\subset\R^2$ and every $m\ge 0$,
$$\lim_{R\to+\infty}\sup_{0\le t\le S}\left(\|(V_R-V^{\rm per})(t)\|_{C^m(\mathcal K)}+\|(\omega_R-\omega^{\rm per})(t)\|_{C^m(\mathcal K)}\right)=0.$$
\item {\em Global control of the derivatives}: for every $m\ge 1$,
\be
\label{globalderivatives}
\limsup_{R\to+\infty}\sup_{0\le t\le S}\|D^m\omega_R(t)\|_{L^\infty(\R^2)}\le\sup_{0\le t\le S}\|D^m\omega^{\rm per}(t)\|_{L^\infty}.
\ee
\item {\em Initial data}: $\omega_R(0)=\Delta\psi_R$, $\|\omega_R(0)\|_{L^\infty}\le\|\omega^{\rm per}(0)\|_{L^\infty}+\frac CR$, and $(V_R(0),\omega_R(0))=(V^{\rm per}(0),\omega^{\rm per}(0))$ on $B_R$.
\end{enumerate}
\end{proposition}

The point of \fref{globalderivatives} is that the derivatives of $\omega_R$ are controlled not only on compact sets, where $\omega_R$ is close to $\omega^{\rm per}$, but also in the transition region $R\le|x|\le 2R$ of the cutoff and beyond. The proof shows that, seen from a point of this region, the solution looks like the periodic solution with a time rescaled by the factor $\chi(x/R)\in[0,1]$.

\begin{proof}
\noindent{\bf Step 1:} Uniform bounds. The constants below depend on the periodic solution and on $S$, but not on $R\ge 1$. The initial vorticity is
$$\omega_R(0)=\chi\left(\frac\cdot R\right)\omega^{\rm per}(0)+\frac2R(\nabla\chi)\left(\frac\cdot R\right)\cdot\nabla\psi^{\rm per}(0)+\frac{1}{R^2}(\Delta\chi)\left(\frac\cdot R\right)\psi^{\rm per}(0),$$
which yields (3) and bounds on all the derivatives of $V_R(0)$ and $\omega_R(0)$. Since $\psi_R$ is smooth and compactly supported, $V_R(0)=K*\omega_R(0)$ belongs to all the Sobolev spaces, and $(\omega_R,V_R)$ is a global smooth solution with compactly supported vorticity. By Lemma \ref{lemmaserfati} (2), it satisfies \fref{serfatiidentity}, and \fref{serfatibound} and Gronwall's lemma give $\|V_R(t)\|_{L^\infty}\le C$ on $[0,S]$, while $\|\omega_R(t)\|_{L^\infty}=\|\omega_R(0)\|_{L^\infty}\le C$. Next, the logarithmic estimate \fref{logestimate} and \fref{gradientgronwall} yield
$$\frac{d^+}{dt}\log\left(e+\|\nabla\omega_R\|_{L^\infty}\right)\le\|\nabla V_R\|_{L^\infty}\le C\left(1+\log\left(e+\|\nabla\omega_R\|_{L^\infty}\right)\right),$$
so that $\|\nabla\omega_R\|_{L^\infty}+\|\nabla V_R\|_{L^\infty}\le C$ on $[0,S]$. The first Schauder estimate \fref{schauder}, the H\"older transport estimate \fref{transportholder} applied to the equations of the derivatives of $\omega_R$, and an induction on the order of the derivatives, yield bounds on $\|\omega_R(t)\|_{C^{m,\gamma}}$ and $\|V_R(t)\|_{C^{m+1,\gamma}}$ on $[0,S]$ for all $m$: at each order, the equation of $D^m\omega_R$ is linear in $D^m\omega_R$, with coefficients involving $D^{l}V_R$, $l\le m$, which are controlled by the previous step since the velocity gains one derivative. Finally, $\pa_t\omega_R=-V_R\cdot\nabla\omega_R$, and the time derivative of \fref{serfatiidentity} expresses $\pa_tV_R$ as $(aK)*\pa_t\omega_R-\nabla\nabla^\perp[(1-a)K]*\cdot(V_R\otimes V_R)$, two convolutions with integrable kernels. By induction, all the space-time derivatives of $V_R$ and $\omega_R$ are bounded on $[0,S]\times\R^2$, uniformly in $R$.\\

\noindent{\bf Step 2:} The periodic solution is a Serfati solution. Let $\beta\in[0,1]$. The function $V_\beta(t,x)=\beta V^{\rm per}(\beta t,x)$ is a smooth periodic solution of \fref{eulervelocity} with vorticity $\omega_\beta(t,x)=\beta\omega^{\rm per}(\beta t,x)$. We claim that it satisfies \fref{serfatiidentity}. Both sides vanish at $t=0$, and we compare their time derivatives on each Fourier mode $\xi\in\frac{2\pi}{L}\Z^2$; recall that the convolution of a periodic function with an integrable kernel $g$ acts on the mode $\xi$ by multiplication by $\widehat g(\xi)=\int g(z)e^{-iz\cdot\xi}dz$. On the mode $\xi=0$, the left-hand side has zero time derivative because $V_\beta$ has zero mean; on the right-hand side, $\pa_t\omega_\beta$ has zero mean, and $\int\pa_l\pa_m^\perp[(1-a)K_j]=0$ because the first derivatives of $(1-a)K_j$ are $O(|z|^{-2})$ at infinity. On a mode $\xi\neq 0$, we write $\pa_t\omega_\beta=-\curl\dv(V_\beta\otimes V_\beta)$ and $\widehat{K}(\xi)=-i\xi^\perp|\xi|^{-2}=\widehat{aK}(\xi)+\widehat{(1-a)K}(\xi)$, where the last term is a smooth function on $\R^2\setminus\{0\}$. Since the periodic Biot--Savart law acts on the mode $\xi$ by $\widehat K(\xi)$, the time derivative of $V_\beta$ is $\widehat{aK}(\xi)$ applied to $\pa_t\omega_\beta$, plus $\widehat{(1-a)K}(\xi)$ applied to $-\curl\dv(V_\beta\otimes V_\beta)$, and transferring the two derivatives of $\curl\dv$ to the kernel yields the time derivative of the last term of \fref{serfatiidentity}. Hence $V_\beta$ is a Serfati solution. See also \cite[Section 7]{AKLNL2015} and \cite{Kelliher2015} for periodic Serfati velocities.\\

\noindent{\bf Step 3:} Translated limits. Let $R_j\to+\infty$ and $y_j\in L\Z^2$. We claim that, up to a subsequence, $\chi(y_j/R_j)\to\beta\in[0,1]$ and
\be
\label{translatedlimit}
(V_{R_j},\omega_{R_j})(\cdot,\cdot+y_j)\to(V_\beta,\omega_\beta)\ \ \mbox{in}\ \ C^m([0,S]\times\mathcal K)
\ee
for all $m\ge 0$ and all compact sets $\mathcal K\subset\R^2$.
By Step 1 and the Arzel\`a--Ascoli theorem, a subsequence converges in these spaces to a smooth bounded pair $(V,\omega)$ with $\dv V=0$, $\curl V=\omega$ and $\pa_t\omega+V\cdot\nabla\omega=0$. The translated solutions satisfy \fref{serfatiidentity}, which is translation invariant. We pass to the limit in \fref{serfatiidentity}: the local term converges because $aK$ is integrable with compact support, and the nonlocal term converges by dominated convergence, since $\nabla\nabla^\perp[(1-a)K]$ is integrable and the velocities are uniformly bounded and converge locally uniformly. Hence $V$ is a Serfati solution. Its initial velocity is the local limit of
$$V_{R_j}(0,x+y_j)=\chi\left(\frac{x+y_j}{R_j}\right)V^{\rm per}(0,x)+\frac{1}{R_j}(\nabla^\perp\chi)\left(\frac{x+y_j}{R_j}\right)\psi^{\rm per}(0,x),$$
where we used the periodicity of $\psi^{\rm per}$, that is $\beta V^{\rm per}(0)$. By Step 2 and the uniqueness of Lemma \ref{lemmaserfati} (1), $V=V_\beta$. The limit being unique, \fref{translatedlimit} holds along the whole subsequence for which $\chi(y_j/R_j)\to\beta$.\\

\noindent{\bf Step 4:} Conclusion. With $y_j=0$, we have $\beta=1$, and \fref{translatedlimit} for every sequence $R_j\to+\infty$ proves (1). Assume now that \fref{globalderivatives} fails for some $m\ge 1$: there exist $\eta>0$, $R_j\to+\infty$, $t_j\in[0,S]$ and $x_j\in\R^2$ with $|D^m\omega_{R_j}(t_j,x_j)|\ge\sup_{t\le S}\|D^m\omega^{\rm per}(t)\|_{L^\infty}+\eta$. We choose $y_j\in L\Z^2$ such that $x_j-y_j\in[0,L]^2$, and we extract subsequences such that $t_j\to t_*$, $x_j-y_j\to x_*$ and $\chi(y_j/R_j)\to\beta$. By \fref{translatedlimit} on $\mathcal K=[0,L]^2$,
$$D^m\omega_{R_j}(t_j,x_j)\to\beta\,D^m\omega^{\rm per}(\beta t_*,x_*),$$
whose absolute value is at most $\sup_{t\le S}\|D^m\omega^{\rm per}(t)\|_{L^\infty}$ since $\beta\le 1$ and $\beta t_*\le S$. This is a contradiction.
\end{proof}

\subsection{Assembly of nondecaying components}

We now place infinitely many compactly supported components far apart from each other. Contrary to Section \ref{sectionfiniteenergy}, the sum of their initial velocities has infinite energy, and the sum of their vorticities is not integrable. The distances must therefore be chosen so as to ensure, in addition to the preservation of the lower bounds, the convergence of the Biot--Savart integral and the smoothness of the limit in space and time.

\begin{proposition}[Assembly of nondecaying components]
\label{propassemblyinfinite}
Let $n_0\ge 1$, $c_1>0$, $\tau_n=\frac{n+n_0}{2}$, and let $k_n\to+\infty$ be integers with $k_n\le\tau_n$. For $n\ge 1$, let $(\omega_n,V_n)$ be a global smooth solution of \fref{euler} on $\R^2$ with initial velocity $V_n(0)=\nabla^\perp\psi_n$, where $\psi_n\in C^\infty_c(\R^2)$ is supported in a disk $B_{\rho_n}$, $\rho_n\ge 1$, such that:
\be
\label{hypassembly}
\left\{\begin{array}{ll}
\|\omega_n(0)\|_{L^\infty}\le 2,&\\[1mm]
\|D^j\omega_n(t)\|_{L^\infty}\le 2^{-n}, & 0\le t\le k_n,\ \ 1\le j\le k_n,\\[1mm]
\|\nabla\omega_n(t)\|_{L^\infty}\ge 2\exp\left(c_1e^{\kappa t}\right), & \tau_{n-1}\le t\le\tau_n.
\end{array}\right.
\ee
Then there exist distances $d_n\to+\infty$ such that the disks $B_{\rho_n}(d_ne_1)$ are pairwise disjoint, and a solution $(\omega,u,p)\in C^\infty([0,+\infty)\times\R^2)$ of \fref{euler} and \fref{eulervelocity} with initial data
$$\omega_0=\sum_{n\ge 1}\omega_n(0,\cdot-d_ne_1),\qquad u_0=\sum_{n\ge 1}V_n(0,\cdot-d_ne_1),$$
such that the following holds.
\begin{enumerate}
\item {\em Biot--Savart law}: for every $(t,x)$, $\int_{\R^2}|K(x-y)||\omega(t,y)|dy<+\infty$ and $u(t,x)=\int_{\R^2}K(x-y)\omega(t,y)dy$. Moreover \fref{loggrowthvelocity} holds.
\item {\em Regularity}: $\|\omega(t)\|_{L^\infty}\le 2$, and $\sup_{0\le t\le S}\|D^j\omega(t)\|_{L^\infty}<+\infty$ for all $S>0$, $j\ge 1$.
\item {\em Growth}: $\|\nabla\omega(t)\|_{L^\infty}\ge\exp(c_1e^{\kappa t})$ for all $t\ge\tau_0$.
\end{enumerate}
\end{proposition}

\begin{proof}
We denote by $M_n=\|\omega_n(0)\|_{L^1}$ the mass of the $n$-th component.\\

\noindent{\bf Step 1:} Choice of the distances. We apply Lemma \ref{lemmaassembly} with the horizons $S_n=\tau_n$, the Sobolev orders $s_n=n+5$, the errors $\eps_n=2^{-n-12}C_0^{-1}$, and the compact sets
$$\Sigma_n=\{(t,x):\tau_{n-1}\le t\le\tau_n,\ x\in\supp\omega_n(t)\}.$$
Each $V_n$ belongs to $C([0,S];H^m)$ for all $S$ and $m$, since $\psi_n$ is smooth and compactly supported. In addition to $d_n\ge\overline d_n$, we impose the following conditions, which hold for $d_n$ large enough since $V_n(t)$ and its derivatives tend to zero at infinity uniformly on compact time intervals:
\begin{align}
\label{condcenters}
&d_n\ge 4d_{n-1}+4\rho_{n-1},\qquad\rho_n\le\frac{d_n}{16},\qquad\frac{M_n}{d_n}+\frac{M_n\log(2+2d_n)^n}{d_n^2}\le\eps_n,\\
\label{condweighted}
&\sup_{0\le t\le k_n}\sup_{x\in\R^2}\frac{|V_n(t,x-d_ne_1)|}{\log(2+|x|)}\le\eps_n,\\
\label{condlocal}
&\sup_{0\le t\le k_n}\sup_{|x|\le k_n}|D^jV_n(t,x-d_ne_1)|\le\eps_n,\ \ 0\le j\le k_n.
\end{align}
For \fref{condweighted}, we split $|x|\le\sqrt{d_n}$, where $|x-d_ne_1|\ge d_n-\sqrt{d_n}$ is large, and $|x|>\sqrt{d_n}$, where $\log(2+|x|)\ge\log\sqrt{d_n}$ is large. Since $d_{n-1}+\rho_{n-1}\le\frac{d_n}{4}+\frac{d_n}{4}$ and $\rho_n\le\frac{d_n}{16}$, the disks $B_{\rho_n}(d_ne_1)$ are then pairwise disjoint. We denote by $u^{(N)}$, $\omega^{(N)}$ and $w_n$ the solutions and the interaction errors defined in \fref{defprefix}, \fref{definteraction}.\\

\noindent{\bf Step 2:} Local convergence. By \fref{definteraction},
\be
\label{telescoping}
u^{(N)}=\sum_{n=1}^NV_n(\cdot,\cdot-d_ne_1)+\sum_{n=1}^Nw_n.
\ee
Let $S>0$, $m\ge 0$ and $\rho>0$. For $n$ large, $\tau_n\ge S$ and $n+5\ge m+2$, so that $\sup_{t\le S}\|w_n(t)\|_{C^m}\lesssim\|w_n\|_{H^{n+5}}\le\eps_n$ by \fref{smallinteraction}; and $k_n\ge\max\{S,m,\rho\}$, so that the terms of the first sum are bounded by $\eps_n$ in $C([0,S];C^m(B_\rho))$ by \fref{condlocal}. Hence $u^{(N)}$ converges in $C([0,S];C^m(B_\rho))$ for all $S,m,\rho$. We denote by $u$ its limit and by $\omega=\curl u$. Passing to the limit in the equations, $\dv u=0$ and $\pa_t\omega+u\cdot\nabla\omega=0$ in the sense of distributions; the smoothness in time is proved in Step 5. At $t=0$, the sums defining $u_0$ and $\omega_0$ are locally finite.\\

\noindent{\bf Step 3:} Bounds on the vorticity and on the velocity. Since the initial vorticities have disjoint supports, $\|\omega^{(N)}(t)\|_{L^\infty}=\|\omega^{(N)}(0)\|_{L^\infty}\le 2$. Let $S>0$ and $j\ge 1$. By \fref{telescoping}, $D^j\omega^{(N)}$ is the sum of the terms $D^j\omega_n(\cdot,\cdot-d_ne_1)$, which are bounded by $2^{-n}$ on $[0,S]$ as soon as $k_n\ge\max\{S,j\}$ by \fref{hypassembly}, and of the terms $D^j\curl w_n$, which are bounded by $C\|w_n\|_{H^{n+5}}\le C\eps_n$ on $[0,S]$ as soon as $\tau_n\ge S$ and $n\ge j$. The finitely many remaining terms are bounded on $[0,S]$. Hence
$$\|\omega^{(N)}(t)\|_{L^\infty}\le 2,\qquad\sup_{0\le t\le S}\|D^j\omega^{(N)}(t)\|_{L^\infty}\le C_{S,j},$$
uniformly in $N$. In the same way, \fref{condweighted} and \fref{smallinteraction} give, uniformly in $N$,
\be
\label{logboundprefix}
|u^{(N)}(t,x)|\le C_S\log(2+|x|),\ \ 0\le t\le S.
\ee
These bounds pass to the limit, which proves (2) and \fref{loggrowthvelocity}.\\

\noindent{\bf Step 4:} The Biot--Savart law. By \fref{logboundprefix}, the flow $X^{(N)}$ of $u^{(N)}$ is defined for all times, and the same holds for $u$. Let $y\in\supp\psi_n(\cdot-d_ne_1)\subset B_{\rho_n}(d_ne_1)$, so that $\frac{15}{16}d_n\le|y|\le\frac{17}{16}d_n$ by \fref{condcenters}. As long as $X^{(N)}(t,y)$ stays in the annulus $\mathcal R_n=\{\frac{d_n}{2}\le|x|\le 2d_n\}$, \fref{logboundprefix} gives
$$|X^{(N)}(t,y)-y|\le C_SS\log(2+2d_n),\ \ 0\le t\le S,$$
which is smaller than $\frac{d_n}{4}$ for $n\ge n_S$ large. By a continuity argument, for $n\ge n_S$ and $N\ge n$, the vorticity carried by the $n$-th component, namely $\omega_n(0,\cdot-d_ne_1)\circ X^{(N)}(t)^{-1}$, is supported in $\mathcal R_n$ for $t\le S$, and its $L^1$ norm is $M_n$. The vorticity carried by the finitely many components $n<n_S$ starts in a fixed disk $B_{r_*}$, and since $\frac{d}{dt}|X^{(N)}|\le C_S\log(2+|X^{(N)}|)$ by \fref{logboundprefix}, it stays in a disk $B_{\rho'}$ for $t\le S$, where $\rho'$ is independent of $N$. Let now $\rho>0$, $x\in B_\rho$, and let $n_1\ge n_S$ be such that $d_n\ge 4\max\{\rho,\rho'\}$ for $n>n_1$. Then the vorticity $\omega^{(N)}(t)$ outside the disk $B_{2d_{n_1}}$ is carried by the components $n>n_1$, and on $\mathcal R_n$ we have $|K(x-y)|\le\frac{2}{\pi d_n}$. Hence, by \fref{condcenters},
$$\int_{|y|\ge 2d_{n_1}}|K(x-y)||\omega^{(N)}(t,y)|dy\le\sum_{n>n_1}\frac{2M_n}{\pi d_n}\le\sum_{n>n_1}\eps_n,$$
uniformly in $N$. Since $\omega^{(N)}(t)\to\omega(t)$ locally uniformly, Fatou's lemma gives the same bound for $\omega(t)$, and in particular $\int|K(x-y)||\omega(t,y)|dy<+\infty$. On the disk $B_{2d_{n_1}}$, the convergence is uniform and $K$ is locally integrable. Letting $N\to+\infty$ and then $n_1\to+\infty$, we obtain from $u^{(N)}=K*\omega^{(N)}$ that $u(t,x)=\int K(x-y)\omega(t,y)dy$. This proves (1).\\

\noindent{\bf Step 5:} Smoothness in space and time, and the pressure. We claim that for all $S>0$ and integers $q,m\ge 0$, there exist $C_{S,q,m}$ and $P_{q}$ such that, uniformly in $N$,
\be
\label{boundtimederivatives}
|\pa_t^qD^mu^{(N)}(t,x)|+|\pa_t^qD^m\omega^{(N)}(t,x)|\le C_{S,q,m}\log(2+|x|)^{P_q},\ \ 0\le t\le S.
\ee
For $q=0$, this follows from Step 3 and the local version (5) of Lemma \ref{lemmabiotsavart}, since $\|u^{(N)}(t)\|_{L^\infty(B_3(x))}\le C_S\log(2+|x|)$ by \fref{logboundprefix}. Assume \fref{boundtimederivatives} up to the order $q-1$ in time. Differentiating the vorticity equation $q-1$ times in time expresses $\pa_t^qD^m\omega^{(N)}$ in terms of derivatives of order at most $q-1$ in time, which proves the bound for $\omega^{(N)}$. For the velocity, we split $K=aK+(1-a)K$ with $a\in C^\infty_c(B_1)$ equal to one near zero. The local part $(aK)*\pa_t^q\omega^{(N)}$ is bounded by $\|aK\|_{L^1}\sup_{B_1(x)}|\pa^q_t\omega^{(N)}|$. For the nonlocal part, we write, by the conservation of the area,
$$((1-a)K)*\omega^{(N)}(t,x)=\sum_{n\le N}\int((1-a)K)\left(x-X^{(N)}(t,\alpha)\right)\omega_n(0,\alpha-d_ne_1)d\alpha.$$
Differentiating $q\ge 1$ times in time produces derivatives of $(1-a)K$ of order between $1$ and $q$, which are bounded by $C(1+|x-y|)^{-2}$, multiplied by time derivatives of the trajectories $\pa_t^iX^{(N)}$, $1\le i\le q$, which are material derivatives of $u^{(N)}$ of order at most $q-1$, and are bounded by a power of $\log(2+|X^{(N)}|)$ by the induction assumption. The vorticity at distance at most $4(1+|x|)$ from $x$ contributes at most a power of $\log(2+|x|)$, since $|\omega^{(N)}|\le 2$ and $\int_{|x-y|\le4(1+|x|)}(1+|x-y|)^{-2}dy\lesssim \log(2+|x|)$. The vorticity farther away is carried by components with $d_n\gtrsim|x|$, in the annuli $\mathcal R_n$, and contributes at most $C\sum_nM_n\log(2+2d_n)^{P_q}d_n^{-2}$, which is finite for every fixed $P_q$ by \fref{condcenters}. This proves the bound on $\pa_t^qu^{(N)}$, and the local Schauder estimates applied to $\pa^q_tu^{(N)}$, whose curl is $\pa_t^q\omega^{(N)}$, give the bounds on its spatial derivatives.

By \fref{boundtimederivatives} and the Arzel\`a--Ascoli theorem, the convergence of Step 2 holds in $C^\infty$ on compact subsets of $[0,+\infty)\times\R^2$, so that $(u,\omega)$ is smooth and solves the vorticity equation. In two dimensions, $\curl(\pa_tu+u\cdot\nabla u)=\pa_t\omega+u\cdot\nabla\omega=0$ for divergence-free fields, and since $\R^2$ is simply connected, there exists a smooth function $p$ such that $\pa_tu+u\cdot\nabla u=-\nabla p$, that is \fref{eulervelocity}.\\

\noindent{\bf Step 6:} Growth. Let $n\ge 1$ and $t\in[\tau_{n-1},\tau_n]$. Since $\nabla\omega_n(t)$ vanishes outside $\supp\omega_n(t)$ and at infinity, its maximal norm is attained at a point $x_t$ with $(t,x_t)\in\Sigma_n$, and $|\nabla\omega_n(t,x_t)|\ge 2e^{c_1e^{\kappa t}}$ by \fref{hypassembly}. By \fref{preservation} and $C_0\sum_j\eps_j\le 1$, $|\nabla\omega^{(N)}(t,x_t+d_ne_1)|\ge 2e^{c_1e^{\kappa t}}-1$ for all $N\ge n$, and the local convergence of Step 2 yields $\|\nabla\omega(t)\|_{L^\infty}\ge e^{c_1e^{\kappa t}}$. Since the intervals $[\tau_{n-1},\tau_n]$ cover $[\tau_0,+\infty)$, this proves (3).
\end{proof}

\subsection{Proof of Theorem \ref{thmdoubleexp}}

\noindent{\bf Step 1:} The components. Let $n_0$, $c_1$, $k_n$, $\sigma_n$, $h_n$ and $(\omega_n^{\rm per},V_n^{\rm per})$ be given by Proposition \ref{propsequence}. The period $L_n=\frac{2\pi}{h_n}$ is larger than $4$. We apply Proposition \ref{propcompactification} to $(\omega_n^{\rm per},V_n^{\rm per})$, both with $S=k_n$ and with $S=\tau_n$, and we choose $R_n\ge L_n$ so large that the solution $(\omega_n,V_n)=(\omega_{R_n},V_{R_n})$ satisfies the following. First, by \fref{globalderivatives} with $S=k_n$ and by \fref{earlysequence}, $\sup_{t\le k_n}\|D^j\omega_n(t)\|_{L^\infty}\le 2^{-n}$ for each of the finitely many $j\in[1,k_n]$. Second, by the local convergence (1) with $S=\tau_n$ on the compact set $[-L_n,L_n]^2$, which contains a period cell, and by \fref{growthwindow}, 
$$\|\nabla\omega_n(t)\|_{L^\infty}\ge\|\nabla\omega_n^{\rm per}(t)\|_{L^\infty}-1\ge 2e^{c_1e^{\kappa t}},\ \ t\in[\tau_{n-1},\tau_n].$$
Third, $\|\omega_n(0)\|_{L^\infty}\le 2$ by (3) of Proposition \ref{propcompactification}, and the stream function $\psi_{R_n}$ is supported in $B_{\rho_n}$ with $\rho_n=2R_n$. Hence the hypotheses of Proposition \ref{propassemblyinfinite} hold, and this proposition provides a solution $(\omega,u,p)$ with $\|\omega(t)\|_{L^\infty}\le 2$ and $\|\nabla\omega(t)\|_{L^\infty}\ge\exp(c_1e^{\kappa t})$ for $t\ge\tau_0$.\\

\noindent{\bf Step 2:} The initial data. Since $R_n\ge L_n$, the initial data of $(\omega_n,V_n)$ coincide with the periodic ones on the cell $\mathcal C_n=[-\frac{L_n}2,\frac{L_n}2]^2\subset B_{\rho_n}$, and since the disks $B_{\rho_n}(d_ne_1)$, which contain the supports of the initial data of the components, are pairwise disjoint, $(u_0,\omega_0)=(V_n(0),\omega_n(0))(\cdot-d_ne_1)$ on $d_ne_1+\mathcal C_n$. By \fref{initialsequence},
$$\int_{\R^2}|u_0|^2\ge\sum_n\frac{c_V^2}{h_n^4}=+\infty,\qquad\int_{\R^2}|\omega_0|^p\ge\sum_n\frac{2\pi^2}{h_n^2}=+\infty\ \ (p=1,2).$$
Moreover, since $|\mathcal C_n|=L_n^2$, $\|u_0\|_{L^\infty}\ge L_n^{-1}c_Vh_n^{-2}=\frac{c_V}{2\pi h_n}\to+\infty$, and $\pa_1u_{0,1}(d_ne_1)=\pa_1V^{\rm per}_{n,1}(0,0)\ge\kappa\log\frac{1}{\sigma_n}-C_2\to+\infty$.\\

\noindent{\bf Step 3:} Normalization. The functions
$$\widetilde\omega(t,x)=\frac12\omega\left(\frac t2,x\right),\qquad\widetilde u(t,x)=\frac12u\left(\frac t2,x\right),\qquad\widetilde p(t,x)=\frac14p\left(\frac t2,x\right)$$
form a smooth solution of \fref{euler} and \fref{eulervelocity}, by the first scaling in \fref{scalings} with $a=\frac12$. The Biot--Savart law, the logarithmic bound \fref{loggrowthvelocity}, the bounds on the derivatives on finite time intervals and the properties of the initial data are preserved, and $\|\widetilde\omega(t)\|_{L^\infty}\le 1$. Since $\frac{\kappa}{2}=\frac1\pi$, we have for $t\ge 2\tau_0$
$$\|\nabla\widetilde\omega(t)\|_{L^\infty}\ge\frac12\exp\left(c_1e^{\frac t\pi}\right)\ge\exp\left(\frac{c_1}{2}e^{\frac{t}{\pi}}\right)$$
as soon as $\frac{c_1}{2}e^{\frac t\pi}\ge\log 2$. This proves Theorem \ref{thmdoubleexp} with $c_0=\frac{c_1}{2}$. \hfill$\Box$

\section*{Acknowledgments}
Z. Zhang is partially supported by National Natural Science Foundation of China under Grant 12288101.

\section*{Disclosure of AI assistance}
Some of key ideas underlying the proofs in this paper emerged from discussions between the authors and an AI system(ChatGPT 6.0). 
The author also used the AI tools to help draft and revise the manuscript. The authors take full responsibility for the accuracy of the entire paper, including the correctness of all mathematical statements and proofs.

\end{document}